\documentclass[a4paper, twoside, 11pt]{article}

\usepackage[english]{babel}
\usepackage{amsmath,amssymb}
\usepackage{stmaryrd}
\usepackage{theorem}
\usepackage[all]{xy}
\usepackage{enumerate}
\usepackage{graphicx}
\usepackage{vmargin}
\usepackage{mathrsfs}

\theoremstyle{break}
\newtheorem{thm}{Theorem}[section]
\newtheorem{prop}{Proposition}[section]
\newtheorem{lem}{Lemma}[section]
\newtheorem{dfn}{Definition}[section]
\newtheorem{cor}{Corollary}[section]

\newcommand\C{\mathbb{C}}
\newcommand\F{\mathbb{F}}
\newcommand\Z{\mathbb{Z}}
\newcommand\N{\mathbb{N}}
\newcommand\Q{\mathbb{Q}}

\newcommand\id{\ensuremath{\mathrm{id}}}

\newcommand\Ind{\ensuremath{\mathrm{Ind}}}
\newcommand\Fil{\ensuremath{\mathrm{Fil}}}
\newcommand\Fonc{\ensuremath{\mathcal F}}
\newcommand\Frac{\ensuremath{\mathrm{Frac }}}
\newcommand\Gal{\ensuremath{\mathrm{Gal }}}
\newcommand\Hom{\ensuremath{\mathrm{Hom}}}
\newcommand\inv{\ensuremath{\mathrm{inv}}}
\newcommand\Res{\ensuremath{\mathrm{Res}}}
\newcommand\Ker{\ensuremath{\mathrm{Ker }}}
\newcommand\Tr{\ensuremath{\mathrm{Tr}}}
\newcommand\Frob{\ensuremath{\mathrm{Frob}}}

\newcommand\fcnl{\ensuremath{\mathscr L}}

\newcommand\Eproof{\hfill{{\large $\square$}}}
\newcommand\Epartproof{\hfill{{\large $\lozenge$}}}

\newcommand\maxim{\ensuremath{\mathfrak{m}}}
\newcommand\maximE{\ensuremath{\maxim_{\tilde{\bold E}}}}

\newcommand{\hooklongrightarrow}{\lhook\joinrel\longrightarrow}

\title{An explicit formula for the Hilbert symbol \\ of a formal group}

\author{Floric Tavares Ribeiro}
\date{}

\begin{document}

\maketitle

\begin{abstract}
In \cite{abr_hilb_forml}, Abrashkin established the
Br\"uckner-Vostokov formula for the Hilbert symbol of a formal group
under the assumption that roots of unity belong to the base field.
The main motivation of this work is to remove this hypothesis. It is
obtained by combining methods of ($\varphi, \Gamma$)-modules and a
cohomological interpretation of Abrashkin's technique. To do this,
we build ($\varphi, \Gamma$)-modules adapted to the false Tate curve
extension and generalize some related tools like the Herr complex
with explicit formulas for the cup-product and the Kummer
map.
\end{abstract}

\tableofcontents

\section*{Introduction}\addcontentsline{toc}{section}{Introduction}

\subsection{($\varphi, \Gamma$)-modules}
Let $p$ be a prime number and $K$ a finite extension of $\Q_p$ with
residue field $k$. Fix $\overline K$ an algebraic closure of $K$ and
note $G_K =\Gal(\overline K/K)$ the absolute Galois group of $K$.
Let us furthermore introduce $K_{\infty} = \cup_nK(\zeta_{p^n})$ the
cyclotomic extension of $K$ and $\Gamma_K= \Gal(K_{\infty}/K)$.

The context of this work is the theory of $p$-adic representations
of the Galois group of a local field, here $G_K$. We are
particularly interested in $\Z_p$-adic representations of $G_K$,
\emph{i.e.} $\Z_p$-modules of finite type endowed with a linear and
continuous action of $G_K$.

In \cite{rep_p_corloc}, Fontaine introduced the notion of a
($\varphi, \Gamma_K$)-module over the ring $\mathbf A_K$. This ring
is, when $K$ is absolutely unramified, the set of power series
$\sum_{n\in \Z}a_n X^n$ with $a_n \in \mathcal O_K$, $a_n$
$p$-adicly converging to $0$ as $n$ goes to $-\infty$ and $X$ a
variable on which $\varphi$ and $\Gamma_K$ act for $\gamma \in
\Gamma_K$ via
$$\varphi(X) = (1+X)^p-1 \\ \ ; \ \ \ \gamma(X) = (1+X)^{\chi(\gamma)}-1$$
where $\chi$ is the cyclotomic character.

A ($\varphi,\Gamma_K$)-module over $\mathbf A_K$ is then a module of
finite type over $\mathbf A_K$ endowed with commuting semi-linear
actions of $\varphi$ and $\Gamma_K$.

Fontaine  defined an equivalence of categories between the category
of $\Z_p$-adic representations of $G_K$ and the category of \'etale
($\varphi, \Gamma_K$)-modules over $\mathbf A_K$. Cherbonnier and
Colmez showed in \cite{CherbonnierColmezSurconvergentes} that any
$p$-adic representation is overconvergent, which established a first
link between the ($\varphi, \Gamma_K$)-module $D(V)$ of a
representation $V$ and its de Rham module which contains the
geometric information on $V$. Berger then, in
\cite{Berger_equadiff}, showed how to recover the de Rham module
$D_{dR}(V)$, the semi-stable module $D_{st}(V)$ or the crystalline
module $D_{crys}(V)$ of Fontaine's theory from $D(V)$. For
absolutely unramified crystalline representations, Wach furnished in
\cite{Wach_pot_cris} another powerful construction which permits to
recover $D_{crys}(V)$ in the ($\varphi, \Gamma_K$)-module $D(V)$.
This construction was studied in details and made more precise by
Berger (\cite{Berger_limitesrepcris}). ($\varphi, \Gamma_K$)-modules
are also intimately linked to Iwasawa theory as was shown in works
by Cherbonnier and Colmez (\cite{CherbonnierColmezIwa}), Benois
(\cite{benois_crys}) or Berger (\cite{Berger_BKexp}).

Let us eventually cite another significant result brought by Herr in
his PhD thesis (\cite{HerrCoho}) who furnished a three terms complex
in the ($\varphi, \Gamma_K$)-module of a representation, whose
homology computes the Galois cohomology of the representation.

\subsection{The false Tate curve extension}

The construction of ($\varphi, \Gamma_K$)-modules lies on the use of
the cyclotomic tower and shows its fundamental role in the study of
$p$-adic representations. But another extension appears as
particularly significant.

Fix $\pi$ a uniformizer of $K$ and $\pi_n$ a system of $p^n$-th
roots of $\pi$:
$$\pi_0 = \pi \ \textrm{ et } \ \forall n \in \N, \ \pi_{n+1}^p = \pi_n.$$
It is then the behavior in extension $K_\pi = \cup_nK(\pi_n)$ which
makes the difference between a crystalline and a semi-stable
representation.

Let us cite moreover the following remarkable result.
\begin{thm}[Breuil, Kisin]
The forgetful functor from the category of $p$-adic crystalline
representations of $G_K$ to the category of $p$-adic representations
of $G_{K_\pi}$ is fully faithful.
\end{thm}
This theorem was conjectured by Breuil in \cite{Breuil_CorpsNormes}
where it was shown under some conditions on the Hodge-Tate weights
of the representation, with the help of objects very similar to
Fontaine's ($\varphi, \Gamma_K$)-modules. Kisin proved this result
unconditionally in \cite{Kisin_crys}. Other results, in particular
by Abrashkin (\cite{abr_hilb_forml, abr_RamFil2}), encourage us to
introduce, like Breuil, ($\varphi, \Gamma$)-modules where the
cyclotomic extension $K_\infty$ is replaced by $K_\pi$. However
$K_\pi/K$ is not Galois and we only get $\varphi$-modules (also
studied by Fontaine in \cite{rep_p_corloc}).

Let us then consider the Galois closure $L$ of $K_\pi$ which is
nothing more than the compositum of $K_\pi$ and $K_\infty$, a
metabelian extension of $K$, \emph{the false Tate curve extension}.
What we lose here is the explicit description of the field of norms
of this extension. Note $G_\infty = \Gal(L/K)$. Our first result can
then, for $\mathbf A' = \mathbf A$ or $\mathbf{\tilde A}$, and
$\mathbf A_L' = \mathbf A'^{G_L}$ (where $\mathbf A$ are
$\mathbf{\tilde A}$ Fontaine rings defined in Paragraph
\ref{corpsE}), be expressed as:

\begin{thm}
The functor
$$\left. \begin{array}{ccc}
\left\{\Z_p-\textrm{adic representations of } G_K \right\} &
\rightarrow &
\left\{\textrm{\'etale }(\varphi,G_{\infty})-\textrm{modules over }\mathbf A_L'\right\}\\
V & \mapsto & D_L(V) = (V\otimes_{\Z_p}\mathbf A')^{G_L}
\end{array} \right.$$
is an equivalence of categories.
\end{thm}
In fact we show that the $(\varphi,G_{\infty})$-module $D_L(V)$ is
nothing but the scalar extension of the usual
$(\varphi,\Gamma_K)$-module $D(V)$ from $\mathbf A_K$ to $\mathbf
A_L'$.

\subsection{Galois cohomology}

We are now able to associate with a representation a ($\varphi,
G_\infty$)-module giving a better control of the behavior of the
representation in the extension $K_\pi$. But we would like to use
tools available in the classical framework, first of all Herr's
complex. Recall that in the usual case of
$(\varphi,\Gamma_K)$-modules, Herr showed in \cite{HerrCoho} that
the homology of the complex
$$\xymatrix{ 0 \ar[r] & D(V) \ar[r]^-{f_1}
& D(V) \oplus D(V) \ar[r]^-{f_2} & D(V) \ar[r] & 0}$$ with maps
$$f_1 =\left(\begin{array}{c} \varphi - 1  \\ \gamma - 1 \end{array}\right) \textrm{ et }
f_2 =(\gamma - 1, 1 - \varphi)$$ computes the Galois cohomology of
the representation $V$.

Since the group $G_\infty$ is now of dimension $2$, the
corresponding complex loses some simplicity. Let $\tau$ be a
topological generator of the sub-group $\Gal(L/K_\infty)$ and
$\gamma$ a topological generator of $\Gal(L/K_\pi)$ satisfying
$\gamma \tau \gamma^{-1} = \tau^{\chi(\gamma)},$ it can be described
as:

\begin{thm}
Let $V$ be a $\Z_p$-adic representation of $G_K$ and $D$ its
($\varphi, G_\infty$)-module. The homology of the complex
$$\xymatrix{ 0 \ar[r] & D \ar[r]^-{\alpha} & D\oplus D\oplus D \ar[r]^-{\beta} & D\oplus D\oplus D
\ar^-{\eta}[r] & D \ar[r] & 0\\}$$ where
$$ \alpha = \begin{pmatrix}
\varphi - 1  \\
\gamma - 1 \\
\tau - 1 \end{pmatrix},  \beta = \begin{pmatrix}
\gamma - 1 & 1 - \varphi & 0 \\
\tau - 1 & 0 & 1 - \varphi \\ 0 & \tau^{\chi(\gamma)} - 1 & \delta -
\gamma
\end{pmatrix},
\eta = \begin{pmatrix} \tau^{\chi(\gamma)} - 1, & \delta - \gamma, &
\varphi - 1 \end{pmatrix}$$ with $\delta = (\tau^{\chi(\gamma)} -
1)(\tau - 1)^{-1} \in \Z_p[[\tau-1]]$, identifies canonically and
functorially with the continuous Galois cohomology of $V$.
\end{thm}

In fact, we get explicit isomorphisms. In particular for the first
cohomology group, let $(x,y,z) \in \ker \beta$, let $b$ be a
solution in $V \otimes \mathbf A'$ of
$$(\varphi-1)b=x,$$ then the above theorem associates with the class of the triple $(x,y,z)$ the class
of the cocycle:
$$ c\ :\ \sigma \mapsto c_{\sigma} = -(\sigma - 1)b + \gamma^n \frac {\tau^m-1}{\tau - 1}z + \frac{\gamma^n - 1 }{\gamma - 1}y$$
where $\sigma_{|_{G_{\infty}}} = \gamma^n\tau^m$.

Moreover, like Herr in \cite{Herr_cup}, we furnish explicit formulas
describing the cup-product in terms of the four terms Herr complex
above.

\subsection{Explicit formulas for the Hilbert symbol}

The Hilbert symbol, for a field $K$ containing the group $\mu_{p^n}$
of $p^n$-th roots of unity is defined as the pairing
\begin{eqnarray*}
(,)_{p^n} : K^*/K^{*p^n} \times  K^*/ K^{*p^n} & \rightarrow & \mu_{p^n}\\
(a, b)_{p^n} & = & \left(\sqrt[p^n]{b}\right)^{r_{ K}(a)-1}
\end{eqnarray*}
where $r_{ K} :  K^* \rightarrow G_{ K}^{\textrm{ab}}$ is the
reciprocity map.

Since 1858 and Kummer's work, many explicit formulas have been given
for the Hilbert symbol. Let us cite the one of Coleman
(\cite{Coleman_dilog}): suppose that $K= K_0(\zeta_{p^n})$ where
$K_0$ is a finite unramified extension of $\Q_p$ and $\zeta_{p^n}$ a
fixed primitive $p^n$-th root of unity. Note $W$ the ring of
integers of $K_0$. If $F \in 1+ (p,X) \subset  W[[X]]$, then
$F(\zeta_{p^n}-1)$ is a principal unit in $K$ and all of them are
obtained in that way. Extend the absolute Frobenius $\varphi$ from $
W$ to $W[[X]]$ by putting $\varphi(X) = (1+X)^p-1$. Denote for $F\in
W[[X]]$
$$\fcnl(F) =  \frac 1 p \log \frac{F(X)^p}{\varphi(F(X))} \in W[[X]].$$
Then for $F\in 1+ (p,X)$,
$$\fcnl(F) = \left(1- \frac{\varphi} p\right) \log F(X).$$

Coleman's formula can then be written as:
\begin{thm}[Coleman]
Let $F,G \in 1 + (p,X) \subset  W[[X]]$, then
$$ (F(\zeta_{p^n}-1), G(\zeta_{p^n}-1))_{p^n} = \zeta_{p^n}^{[F,G]_n}$$
where
$$[F,G]_n = \Tr_{K_0/\Q_p}\circ \Res_X \frac 1 {\varphi^n(X)} \left(\fcnl(G)d\log F - \frac 1 p \fcnl(F)d\log G^\varphi\right).$$
\end{thm}

Let us furthermore cite the Br\"uckner-Vostokov formula: suppose now
that $p\neq 2$, let $\zeta_{p^n} \in K$, let $W$ be the ring of
integers of $K_0$, the maximal unramified extension of $K/\Q_p$.
Extend the Frobenius $\varphi$ from $W$ to $W[[Y]][1/Y]$ via
 $\varphi(Y) = Y^p$. Fix moreover $\pi$ a uniformizer of $K$.
\begin{thm}[Br\"uckner-Vostokov]
Let  $F, G \in (W[[Y]][1/Y])^\times$, then
$$ (F(\pi), G(\pi))_{p^n} = \zeta_{p^n}^{[F,G]_n}$$
where
$$[F,G]_n = \Tr_{K_0/\Q_p}\circ \Res_Y \frac 1 {s^{p^n}-1} \left(\fcnl(G)d\log F - \frac 1 p \fcnl(F)d\log G^\varphi\right)$$
with $s \in W[[Y]]$ such that $s(\pi) = \zeta_{p^n}$.
\end{thm}

The purpose of the second part of this work is to show a
generalization of this formula to the case of formal groups.

Remark that there are other types of formulas, in particular the one
of Sen (\cite{Sen_ExplReciLaws}), generalized to formal groups by
Benois in \cite{benois_periodes}.

We refer interested readers to Vostokov's \cite{Vostokov_formulas}
which provides a comprehensive background on explicit formulas for
the Hilbert symbol.

\subsection{An explicit formula for formal groups}

Let $G$ be a connected smooth formal group of dimension $d$ and of
finite height $h$ over the ring of Witt vectors $W=W(k)$ with
coefficients in a finite field $k$. Let $K_0$ be the fraction field
of $W$ and $K$ a finite extension of $K_0$ containing the
$p^M$-torsion $G[p^M]$ of $G$. Define then the Hilbert symbol of $G$
to be the pairing
\begin{eqnarray*}
(,)_{G,M} : K^* \times G(\maxim_K) & \rightarrow & G[p^M]\\
(x, \beta)_{G,M} & = & r_K(x)(\beta_1)-_G\beta_1
\end{eqnarray*}
where $r_K : K^* \rightarrow G_K^{\textrm{ab}}$ is the reciprocity
map and $\beta_1$ satisfies
$$p^M\id_G \beta_1 = \beta.$$

\bigskip

Fix a basis of logarithms of $G$ under the form of a vectorial
logarithm $l_G \in K_0[[\mathbf X]]^d$ where $\mathbf X =
(X_1,\dots,X_d)$ such that one has the formal identity
$$l_G(\mathbf X +_G \mathbf Y) = l_G(\mathbf X) + l_G(\mathbf Y).$$
Complete $l_G$ with almost-logarithms $m_G\in K_0[[\mathbf
X]]^{h-d}$ in a basis $\begin{pmatrix} l_G \\ m_G \end{pmatrix}$ of
the Dieudonn\'e module of $G$.

Fontaine defined in \cite{FontaineDivisible} (see also
\cite{ColmezPeriodesVarAb} for an explicit description) a pairing
between the Dieudonn\'e module and the Tate module of $G$
$$T(G) = \mathop{\lim}_{\longleftarrow}G[p^n].$$

Honda showed in \cite{HondaFG} the existence of a formal power
series of the form $\mathcal A = \sum_{n \geq 1} F_n\varphi^n$ with
$F_n \in M_d(W)$ such that
$$\left(1- \frac{\mathcal A} p\right)\circ l_G( \mathbf X) \in M_d(W[[\mathbf X]]).$$
Let us introduce moreover the approximated period matrix. Fix $(o^1,
\dots, o^h)$ a basis of $T(G)$ where $o^i = (o^i_n)_{n\geq 1}$ such
that $p\id_Go^i_n = o^i_{n-1}$. Approach $(o^1 = (o^1_n)_n, \dots,
o^h)$ by a basis $(o^1_M, \dots, o^h_M)$ of $G[p^M]$. Then for all
$i$, choose $\hat o^i_M \in F(YW[[Y]])$ such that $\hat o^i_M(\pi) =
o^i_M$. The matrix $\mathcal V_Y$ is then
$$\mathcal V_Y = \begin{pmatrix} p^Ml_G(\hat o^1_M) & \dots & p^Ml_G(\hat o^h_M)\\ p^Mm_G(\hat o^1_M) & \dots & p^Mm_G(\hat o^h_M) \end{pmatrix}.$$
It is an approximation of the period matrix $\mathcal V$.

Now we can state the reciprocity law which generalizes the
Br\"uckner-Vostokov law and which constitutes the goal of the second
part of this work:
\begin{thm}
Let $\alpha \in (W[[Y]][\frac 1 Y])^\times$ and $\beta \in
G(YW[[Y]])$. Coordinates of the Hilbert symbol
$(\alpha(\pi),\beta(\pi))_{G,M}$ in the basis $(o^1_M, \dots,
o^h_M)$ are
$$(\Tr_{W/\Z_p} \circ \Res_Y)\mathcal V_Y^{-1}\left(\begin{pmatrix} (1 - \frac {\mathcal A} p)\circ
l_G(\beta) \\ 0 \end{pmatrix}d_{\log}\alpha(Y)-\fcnl(\alpha) \frac d
{dY}
\begin{pmatrix} \frac {\mathcal A} p \circ l_G(\beta) \\ m_G(\beta)\end{pmatrix}\right).$$
\end{thm}
This formula was shown by Abrashkin in \cite{abr_hilb_forml} under
the assumption that $K$ contains $p^M$-th roots of unity. Vostokov
and Demchenko proved it in \cite{VostokovDemchenko} without any
condition on $K$ for formal groups of dimension $1$.

\subsection{The strategy}
The main idea of the proof is due to Benois who carried it out in
\cite{benois_crys} to show Coleman's reciprocity law. Let us recall
what it consists in.

The Hilbert symbol can be seen as a cup-product via the following
commutative diagram
$$\xymatrix{
K^*\times K^* \ar[r]^-{(,)_{p^n}}  \ar[d]_-{\kappa\times \kappa} & \mu_{p^n} \\
H^1(K,\mu_{p^n}) \times H^1(K,\mu_{p^n}) \ar[r]^-{\cup} &
H^2(K,\mu_{p^n}^{\otimes 2}) \ar[u]_-{\textrm{inv}_K} }$$ where
$\kappa$ is Kummer's map. He first explicitly computed Kummer's map
in terms of the Herr complex associated with the representation
$\Z_p(1)$, then he used Herr's cup-product explicit formulas and he
finally computed the image of the couple he obtained via the
isomorphism $\textrm{inv}_K$.

For a formal group, the situation is rather similar, we get the
diagram
$$\xymatrix{
K^*\times G(\maxim_K) \ar[r]^-{(,)_{G,M}}  \ar[d]_-{\kappa\times \kappa_G} & G[p^M]\\
H^1(K,\mu_{p^M}) \times H^1(K,G[p^M]) \ar[r]^-{\cup} &  H^2(K,\mu_{p^M}\otimes G[p^M]) \ar[u]_-{\textrm{inv}_K} \\
}$$ with
$$G[p^M]\simeq (\Z/p^M\Z)^h,$$
$$\textrm{and } \ H^2(K,\mu_{p^M}\otimes G[p^M])\simeq H^2(K,\Z/p^M\Z(1))\otimes_{\Z/p^M\Z} G[p^M].$$

Formulas for the Kummer map and the cup-product are shown in the
section on $(\varphi, \Gamma)$-modules. The computation of the
explicit formula for the map $\kappa_G : G(\maxim_K) \rightarrow
H^1(K,G[p^M])$ constitutes the technical axis of this work.

Abrashkin made use of the Witt symbol, and to conclude via the field
of norms of extension $K_\pi/K$, he used the compatibility of the
reciprocity map between the field of norms of an extension and the
basis field. Some of his intermediate results (\cite[Propositions
3.7 and 3.8]{abr_hilb_forml}) can be directly translated in the
language of ($\varphi, G_\infty$)-modules. Indeed, we want to
compute a triple $(x,y,z)$ in the first homology group of the Herr
generalized complex associated with the representation $G[p^M]$.
Abrashkin's results permit to obtain $x$, the vanishing of $y$ and
the belonging of $z$ to $W(\maximE)$ (where $\mathbf{\tilde E}$ is a
Fontaine ring, cf. \ref{corpsE} below). However we need to know $z$
modulo $XW(\maximE)$ and then we have to carry Abrashkin's
computations to the higher order.

\subsection{Organization of the paper}

This work splits in two parts. In the first one, we introduce
$(\varphi, G_\infty)$-modules and give the associated Herr complex
with explicit formulas between its homology and the cohomology of
the representation. Then we provide explicit formulas for the
cup-product and the Kummer map in terms of the Herr complex.\\
The second part is devoted to the proof of the Br\"uckner-Vostokov
formula for formal groups. The main difficulty lie in the fact that
the period matrix does not live in the right place: we introduce an
approximated period matrix and show that it enjoys similar
properties as the original matrix modulo suitable rings. Then, we
carry out the computation of the Hilbert symbol in terms of the Herr
complex.

\bigskip

\emph{Acknowledgements} This work is based on my PhD thesis under
the supervision of Denis Benois. I wish to thank him for the
precious ideas he shared with me and the time and energy he offered
me. I am also very grateful to Laurent Berger. He carefully read an
earlier version of this paper, some of his remarks allowed me to
improve it.

\section{$(\varphi, \Gamma)$-modules and cohomology}

\subsection{Notation}\label{notations}

Let $p$ be a prime.

Let us recall (cf. \cite{corloc}) that if $\mathbb K$ is a perfect
field of characteristic $p$, one can endow the space $\mathbb
K^{\N}$ of sequences of elements in $\mathbb K$ with a structure of
a local ring of characteristic $0$ absolutely unramified and with
residue field $\mathbb K$. It is called the ring of Witt vectors
over $\mathbb K$ and is denoted by $W(\mathbb K)$. Recall moreover
that this construction permits to define a multiplicative section of
the canonical surjection
$$ W(\mathbb K) \rightarrow \mathbb K, $$
called the Teichm\"uller representative and denoted by $[\ ]$. If
$R$ is a (unitary or not) subring of $\mathbb K$, we still denote by
$W(R)$ the Witt vectors with coefficients in $R$. It is then a
subring of $W(\mathbb K)$.

\bigskip

Fix $K$ a finite extension of $\Q_p$ with residue field $k$.\\
Denote $W=W(k)$ the ring of Witt vectors over $k$. Then $K_0 =
W\otimes_{\Z_p}\Q_p$ identifies with the maximal unramified
sub-extension of $\Q_p$ in $K$.\\
Fix $\overline K$ an algebraic closure of $K$ and denote
$$G_K =\Gal(\overline K/K)$$
the absolute Galois group of $K$ and $\C_p$ the $p$-adic completion
of $\overline K$. Endow $\C_p$ with the $p$-adic valuation $v_p$
normalized by
$$v_p(p)=1.$$
Recall that the action of $G_K$ on $\overline K$ extends by
continuity to $\C_p$.\\
Let us fix $\varepsilon = (\zeta_{p^n})_{n \geq 0}$ a coherent
system of $p^n$-th roots of unity, \emph{i.e.}
$\zeta_{p^n}^p=\zeta_{p^{n-1}}$ for all $n$, $\zeta_1 = 1$ and
$\zeta_p \neq 1$. Then
$$K_{\infty} := \bigcup_{n\in \N}K(\zeta_{p^n})$$
is the cyclotomic extension of $K$. Denote $G_{K_{\infty}} =
\Gal(\overline K/K_{\infty})$ its absolute Galois group
and $\Gamma_K = \Gal(K_{\infty}/K)$ the quotient.\\
Let us fix as well $\pi$ a uniformizer of $K$ and $\rho =
(\pi_{p^n})_{n \geq 0}$ a coherent system of $p^n$-th roots of
$\pi$. Denote
$$ K_\pi = \bigcup_{n \geq 0} K(\pi_{p^n}).$$
The extension $K_\pi/K$ is not Galois, so put
$$L =\bigcup_{n \geq 0} K(\zeta_{p^n}, \pi_{p^n})$$
its Galois closure. It is the compositum of $K_\pi$ and
$K_{\infty}$. Denote $G_L = \Gal(\overline K/L)$ its absolute Galois
group and $G_{\infty} = \Gal(L/K)$ the quotient. The cyclotomic
character $\chi : G_K \rightarrow \Z_p^*$ factorizes through
$G_{\infty}$ (even through $\Gamma_K$) ; it is also true for the map
$\psi : G_K \rightarrow \Z_p$ defined by
$$\forall g \in G_K \ \ \ \ g(\pi_{p^n})=\pi_{p^n}\zeta_{p^n}^{\psi(g)}.$$
Moreover, the group $G_{\infty}$ identifies with the semi-direct
product $\Z_p \rtimes \Gamma_K$. So $G_{\infty}$ is topologically
generated by two elements, $\gamma$ and $\tau$ satisfying:
$$\gamma \tau \gamma^{-1} = \tau^{\chi(\gamma)}.$$
Let us fix $\gamma$ and choose $\tau$ such that $\psi (\tau) = 1$,
\emph{i.e.} with
$$\tau (\rho) = \rho \varepsilon.$$

\bigskip

We adopt the convention that complexes have their first term in
degree $-1$ if this term is $0$, and otherwise in degree $0$.

\paragraph{Remark}
The group $G_{\infty}$ is a $p$-adic Lie group so that the extension
$L/K$ is arithmetically profinite (cf \cite{WinCorNorm,
Venjakob_WeierPrep}).

\subsection{The field $\mathbf{\tilde E}$, the ring $\mathbf{\tilde
A}$ and some of their subrings.}\label{corpsE} We refer to
\cite{rep_p_corloc} for results of this section. However we adopt
Colmez' notation. Rings $R$, $W(\Frac R)$ and $\mathcal
O_{\widehat{\mathcal E^{\textrm{nr}}}}$ of \cite{rep_p_corloc}
become $\mathbf{\tilde E}^+$, $\mathbf{\tilde A}$ and $\mathbf A$.

\bigskip

Define $\mathbf{\tilde E}$ as the inverse limit
$$\mathbf{\tilde E} = \mathop{\lim_{\leftarrow}}_n \ \C_p$$
where transition maps are exponentiation to the power $p$. An
element of $\mathbf{\tilde E}$ is then a sequence $x=(x^{(n)})_{n\in
\N}$ satisfying
$$(x^{(n+1)})^p=x^{(n)} \ \forall n\in \N.$$
Endow $\mathbf{\tilde E}$ with the addition
$$x+y = s \textrm{ where } s^{(n)}=\displaystyle{\lim_{m \to
+\infty}}(x^{(n+m)}+y^{(n+m)})^{p^m}$$ and the product
$$x.y = t \textrm{ where } t^{(n)}=x^{(n)}.y^{(n)}.$$
These operations make $\mathbf{\tilde E}$ into a field of
characteristic $p$, algebraically closed and complete for the
valuation
$$v_{\mathbf E}(x) :=v_p(x^{(0)}).$$
The ring of integers of $\mathbf{\tilde E}$, denoted by
$\mathbf{\tilde E}^+$, identifies then with the inverse limit
$\displaystyle{\lim_{\leftarrow}} \ \mathcal O_{\C_p}$. It is a
local ring whose maximal ideal, denoted by $\maximE$ identifies with
$\displaystyle{\lim_{\leftarrow}} \ \maxim_{\C_p}$ and whose residue
field is isomorphic to $\overline k$.\\
The field $\mathbf{\tilde E}$, as well as its ring of integers
$\mathbf{\tilde E}^+$, still has a natural action of $G_K$ which is
continuous with respect to the $v_{\mathbf E}$-adic topology. Define
the Frobenius
$$\varphi \ : \ x \mapsto x^p$$
which acts continuously, commutes with the action of $G_K$ and
stabilizes $\mathbf{\tilde E}^+$.

\bigskip

Let $\tilde {\mathbf A} = W(\tilde {\mathbf E})$  be the ring of
Witt vectors on $\mathbf{\tilde E}$ and $\tilde {\mathbf A}^+ =
W(\tilde {\mathbf E}^+)$.\\
Any element of $\tilde {\mathbf A}$ (respectively  $\tilde {\mathbf
A}^+$) can be written uniquely as
$$ \sum_{n\in \N}p^n[x_n]$$
where $(x_n)_{n\in \N}$ is a sequence of elements in $\mathbf{\tilde
E}$ (respectively in $\tilde {\mathbf E}^+$).\\
The topology on $\tilde {\mathbf A}$ comes from the product topology
on $W(\tilde {\mathbf E}) = \tilde {\mathbf E}^{\N}$. This topology
is compatible with the ring structure on $\tilde {\mathbf A}$. It is weaker than the $p$-adic topology.\\
Let us remark that the sequences $\varepsilon$ and $\rho$ introduced
below define elements in $\mathbf{\tilde E}^+$. Denote
$$X = [\varepsilon]-1 \textrm{ and }Y = [\rho].$$
These are elements of $\tilde {\mathbf A}^+$ and even of
$W(\maximE)$. They are topologically nilpotent. We also have a basis
of neighborhoods of $0$ in $\tilde {\mathbf A}$:
$$ \{p^n\tilde {\mathbf A} + X^m\tilde {\mathbf A}^+\}_{(n,m)\in \N^2} \textrm{ and } \{p^n\tilde {\mathbf A} + Y^m\tilde {\mathbf A}^+\}_{(n,m)\in \N^2}.$$
Let $W[[X,Y]]$ denote the subring of $\tilde {\mathbf A}^+$
consisting in sequences in $X$ and $Y$ ; it is stable under the
action of $G_K$ which is given by:
$$g (1+X) = (1+X)^{\chi(g)} \textrm{ and } g(Y) = Y(1+X)^{\psi(g)}$$
and the one of $\varphi$:
$$\varphi(X) = (1+X)^p-1 \textrm{ and } \varphi(Y) = Y^p.$$

\paragraph{Remark}$\ $\\
The specialization morphism for polynomials
$$\left. \begin{array}{ccc}
W[X_1,X_2] & \rightarrow & \mathbf{\tilde A}^+\\
X_1, \ X_2 & \mapsto & X, \ Y
\end{array}
\right.$$ is injective. However, the one for formal power series
$$\left. \begin{array}{ccc}
W[[X_1,X_2]] & \rightarrow & \mathbf{\tilde A}^+\\
X_1, \ X_2 & \mapsto & X, \ Y
\end{array}
\right.$$ is not a priori.

\bigskip

Let $\mathbf A_{\Q_p}$ denote the $p$-adic completion of
$\Z_p[[X]][\frac 1 X]$, it consists in the set
$$\mathbf A_{\Q_p} = \left\{\sum_{n \in \Z} a_nX^n | \ \forall n \in \Z, \ a_n \in \Z_p \textrm{ and }
a_n \mathop {\longrightarrow }\limits_{n \to -\infty } 0\right\}.$$
It is a local $p$-adic, complete subring of $\tilde {\mathbf A}$,
with residue field $\F_p((\varepsilon -1))$. Define $\mathbf A$ the
$p$-adic completion of the maximal unramified extension of $\mathbf
A_{\Q_p}$ in $\mathbf{\tilde A}$. Its residue field is then the
separable closure of $\F_p((\varepsilon -1))$ in $\tilde {\mathbf
E}$. Denote this field by $\mathbf E$. It is a dense subfield of
$\tilde {\mathbf E}$.

\subsection{Rings of $p$-adic periods.}

\subsubsection{$B_{dR}$ and some of its subrings}
We refer to \cite{cor_perio} for further details on these rings.\\
The map
$$\theta : \left\{ \begin{array}{ccl}
\tilde {\mathbf A}^+ & \rightarrow & \mathcal O_{\C_p}\\
\sum_{n\geq 0} p^n[r_n] & \mapsto & \sum_{n\geq 0} p^nr_n^{(0)}
\end{array} \right.$$
is surjective, with kernel $W^1(\tilde {\mathbf E}^+)$ which is a
principal ideal of $\tilde {\mathbf A}^+$ generated, for instance,
by $\omega = X/\varphi^{-1}(X)$. Denote
$$ B_{dR}^+ = \mathop{\lim_{\longleftarrow}}\limits_n (\tilde {\mathbf A}^+\otimes \Q_p)/(W^1(\tilde {\mathbf
E}^+)\otimes \Q_p)^n$$ the completion of $\tilde {\mathbf
A}^+\otimes \Q_p$ with respect to the $W^1(\tilde {\mathbf
E}^+)$-adic topology. The action of $G_K$ on $\tilde {\mathbf A}^+$
extends by continuity to $B_{dR}^+$. Yet it is not the case of the
Frobenius $\varphi$ which is not continuous with respect to the
$W^1(\tilde {\mathbf E}^+)$-adic topology. The sequence
$$\log [\varepsilon] = \sum_{n\geq1}(-1)^{n+1}\frac{X^n}{n}$$
converges in $B_{dR}^+$ towards an element denoted by $t$. Define
then
$$B_{dR} = B_{dR}^+[1/t].$$
It is the fraction field of $B_{dR}^+$. It is still endowed with an
action of $G_K$ for which
$$ B_{dR}^{G_K} = K$$ and with a
compatible, decreasing, exhaustive filtration
$$\Fil^kB_{dR} = t^kB_{dR}^+.$$

\bigskip

Define now the ring $A_{crys}$ to be the $p$-adic completion of the
divided powers envelop of $\tilde {\mathbf A}^+$ with respect to
$W^1(\tilde {\mathbf E}^+)$. It consists in the sequences
$$\left\{ \sum_{n\geq 0} a_n \frac{\omega^n}{n!}  \textrm{ such that }
 a_n \in \tilde {\mathbf A}^+ \textrm{ and } a_n \to 0 \ p\textrm{-adically.} \right\}.$$
This ring is naturally a subring of $B_{dR}$. Moreover, the sequence
defining $t$ still converges in $A_{crys}$ and we set
$$B_{crys}^+ = A_{crys}\otimes \Q_p \textrm{ and } B_{crys}=B_{crys}^+[1/t] = A_{crys}[1/t].$$
Moreover, if one chooses $\tilde p =(p_0,p_1, \dots) \in
\mathbf{\tilde E}$ with $p_0 = p$, then the series $\log \frac
{[\tilde p]} p$ converges in $B_{dR}$ towards a limit denoted by
$\log [\tilde p]$ (with the implicit convention $\log p =0$). Define
then
$$B_{st} = B_{crys}[\log [\tilde p]].$$
It is still a subring of $B_{dR}$.\\
All these rings, endowed with their $p$-adic topology, come with a
continuous action of $G_K$, the filtration induced by the one on
$B_{dR}$, and a Frobenius $\varphi$ extending by continuity the one
on $\tilde {\mathbf A}^+$. Note that
$$B_{crys}^{G_{K}} = K_0 \ \ \textrm{ and } \ \ B_{st}^{G_{K}} = K_0.$$

\subsubsection{A classification of $G_K$-representations}
We call a \emph{$\Z_p$-adic representation of $G_K$} any finitely
generated $\Z_p$-module with a linear, continuous action of $G_K$
and a \emph{$p$-adic representation of $G_K$} any finite dimensional
$\Q_p$-vector space with a linear, continuous action of $G_K$. A
$\Z_p$-adic representation is then turned into a $p$-adic
representation by tensorizing by $\Q_p$.\\
Let $V$ be a $p$-adic representation of $G_K$. Note
\begin{eqnarray*}
D_{dR}(V) & := & (V \otimes_{\Q_p} B_{dR})^{G_K}\\
D_{st}(V) & := & (V \otimes_{\Q_p} B_{st})^{G_{K}}\\
D_{crys}(V) & := & (V \otimes_{\Q_p} B_{crys})^{G_{K}}.
\end{eqnarray*}
$D_{dR}(V)$ (respectively $D_{st}(V)$, $D_{crys}(V)$) is a $K$
(respectively $K_0$, $K_0$)-vector space of dimension lower or equal
to the dimension of $V$ on $\Q_p$. The representation $V$ is said to
be \emph{de Rham} (respectively \emph{semi-stable},
\emph{crystalline})
when these dimensions are equal.\\
One immediately sees that crystalline representations are
semi-stable and semi-stable representations are de Rham.

We say as well that a $\Z_p$-adic representation $V$, free over
$\Z_p$, is de Rham, semi-stable or crystalline when so is the
$p$-adic representation $V\otimes_{\Z_p}\Q_p$.

\paragraph{Example} \emph{The false Tate curve}\\
Let us define the false Tate curve (or Tate's representation) by
$$ V_{Tate} = \Z_pe_1 + \Z_p e_2$$
with the action of $G_K$:
$$\left\{ \begin{array}{l}
g(e_1) = \chi(g)e_1\\
g(e_2) = \psi(g) e_1 + e_2
\end{array}\right.$$
for all $g\in G_K$, where $\chi$ is the cyclotomic character and
$\psi$ is defined in Paragraph \ref{notations}. This representation
is an archetypic semi-stable representation and will be an important
reference. We will confront our approach with it, in particular
modified $(\varphi, \Gamma)$-modules. For the moment, just note that
the action of $G_K$ on $V_{Tate}$ factorizes through $G_{\infty}$.

The name "false Tate curve" comes from the similarity of this module
with the Tate module of an elliptic curve with split multiplicative
reduction at $p$.

\subsection{Fontaine's theory}

Let $R$ be a topological ring with a linear, continuous action of
some group $\Gamma$ and a continuous Frobenius $\varphi$ commuting
with the action of $\Gamma$. Call a \emph{$(\varphi, \Gamma)$-module
on $R$} any finitely generated $R$-module $M$ with commuting
semi-linear actions of  $\Gamma$ and $\varphi$. A $(\varphi,
\Gamma)$-module on $R$ is moreover said \emph{\'etale} if the image
of $\varphi$ generates $M$ as an $R$-module:
$$R\varphi(M) = M.$$

\subsubsection{The classical case}

Let us recall the theory of $(\varphi, \Gamma)$-modules introduced
by Fontaine in \cite{rep_p_corloc}.

\bigskip

Set $\mathbf A_K = \mathbf A^{G_{K_{\infty}}}.$

Define the functors
$$D : V \mapsto D(V) = (\mathbf A \otimes_{\Z_p} V)^{G_{K_{\infty}}}$$ from the category of $\Z_p$-adic representations of $G_K$ to the one of
$(\varphi, \Gamma_K)$-modules on $\mathbf A_K$ and
$$ V : M \mapsto V(M) = (\mathbf A \otimes_{\mathbf A_K} M)^{\varphi = 1}$$
from the category of \'etale $(\varphi, \Gamma_K)$-modules on
$\mathbf A_K$ to the one of $\Z_p$-adic representations of $G_K$.
The following theorem was shown by Fontaine (\cite{rep_p_corloc}):

\begin{thm}\label{thm_Fontaine}
Natural maps
$$ \mathbf A \otimes_{\mathbf A_K} D(V) \rightarrow \mathbf A \otimes_{\Z_p} V $$
$$ \mathbf A \otimes_{\Z_p} V(M) \rightarrow \mathbf A \otimes_{\mathbf A_K} M$$
are isomorphisms. In particular, $D$ and $V$ are quasi-inverse
equivalences of categories between the category of $\Z_p$-adic
representations of $G_K$ and the one of \'etale $(\varphi,
\Gamma_K)$-modules on $\mathbf A_K$.
\end{thm}

\paragraph{Example} The $(\varphi, \Gamma_K)$-module of the false
Tate curve admits a basis of the form $(1\otimes e_1, b\otimes e_1 +
1 \otimes e_2)$ where $b\in \mathbf A_L$ satisfies $(\tau - 1) b =
-1$. However $V_{Tate}$ is not potentially crystalline, and then, by
a theorem of Wach (cf. \cite{Wach_pot_cris}), not \emph{of finite
height}, which means $b \notin \mathbf A_L^+ = \mathbf A_L \bigcap
\mathbf{\tilde A}^+$.

\bigskip

We want to build a $(\varphi, \Gamma)$-module which furnishes more
information (which will then be redundant but easier to use) on the
behavior of the associated representation in the extension $K_\pi/K$
or in its Galois closure $L/K$. For this, we want $\Gamma =
G_{\infty}$.

\subsubsection{The metabelian case}\label{hypoA}
Suppose $\mathbf A' =  \mathbf A$ or $\mathbf A' =\mathbf {\tilde
A}$. Then, $\mathbf A'$ is a complete $p$-adic valuation ring,
stable under both $G_K$ and $\varphi$. Its residue field $\mathbf E'
= \mathbf E$ or $\mathbf{\tilde E}$ is separably closed.

\bigskip

Set $\mathbf A_L' = \mathbf A'^{G_L}$ ; if $\mathbf E'_L = \mathbf
E'^{G_L}$.
Then $\mathbf A_L'$ is a complete $p$-adic valuation ring with residue field $\mathbf E'_L$.\\
For any $\Z_p$-adic representation $V$ of $G_K$, define
$$D_L'(V) = (\mathbf A' \otimes_{\Z_p} V)^{G_L}$$
and for any $(\varphi, G_{\infty})$-module $D$, \'etale over
$\mathbf A_L'$,
$$V_L'(D) = (\mathbf A' \otimes_{\mathbf A_L'} D)^{\varphi =1}.$$
Denote these functors by $D_L$ and $V_L$ when $\mathbf A' = \mathbf
A$ and by $\widetilde D_L$ and $\widetilde V_L$ when $\mathbf A' =
\mathbf {\tilde A}$.

Remark that $D_L'(V)$ and $D(V)\otimes_{\mathbf A_K} \mathbf A_L'$
are $(\varphi, G_{\infty})$-modules over $\mathbf A_L'$, the latter
being \'etale. The following theorem shows that they are indeed
isomorphic and assures that $D_L'$ is a good equivalent for $D$ in
the metabelian case.

\begin{thm}\label{thm_isom}
\begin{enumerate}
\item
The natural map
$$\iota \ : \ D(V)\otimes_{\mathbf A_K} \mathbf A_L' \rightarrow D_L'(V)$$
is an isomorphism of $(\varphi,G_{\infty})$-modules \'etale over
$\mathbf A_L'$.
\item
Functors $D_L'$ and $V_L'$ are quasi-inverse equivalences of
categories between the category of $\Z_p$-adic representations of
$G_K$ and the one of \'etale $(\varphi,G_{\infty})$-modules on
$\mathbf A_L'$.
\end{enumerate}
\end{thm}
\emph{Proof:} First, remark that, because of Theorem
\ref{thm_Fontaine}, and after extending scalars, the natural map
$$D(V) \otimes_{\mathbf A_K} \mathbf A' \rightarrow V\otimes_{\Z_p} \mathbf A'$$
is an isomorphism.\\
Taking Galois invariants, we get an isomorphism
$$D(V) \otimes_{\mathbf A_K} \mathbf A_L'  = (D(V) \otimes_{\mathbf A_K} \mathbf A')^{G_L} \mathop{\longrightarrow}^{\sim}
(V\otimes_{\Z_p} \mathbf A')^{G_L} = D_L'(V)$$ as desired.

\bigskip

We immediately deduce that the functor $D_L'$ from the category of
$\Z_p$-adic representations of $G_K$ to the one of \'etale
$(\varphi,G_{\infty})$-modules over $\mathbf A_L'$ is exact and
faithful.

In fact, this result and the expression of the quasi-inverse of
$D_L'$ (seen as an equivalence of categories on its essential image)
suffice for our use of $(\varphi,G_{\infty})$-modules. This
quasi-inverse is obtained with the help of the comparison
isomorphism after extending scalars:
$$ D_L'(V) \otimes_{\mathbf A_L'} \mathbf A' \simeq D(V) \otimes_{\mathbf A_K} \mathbf A' \simeq V \otimes_{\Z_p} \mathbf A'$$
so that
$$ V_L'(D_L'(V)) \simeq V $$
and $V_L'$ is the quasi-inverse of $D_L'$.

Fontaine's computation (cf. \cite[Proposition 1.2.6.]{rep_p_corloc})
still applies here and permits to compute the essential image which
is still the category of \'etale $(\varphi,G_{\infty})$-modules over
$\mathbf A_L'$. It consists in proving that any $p$-torsion \'etale
$(\varphi,G_{\infty})$-module, which is then an $\mathbf E'$-vector
space, has a $\varphi$-invariant basis, by showing that for any
matrix $(a_{j,l}) \in GL_d(\mathbf E')$, the system
$$x_j^p = \sum a_{j,l}x_l$$
admits $p^d$ solutions in $\mathbf E'^d$, generating $\mathbf E'^d$.
The general case is deduced by \emph{d\'evissage} and passing to the
limit. \Eproof

\begin{cor}
The functor
$$\left. \begin{array}{ccc}
\left\{\textrm{\'etale}(\varphi,\Gamma_K)-\textrm{modules over
}\mathbf A_K\right\} & \rightarrow &
\left\{\textrm{\'etale}(\varphi,G_{\infty})-\textrm{modules over }\mathbf A_L'\right\}\\
D & \mapsto & D\otimes_{\mathbf A_K} \mathbf A_L'
\end{array} \right.$$
is an equivalence of categories.
\end{cor}

\paragraph{Example}
The ($\varphi,G_{\infty}$)-module associated with the false Tate
curve admits a trivial basis $(1\otimes e_1, 1 \otimes e_2)$. This
module is then of finite height over $\mathbf A_L'$. It would be
interesting to know whether this remains true or not for any
semi-stable representation.\\
When $\mathbf A' = \mathbf{\tilde A}$, it follows from a result of
Kisin (\cite[Lemma 2.1.10]{Kisin_crys}). He builds the
$\varphi$-module associated with the extension $K_{\pi}$:
$$ (V\otimes_{\Z_p}\mathbf A_Y)^{\Gal(\overline K/ K_{\pi})}$$
where $\mathbf A_Y$ is the $p$-adic completion of the maximal
unramified extension of $W[[Y]][\frac 1 Y]$ in $\mathbf{\tilde A}$.
He shows that semi-stable representations are of finite height in
this framework, which means that the $W[[Y]]$-module
$$ (V\otimes_{\Z_p}W[[Y]]^{nr})^{\Gal(\overline K/ K_{\pi})},$$ with
$W[[Y]]^{nr} = \mathbf A_Y \bigcap \mathbf{\tilde A}^+$, has the
same rank as $V$.

\subsubsection{Remark: the field of norms of $L/K$}\label{remark}

As previously remarked, the extension $L/K$ is arithmetically
profinite ; consider then its field of norms $\mathbf E_{L/K}$ which
can be explicitly described. Indeed if $k_L$ is the residue field of
$L$, then there exists $z \in \mathbf{\tilde E}$ such that $\mathbf
E_{L/K}$ identifies with $k_L((z)) \subset \mathbf{\tilde E}$. We
would then like to reproduce the classical construction of
($\varphi$, $\Gamma$)-modules by substituting $\mathbf E_{L/K}$ to
the field of norms of the cyclotomic extension $K_{\infty}/K$.
However, we then have to build a characteristic $0$ lift (in
$\mathbf{\tilde A}$) of $\mathbf E_{L/K}$ stable under both actions
of $G_K$ and $\varphi$, that we are not able to do. This problem is
linked to the fact that we cannot make explicit a norm coherent
sequence of uniformizers in the tower $K(\zeta_{p^n}, \pi_{p^n})$.

\subsection{Galois Cohomology}\label{section_coho_gal}
\subsubsection{Statement of the theorem}
First recall the classical case. Let $D(V)$ be the \'etale
($\varphi$, $\Gamma_K$)-module over $\mathbf A_K$ associated with a
$\Z_p$-adic representation $V$. Fix $\gamma$ a topological generator
of $\Gamma_K$. Herr introduced in \cite{HerrCoho} the complex
$$\xymatrix{ 0 \ar[r] & D(V) \ar[r]^-{f_1}
& D(V) \oplus D(V) \ar[r]^-{f_2} & D(V) \ar[r] & 0}$$ with maps
$$f_1 =\left(\begin{array}{c} \varphi - 1  \\ \gamma - 1 \end{array}\right) \textrm{ and }
f_2 =(\gamma - 1, 1 - \varphi).$$ He showed that the homology of
this complex canonically and functorially identifies with the Galois cohomology of the representation $V$.\\
This identification was explicitly given in
\cite{CherbonnierColmezIwa} and \cite{benois_crys} for the first
cohomology group by associating with the class of a pair $(x,y)$ of
elements in $D(V)$ satisfying $(\gamma-1)x = (\varphi -1)y$ the
class of the cocycle
$$\sigma \mapsto -(\sigma -1)b + \frac{\gamma^n -1}{\gamma -1}y$$
where $b\in V\otimes_{\Z_p}\mathbf A$ is a solution of $(\varphi
-1)b=x$ and $\sigma_{|\Gamma_K} =\gamma^n$ for some $n\in \Z_p$.\\
We will show that there still exists such a complex in the
metabelian case. However, in order to take into account that
$G_{\infty}$ has now two generators, we will modify it a little.\\
Let $M$ be a given \'etale ($\varphi,G_{\infty}$)-module  over
$\mathbf A_L'$. Associate with $M$ the four terms complex
$C_{\varphi, \gamma, \tau}(M)$:
$$\xymatrix{ 0 \ar[r] & M \ar[r]^-{\alpha} & M\oplus M\oplus M \ar[r]^-{\beta} & M\oplus M\oplus M
\ar^-{\eta}[r] & M \ar[r] & 0\\}$$ where
$$ \alpha = \begin{pmatrix}
\varphi - 1  \\
\gamma - 1 \\
\tau - 1 \end{pmatrix},  \beta = \begin{pmatrix}
\gamma - 1 & 1 - \varphi & 0 \\
\tau - 1 & 0 & 1 - \varphi \\ 0 & \tau^{\chi(\gamma)} - 1 & \delta -
\gamma
\end{pmatrix},
\eta = \begin{pmatrix} \tau^{\chi(\gamma)} - 1, & \delta - \gamma, &
\varphi - 1 \end{pmatrix}$$ with $\delta = (\tau^{\chi(\gamma)} -
1)(\tau - 1)^{-1} \in \Z_p[[\tau-1]]$ defined as follows: set
$$\binom{u}{n} = \frac {u.(u-1)...(u-n+1)}{n!} \in \Z_p \textrm{ for all } u \in \Z_p \textrm{ and all } n \in \N.$$
Then:
$$\tau^{\chi(\gamma)} = \sum_{n \geq 0} \binom{\chi(\gamma)}{n} (\tau -1)^n$$
for $\tau^{p^n}$ converges to $1$ in $G_\infty$, and thus $\tau-1$
is topologically nilpotent in $\Z_p[[G_{\infty}]]$. So
$$\delta = \frac{\tau^{\chi(\gamma)} - 1}{\tau - 1} = \sum_{n \geq 1} \binom{\chi(\gamma)}{n} (\tau -1)^{n-1}.$$
The purpose of this paragraph is to show
\begin{thm}\label{coho_gal}
Let $V$ be a $\Z_p$-adic representation of $G_K$.
\begin{enumerate}[i)]
\item  The homology of the complex $C_{\varphi, \gamma, \tau}(D_L(V))$
canonically and functorially identifies with the continuous Galois
cohomology of $V$.
\item Explicitly, let $(x,y,z) \in Z^1(C_{\varphi, \gamma,
\tau}(D_L(V)))$, let $b$ be a solution in $V \otimes \mathbf A'$ of
$$(\varphi-1)b=x,$$ then the identification above associates with the class of the triple $(x,y,z)$ the class
of the cocycle:
$$ c\ :\ \sigma \mapsto c_{\sigma} = -(\sigma - 1)b + \gamma^n \frac {\tau^m-1}{\tau - 1}z + \frac{\gamma^n - 1 }{\gamma - 1}y$$
where $\sigma_{|_{G_{\infty}}} = \gamma^n\tau^m$.
\end{enumerate}
\end{thm}

\subsubsection{Proof of Theorem \ref{coho_gal} $i)$} The functor $F^\bullet$ which
associates with a $\Z_p$-adic representation $V$ the homology of the
complex $C_{\varphi, \gamma, \tau}(D_L(V))$ is a cohomological
functor coinciding in degree $0$ with the continuous Galois
cohomology of $V$:
$$H^0(C_{\varphi, \gamma, \tau}(D_L(V)) = D_L(V)_{\varphi=1, \gamma = 1, \tau = 1} = V^{G_K}.$$
The proof consists then in showing that it is effaceable. In order
to do that, we would like to work with a category with sufficiently
many injectives and to see $V$ as a submodule of an explicit
injective, its induced module, which is known to be cohomologically
trivial. But the category of $\Z_p$-adic representations of $G_K$
doesn't admit induced modules. We will then work modulo $p^r$ for a
fixed $r$, and even in the category of direct limits of
$p^r$-torsion representations and then deduce the result by passing
to the limit. We have then to show that the homology of the complex
associated with an induced module concentrates in degree $0$, which
shows \emph{a fortiori} the effaceability of $F^\bullet$. We will
yet write this in an explicit manner, which will let us get the
second part of the theorem, and, in the next paragraph, an explicit
description of the cup-product in terms of the Herr complex.

\bigskip

Let $M_{G_K,p^r-tor}$ be the category of discrete $p^r$-torsion
$G_K$-modules, it is also the category of direct limits of finite
$p^r$-torsion $G_K$-modules or also the one of discrete
$\Z/p^r\Z[[G_K]]$-modules. Let us remark that the functor $D_L$
extends to an equivalence of categories from this category to the
one of direct limits of $p^r$-torsion \'etale
$(\varphi,G_{\infty})$-modules over $\mathbf A_L'$.\\
Note finally that this category is stable under passing to the
induced module:
\begin{lem}\label{induit}
Let $V$ be an object of $M_{G_K,p^r-tor}$, define the induced module
associated with $V$ by:
$$\Ind_{G_K}(V):=\Fonc_{cont}(G_K,V)$$
the set of all continuous maps from $G_K$ to $V$.\\
Endow $\Ind_{G_K}(V)$ with the discrete topology and the action of
$G_K$:
$$\left. \begin{array}{ccl}
G_K \times \Ind_{G_K}(V) & \rightarrow & \Ind_{G_K}(V)\\
 g.\eta & = & [x \mapsto \eta(x.g)].
\end{array} \right.$$
Then $\Ind_{G_K}(V)$ is an object of $M_{G_K,p^r-tor}$ and $V$
canonically injects in $\Ind_{G_K}(V)$.
\end{lem}
\emph{Proof:} The first part of the lemma is well-known. See
\cite{mathese} for details. The injection of $V$ in its induced
module is given by sending $v \in V$ on $\eta_v \in \Ind_{G_K}(V)$
such that
$$\forall g \in G_K \ \ \eta_v(g)=g(v).$$
\Eproof\\
Let $F^i$ denote the composed functor $H^i(C_{\varphi, \gamma,
\tau}(D_L(-)))$. The snake lemma gives for any short exact sequence
in $M_{G_K,p^r-tor}$
$$0 \rightarrow V \rightarrow V'' \rightarrow V' \rightarrow 0$$
a long exact sequence
$$0 \rightarrow F^0(V) \rightarrow F^0(V'') \rightarrow F^0(V') \rightarrow F^1(V) \rightarrow F^1(V'') \rightarrow \cdots$$
which shows that $F^{\bullet}$ is a cohomological functor.\\
Let us show that it coincides with the long exact cohomology
sequence when $V''=\Ind_{G_K}(V)$. We use the following result:
\begin{prop} \label{effacable}
Let $U=\Ind_{G_K}(V)$ be an induced module in the category
$M_{G_K,p^r-tor}$, then
$$F^i(U)=H^i(K,U)=0 \textrm{ for all } i>0.$$
\end{prop}

Let us first deduce Point $i)$ of the theorem from this result. The
commutative diagram
$$\xymatrix{
 0 \ar[r] &  F^0(V) \ar[r] \ar@{=}[d] & F^0(\Ind_{G_K}(V)) \ar[r] \ar@{=}[d] & F^0(V') \ar[r] \ar@{=}[d] & F^1(V) \ar[r] & 0\\
 0 \ar[r] &  H^0(K,V) \ar[r] &  H^0(K,\Ind_{G_K}(V)) \ar[r] &  H^0(K,V') \ar[r] &  H^1(K,V)  \ar[r] & 0 \\}
$$
shows that $H^1(K,V) \simeq F^1(V)$.

And in higher dimension vanishing of $F^i(\Ind_{G_K}(V))$ and
$H^i(K,\Ind_{G_K}(V))$ prove both that $F^k(V')=F^{k+1}(V)$ and
$H^k(K,V')=H^{k+1}(K,V)$. Thus, by induction, $F^i(V)=H^i(K,V)$
holds for all $i \in \N$ and for any module $V$ in
$M_{G_K,p^r-tor}$.

\bigskip

{\it Proof of the proposition:}

\bigskip

The Galois cohomology part is a  classical result (cf. \cite{corloc,
cohogal} or \cite{mathese}).

For the second part, we will use
\begin{lem}
For any $V \in M_{G_K,p^r-tor}$, there is a short exact sequence:
$$\xymatrix{ 0 \ar[r] &
\Ind_{G_{\infty}}(V) \ar[r] & D_L(\Ind_{G_K}(V)) \ar[r]^-{\varphi-1}
& D_L(\Ind_{G_K}(V)) \ar[r] & 0.}$$ Moreover, for any $\alpha \in
\Z_p^*$, there is a short exact sequence:
$$\xymatrix{ 0 \ar[r] & \Ind_{\Gamma_K}(V) \ar[r] & \Ind_{G_{\infty}}(V) \ar[r]^-{\tau^{\alpha}-1} &
\Ind_{G_{\infty}}(V) \ar[r] & 0.}$$ Finally, there is a short exact
sequence:
$$\xymatrix{ 0 \ar[r] & V^{G_K} \ar[r] & Ind_{\Gamma_K}(V) \ar[r]^-{\gamma-1} &
Ind_{\Gamma_K}(V) \ar[r] & 0.}$$
\end{lem}
\emph{Proof of the lemma:} Consider the short exact sequence
$$\xymatrix{ 0 \ar[r] & \Z_p \ar[r] & \mathbf A' \ar[r]^-{\varphi-1} & \mathbf A' \ar[r] & 0\\}$$
and tensorize it with $\Ind_{G_K}(V)$. The existence of a continuous
section of $\varphi-1$ (cf. \cite{Scholl_higher}) permits, taking 
Galois invariants, to get a long exact sequence beginning with
$$\xymatrix@C=19pt{ 0 \ar[r] & \Ind_{G_K}(V)^{G_L} \ar[r] &
D_L(\Ind_{G_K}(V)) \ar[r]^-{\varphi-1} & D_L(\Ind_{G_K}(V)) \ar[r] &
H^1(L,\Ind_{G_K}(V)) \\}$$
The kernel is given by $\Ind_{G_K}(V)^{G_L} = \Ind_{G_{\infty}}(V)$.\\
It remains to show the nullity of $H^1(G_L,\Ind_{G_K}(V))$. Remark
(cf.\cite[Chapitre I, Proposition 8]{cohogal}):
$$ H^1(G_L,\Ind_{G_K}(V)) = \mathop{\lim_{\rightarrow}}  H^1(G_M,\Ind_{G_K}(V)) $$
where the direct limit is taken  over the set of all finite Galois
sub-extensions $M$ of $L/K$. Indeed, the sub-Galois groups $G_M$ of
$G_K$ form, for inclusion, a projective system with limit
$$\mathop{\lim_{\leftarrow}} G_M = \bigcap G_M = G_L$$
and this system is compatible with the inductive system formed by
the $G_M$-modules by restriction
$\Ind_{G_K}(V)$ whose limit is the $G_L$-module by restriction $\Ind_{G_K}(V)$.\\
To prove the lemma, it suffices then to show for any finite Galois
extension $M/K$ included in $L$ the vanishing of $H^1(G_M,\Ind_{G_K}(V))$.\\
But, $G_M$ being open in $G_K$, we have the finite decomposition
$$ G_K = \bigcup_{\overline g \in \Gal(M/K)} gG_M$$
from which we deduce that, as a $G_M$-module, $\Ind_{G_K}(V)$ admits
a decomposition as a direct sum
$$\Ind_{G_K}(V) = \bigoplus_{\overline g \in \Gal(M/K)} \Fonc_{cont}(gG_M,V) \simeq \bigoplus_{\Gal(M/K)}\Ind_{G_M}(V).$$
So that
$$H^1(G_M,\Ind_{G_K}(V)) \simeq \bigoplus_{\Gal(M/K)}H^1(G_M,\Ind_{G_M}(V))$$
and any of the $H^1(G_M,\Ind_{G_M}(V))$ is zero, because of the
first part of the proposition.

\bigskip

On the other hand, $\tau^\alpha$ topologically generates
$\Gal(L/K_{\infty})$, so that the complex
$$\mathbf \Ind_{G_{\infty}}(V) \mathop{\longrightarrow}^{\tau^\alpha-1} \Ind_{G_{\infty}}(V)$$
computes the cohomology $H^\bullet(\Gal(L/K_{\infty}),
\Ind_{G_{\infty}}(V))$. We get the kernel
$$\Ind_{G_{\infty}}(V)^{\Gal(L/K_{\infty})}\simeq \Ind_{\Gamma_K}(V).$$
The vanishing of $H^1(\Gal(L/K_{\infty}), \Ind_{G_{\infty}}(V))$
follows from the same arguments as for $H^1(G_L,\Ind_{G_K}(V))$
above.

\bigskip

Finally, the complex
$$\mathbf \Ind_{\Gamma_K}(V) \mathop{\longrightarrow}^{\gamma-1} \Ind_{\Gamma_K}(V)$$
computes the cohomology $H^\bullet(\Gamma_K, \Ind_{\Gamma_K}(V))$.
The surjectivity of $\gamma-1$ still comes from the nullity of
$H^1(\Gamma_K, \Ind_{\Gamma_K}(V))$ which is proved as before.
\Epartproof

\bigskip

From the surjectivity of $(\varphi-1)$ on $D_L(U)$, we immediately deduce that $F^3(U)=0$. \\
We also get the kernel of $\eta$:
$$\Ker \ \eta = \{(x,y,z) ; x,y \in D_L(U) \textrm{ and } z \in (1 -
\varphi)^{-1}((\tau^{\chi(\gamma)} - 1)(x) + (\delta - \gamma)(y))
\}.$$ Let $x,y \in D_L(U)$ and fix $x',y' \in D_L(U)$ such that
$$(1 - \varphi)(x')=x \textrm{ and } (1 - \varphi)(y')=y \ ;$$
proving that $F^2(U)=0$ consists then in proving
$$ \forall u \in \Ind_{G_{\infty}}(V), \
(x,y,(\tau^{\chi(\gamma)} - 1)(x') + (\delta - \gamma)(y') +u
\otimes 1) \in \textrm{Im } \beta.$$ But $(\tau^{\chi(\gamma)} - 1)$
is surjective on $\Ind_{G_{\infty}}(V)$, thus it suffices to
consider $\beta ( 0, x'+u',y')$ with $u'$ chosen so that
$(\tau^{\chi(\gamma)} - 1)(u')=u$.\\
Let $(u,v,w) \in \Ker(\beta)$, {\it i.e.} satisfying:
$$\left\{ \begin{array}{l}
(\gamma - 1)u = (\varphi - 1)v \\
(\tau - 1)u = (\varphi - 1)w \\
(\tau^{\chi(\gamma)} - 1)v = (\gamma - \delta)w \\
\end{array} \right.$$
Fix $x_0 \in D_L(U)$ such that $(\varphi - 1)x_0 = u$. Then the
first two relations show that
$$v_0 := v - (\gamma - 1)x_0 \textrm{ and } w_0 := w - (\tau - 1)x_0$$
lie in the kernel of $\varphi -1$ thus in $\Ind_{G_{\infty}}(V)$,
and satisfy furthermore:
$$(\tau^{\chi(\gamma)} - 1)v_0 = (\gamma - \delta)w_0.$$
Choose now $\eta \in \Ind_{G_{\infty}}(V)$ such that $(\tau - 1)\eta
= w_0.$ Then
$$(\tau^{\chi(\gamma)} - 1)(\gamma-1)\eta = (\gamma - \delta)(\tau - 1)\eta = (\tau^{\chi(\gamma)} - 1)v_0$$
so that $v_0 - (\gamma-1)\eta \in \Ind_{\Gamma_K}(V)$ and then there
exists $\varepsilon \in\Ind_{\Gamma_K}(V)$ such that
$$(\gamma-1)\varepsilon =v_0 - (\gamma-1)\eta$$
so that
$$(\gamma-1)(\eta+\varepsilon) = v_0$$
and
$$(\tau-1)(\eta+\varepsilon) = w_0.$$
Define then $x := x_0 + \eta + \varepsilon$ and let us verify
$\alpha(x) = (u,v,w)$:
$$\left.\begin{array}{ccccccc}
(\varphi - 1) x & = & (\varphi - 1)x_0 + (\varphi - 1)(\eta + \varepsilon) & = & (\varphi - 1)x_0 & = & u \\
(\gamma - 1) x & = & (\gamma - 1)x_0 + (\gamma - 1)(\eta + \varepsilon) & = &  v-v_0+v_0 & =  & v \\
(\tau - 1) x & = & (\tau - 1)x_0 + (\tau - 1)(\eta + \varepsilon) & =  & w-w_0+w_0 & = &  w \\
\end{array}\right.$$
which proves the proposition. \Eproof

\subsubsection{Explicit Formulas}
\paragraph{Proof of Theorem \ref{coho_gal} $ii)$} $ $\\
In order to make the isomorphism explicit, it suffices to do a
diagram chasing following the snake lemma: let
$$(x,y,z) \in Z^1(C_{\varphi, \gamma, \tau}(D_L(V))),$$
then through the injection $D_L(V) \hookrightarrow
D_L(\Ind_{G_K}(V)),$ we can see
$$(x,y,z) \in Z^1(C_{\varphi, \gamma, \tau}(D_L(\Ind_{G_K}(V)))).$$
From the nullity of $H^1(C_{\varphi, \gamma,
\tau}(D_L(\Ind_{G_K}(V))))$ we deduce the existence of a $b' \in
D_L(\Ind_{G_K}(V))$ such that
$$\alpha (b') = (x,y,z).$$
Consider now $\overline{b'} \in D_L(\Ind_{G_K}(V)/V)$ the reduction
of $b'$ modulo $D_L(V)$, then
$$\overline{b'} \in H^0(C_{\varphi, \gamma, \tau}(D_L(\Ind_{G_K}(V)/V))) =
(\Ind_{G_K}(V)/V)^{G_K}.$$ Thus, if $\tilde b \in \Ind_{G_K}(V)$
lifts  $\overline{b'}$, the image of $(x,y,z)$ in $H^1(K,V)$ is the
class of the cocycle
$$c : \sigma \mapsto c_{\sigma}=(\sigma -1)\tilde b.$$
But we can choose $\tilde b = b'-b$  since
$$(\varphi - 1)(b'-b) = x - x = 0$$
so that $b'-b \in \Ind_{G_K}(V)$ and then $b'-b$ lifts
$\overline{b'}$. So if
$$\sigma_{|_{G_{\infty}}} = \gamma^n\tau^m,$$
write
$$c_{\sigma} = (\sigma - 1)(b'-b) = -(\sigma - 1)(b) + (\gamma^n\tau^m-1) b' = -(\sigma - 1)b + \gamma^n
\frac {\tau^m-1}{\tau - 1}z + \frac{\gamma^n - 1 }{\gamma - 1}y$$
which concludes the proof of the theorem.\\
Let us finally show how to pass to the limit in order to get the
result for a representation which is not necessarily torsion. Let
$V$ be a $\Z_p$-adic representation of $G_K$. For all $r\geq 1$,
$$V_r = V\otimes \Z/p^r\Z$$
is a $p^r$-torsion representation such that
$$V = \lim_{\leftarrow} V_r.$$
Then we know that the continuous cohomology of $V$ can be expressed
as the limit:
$$\forall i\geq 0, \ H^i(K,V) = \lim_{\leftarrow} H^i(K,V_r) = \lim_{\leftarrow} F^i(V_r).$$
It suffices then to show
$$\forall i\geq 0,\  F^i(V) = \lim_{\leftarrow} F^i(V_r).$$
Let $H^i_r$ (respectively $B^i_r$, $Z^i_r$) denote the homology
group $H^i(C_{\varphi, \gamma, \tau}(D_L(V_r)))$ (respectively
$B^i(C_{\varphi, \gamma, \tau}(D_L(V_r)))$, $Z^i(C_{\varphi, \gamma,
\tau}(D_L(V_r)))$). The maps in the Herr complex are $\Z_p$-linear
so that in the category of $\Z_p$-modules there is an exact sequence
$$ 0 \rightarrow B^i_r \rightarrow  Z^i_r \rightarrow H^i_r \rightarrow 0$$
from which is obtained the exact sequence
$$ 0 \rightarrow \lim_{\leftarrow} B^i_r \rightarrow  \lim_{\leftarrow} Z^i_r \rightarrow \lim_{\leftarrow} H^i_r
\rightarrow {\lim_{\leftarrow}}^1 B^i_r $$ where
$\displaystyle{\lim_{\leftarrow}}^1$ is the first derived functor of
the functor $\displaystyle{\lim_{\leftarrow}}$. But for all $r$,
$$B^i_r \simeq B^i(C_{\varphi, \gamma, \tau}(D_L(V)))\otimes \Z/p^r\Z$$
so that the transition maps in the projective system $(B^i_r)$ are
surjective, and then this system satisfies Mittag-Leffler
conditions. Thus
$$ {\lim_{\leftarrow}}^1 B^i_r =0$$
shows that the homology of the inverse limit is equal to the inverse
limit of the homology, as desired.

\paragraph{The explicit formula for $H^2$}

The isomorphism from $H^2(C_{\varphi, \gamma, \tau}(D_L(V)))$ to
$H^2(K,V)$ can as well be made explicit:
\begin{prop}
The identification of Theorem \ref{coho_gal} between the homology of
$C_{\varphi, \gamma, \tau}(D_L(V))$ and the Galois cohomology of $V$
associates with a triple $(a,b,c) \in Z^2(C_{\varphi, \gamma,
\tau}(D_L(V)))$ the class of the $2$-cocycle:
$$ (g,h) \mapsto s_g -s_{gh}+ gs_h + \gamma^{n_1}\frac {\tau^{m_1}-1}{\tau - 1}
\frac{(\delta^{-1}\gamma)^{n_2} - 1 }{\delta^{-1}\gamma -
1}\delta^{-1}c$$ where $g_{|_{G_{\infty}}} =
\gamma^{n_1}\tau^{m_1}$, $h_{|_{G_{\infty}}} =
\gamma^{n_2}\tau^{m_2}$ and $s$ is a map $G_K \rightarrow \mathbf
A'\otimes V$ such that
$$s_{\sigma} = \phi\left(\frac{\gamma^n - 1 }{\gamma - 1}a +  \gamma^n \frac
{\tau^m-1}{\tau - 1}b\right)$$ where $\sigma_{|_{G_{\infty}}} =
\gamma^n\tau^m$ and $\phi$ is a continuous section of $\varphi-1$.
\end{prop}
\emph{Proof:} The proof is, mutatis mutandis, the same as the above
one and can be found in \cite{mathese}.

\paragraph{Remark} $\ $ \\
In the classical Herr complex case, with the class of $a$ is
associated the class of the $2$-cocycle:
$$ (g_1,g_2) \mapsto \tilde{\gamma}^{n_1}(h-1)\frac{\tilde {\gamma}^{n_2}-1}{\tilde{\gamma}-1}\tilde a$$
where $(\varphi - 1)\tilde a =a$, $\tilde \gamma$ is a fixed lift of
$\gamma$ in $G_K$ and $g_1 = \tilde \gamma ^{n_1}h$, $g_2 = \tilde
\gamma ^{n_2}h'$ with $h, h' \in G_{K_{\infty}}$ and $n_1, n_2 \in
\Z_p$.

\subsection{The cup-product}

\subsubsection{Explicit formulas for the cup-product}
In \cite{Herr_cup}, Herr gave explicit formulas for the cup-product
in terms of the complex associated with the representation. The
following theorem gives the formulas obtained in the metabelian
case:

\begin{thm}\label{cupprod}
Let $V$ and $V'$ be two $\Z_p$-adic representations of $G_K$, then
the cup-product induces maps:
\begin{enumerate}
\item Let $(a) \in H^0(C_{\varphi, \gamma, \tau}(D_L(V)))$ and $(a') \in
H^0(C_{\varphi, \gamma, \tau}(D_L(V')))$,
$$(a) \cup (a') = (a \otimes a') \in H^0(C_{\varphi, \gamma, \tau}(D_L(V\otimes
V'))),$$
\item let $(x,y,z) \in H^1(C_{\varphi, \gamma, \tau}(D_L(V)))$ and $(a')
\in H^0(C_{\varphi, \gamma, \tau}(D_L(V')))$,
$$(x,y,z) \cup (a') = (x \otimes a', y \otimes a', z \otimes a') \in H^1(C_{\varphi, \gamma, \tau}(D_L(V\otimes V'))),$$
\item let $(a) \in H^0(C_{\varphi, \gamma, \tau}(D_L(V)))$ and $(x',y',z')
\in H^1(C_{\varphi, \gamma, \tau}(D_L(V')))$,
$$(a) \cup (x',y',z') = (a \otimes x',a \otimes y',a \otimes z') \in H^1(C_{\varphi, \gamma, \tau}(D_L(V\otimes V')))$$
\item let $(x,y,z) \in H^1(C_{\varphi, \gamma, \tau}(D_L(V)))$ and
$(x',y',z') \in H^1(C_{\varphi, \gamma, \tau}(D_L(V')))$,\\
$(x,y,z)\cup (x',y',z')\in H^2(C_{\varphi, \gamma,\tau}(D_L(V\otimes
V')))$ can be written as:
$$( y \otimes\gamma x' - x \otimes \varphi y' \ , \
z \otimes \tau x' - x \otimes \varphi z' \ , \ \delta z \otimes
\tau^{\chi(\gamma)}y' - y \otimes \gamma z' + \Sigma_{z,z'})$$ where
$$ \Sigma_{z,z'} = \sum_{n \geq 1} \binom{\chi(\gamma)}{n+1} \sum_{k=1}^{n} \binom n k (\tau
-1)^{k-1}z \otimes \tau^k (\tau -1)^{n-k}z').$$
\end{enumerate}
\end{thm}

\subsubsection{Proof of Theorem \ref{cupprod}}

The only non trivial identity is the last one. We will use the
construction of the previous paragraph and we can then suppose that
$V$ and $V'$ are objects of $M_{G_K,p^r-tor}$. We will use the exact
sequences
$$0 \rightarrow V \rightarrow \Ind_{G_K}(V) \rightarrow V'' \rightarrow 0$$
and
$$0 \rightarrow F^0(V) \rightarrow F^0(\Ind_{G_K}(V)) \rightarrow F^0(V'') \rightarrow F^1(V) \rightarrow 0$$
and the cup-product property $da \cup b = d(a\cup b)$.\\
More precisely, fix $(x,y,z)$ and $(x',y',z')$ as in the theorem.
Then there exists an element $a \in D_L(\Ind_{G_K}(V))$ satisfying
$\alpha(a) = (x,y,z)$ and $\overline a \in (\Ind_{G_K}(V)/V)^{G_K}$.
Then $(x,y,z) \cup (x',y',z')$ is equal to
\begin{eqnarray*}
\alpha(a) \cup (x',y',z') & = & d(\overline a\otimes x',\overline
a\otimes y',\overline
a\otimes z') \ = \ \beta (a\otimes x',a\otimes y',a\otimes z') \\
& = & ((\gamma - 1)(a\otimes x') - (\varphi -1) (a\otimes y'),\\
&& (\tau -1)(a\otimes x') - (\varphi -1) (a\otimes z'),\\
&& (\tau^{\chi(\gamma)} - 1)(a\otimes y') - (\gamma -\delta)
(a\otimes z'))
\end{eqnarray*}

Now we use the formal identity
$$ (\sigma -1) (a\otimes b) = (\sigma-1)a \otimes \sigma b + a \otimes (\sigma -1)b.$$

The first term can be written as
\begin{eqnarray*}
(\gamma - 1)a\otimes x' &-& (\varphi -1) a\otimes y'\\
 & = & (\gamma - 1)a\otimes \gamma x' + a \otimes (\gamma - 1)x' -
(\varphi -1) a\otimes y' - a \otimes (\varphi - 1)y'\\
& = & y \otimes \gamma x' +  a \otimes ((\gamma - 1)x' - (\varphi - 1)y') - x \otimes y' \\
& = & y \otimes \gamma x' - x \otimes y'.
\end{eqnarray*}

From a similar computation, we get for the second one
\begin{eqnarray*}
(\tau -1)(a\otimes x') & - & (\varphi -1) (a\otimes z')\\
& = & (\tau - 1)a\otimes \tau x' + a \otimes (\tau - 1)x' -
(\varphi -1) a\otimes z' - a \otimes (\varphi - 1)z'\\
& = & z \otimes \tau x' +  a \otimes ((\tau - 1)x' - (\varphi - 1)z') - x \otimes z' \\
& = & z \otimes \gamma x' - x \otimes z'.
\end{eqnarray*}

Let us finally write the computation of the third term.

Iterating the identity
$$ (\sigma -1) (a\otimes b) = (\sigma-1)a \otimes \sigma b + a \otimes (\sigma -1)b,$$
we get by induction:
$$ (\sigma -1)^n (a\otimes b) = \sum_{k=0}^n \binom n k
(\sigma -1)^ka \otimes \sigma^k (\sigma -1)^{n-k}b.$$ First:
$$ (\tau^{\chi(\gamma)} - 1)a\otimes y' =
(\tau^{\chi(\gamma)} - 1)a\otimes \tau^{\chi(\gamma)}y' + a \otimes
(\tau^{\chi(\gamma)} - 1)y' = \delta z \otimes \tau^{\chi(\gamma)}y'
+ a \otimes (\gamma - \delta)z' $$ and
$$ (\gamma - 1)a\otimes z' = (\gamma  - 1)a\otimes \gamma z' + a \otimes (\gamma  - 1)z' =
y \otimes \gamma z' + a \otimes (\gamma  - 1)z'.$$ It remains to
compute $\delta (a\otimes z')$. Recall that
$$\delta = \frac{\tau^{\chi(\gamma)} - 1}{\tau - 1} = \sum_{n \geq 1} \binom{\chi(\gamma)}{n} (\tau -1)^{n-1}.$$
So
\begin{eqnarray*}
\delta (a\otimes z') & = & \sum_{n \geq 1} \binom{\chi(\gamma)}{n} (\tau -1)^{n-1}(a\otimes z')\\
& = & \sum_{n \geq 1} \binom{\chi(\gamma)}{n} \sum_{k=0}^{n-1} \binom{n-1} k (\tau-1)^ka \otimes \tau^k (\tau -1)^{n-1-k}z'\\
& = & a\otimes \delta z' + \sum_{n \geq 1} \binom{\chi(\gamma)}{n}
\sum_{k=1}^{n-1} \binom {n-1} k (\tau -1)^{k-1}z \otimes \tau^k
(\tau -1)^{n-1-k}z'.
\end{eqnarray*}
Which gives the result. \Eproof

\subsection{Kummer's map}\label{KummersMap}

In this paragraph, we suppose $p$ is odd and $\mathbf A' = \mathbf{\tilde A}$.\\
The purpose is to compute, in terms of the Herr complex, Kummer's
map
$$\kappa : K^* \rightarrow H^1(K, \Z_p(1)).$$
More precisely, let
$$F(Y) \in \left(W[[Y]][\frac 1 Y]\right)^\times,$$
we will compute a triple $(x,y,z) \in Z^1(C_{\varphi, \gamma,
\tau}(\mathbf{\tilde A}_L(1)))$ corresponding to the image $\kappa o
\theta (F(Y))$ of
$$\theta (F(Y)) = F(\pi) \in K^*.$$
Remark that there exist $d \in \Z$ and $G(Y) \in
\left(W[[Y]]\right)^\times$ such that
$$F(Y) = Y^{d}G(Y).$$
In fact $G(Y)$ can be written as the product of a $p$th root of
unity (which doesn't play any role) and a series in $1+(p) \subset
W[[Y]]$.

Denote
$$ \alpha = \theta (F(Y)) \in K^*.$$
Choose
$$\tilde \alpha = (\alpha_0,\alpha_1, \ldots , \alpha_n, \ldots) \in \mathbf{\tilde E}$$
such that $\alpha_0 =\alpha$. Then
$$ \frac {\tilde \alpha} {\rho^{d}} \in \mathbf{\tilde E}^+$$
thus
$$\frac {[\tilde \alpha]} {Y^{d}} \in \mathbf{\tilde A}^+$$
and for all $\sigma \in G_K,$ there exists $\psi_{\alpha}(\sigma)
\in \Z_p$ such that
$$ \sigma(\alpha) = \alpha \varepsilon^{\psi_{\alpha}(\sigma)}.$$
The map $ \sigma \mapsto \varepsilon^{\psi_{\alpha}(\sigma)}$ is in
fact a cocycle computing $\kappa(\alpha)$. So
$$ \sigma ([\tilde \alpha] ) = [\tilde \alpha] (1+X)^{\psi_{\alpha}(\sigma)} \textrm{ where  } \kappa(\alpha) =
\varepsilon^{\psi_{\alpha}} \in H^1(K, \Z_p(1)).$$

On the other hand, the series $\log \frac{[\tilde \alpha] }{F(Y)}$
converges in $B_{crys}$ and even in $\Fil^1B_{crys}$, namely $\frac
{[\tilde \alpha] }{F(Y)} \in \mathbf{\tilde A}^+$ and
$\theta\left(\frac {[\tilde \alpha] }{F(Y)}\right)=1$.

For all $h \in G_L,$
$$(h-1)\log \frac {[\tilde \alpha] }{F(Y)} = \psi_{\alpha}(h) t \textrm{ where } t =\log (1+X).$$
Define
$$ \tilde b = \log \frac {[\tilde \alpha] }{F(Y)} / t \in
\Fil^0B_{crys}.$$ Then
$$\psi_{\alpha}(h) = (h-1) (\tilde b) \ \forall h \in G_L.$$
Let
$$f(Y) = \fcnl(F) = \frac 1 p \log \frac {F(Y)^p}{\varphi(F(Y))}  \in W[[Y]]$$
then
$$(\varphi -1)(\tilde b) = \frac 1 t f(Y).$$
Choose $b_1 \in \mathbf{\tilde A}$ a solution of
$$(\varphi - 1)b_1 = -\frac {f(Y)} X.$$
Let $X_1 = \varphi^{-1}(X)=[\varepsilon^{\frac 1 p}]-1$, and $\omega
= \frac X {X_1} \in \mathbf{\tilde A}^+$ then
$$(\varphi - \omega)(b_1X_1) = -f(Y).$$
But reducing modulo $p$ this identity yields to an equation of the
form
$$T^p - \overline \omega T = -\overline{f(Y)}$$
and then by successive approximations modulo $p^m$, and because
$\mathbf{\tilde E}^+$ is integrally closed, $b_1X_1 \in
\mathbf{\tilde A}^+.$ But $\frac 1 {X_1} \in \Fil^0B_{crys}$, namely
the series
$$\frac t {X_1} = \sum_{n >0} (-1)^{n+1}\frac{\omega X^{n-1}}n = \sum_{n >0} (-1)^{n+1}\frac{\omega^n X_1^{n-1}}n $$
converges in $\Fil^1A_{crys}$, and thus
$$\frac 1 {X_1} = \frac t {X_1} \frac 1 t \in \Fil^0B_{crys}.$$
So
$$b_1=(b_1X_1) . \frac 1 {X_1} \in \Fil^0B_{crys}.$$
Moreover, $(\varphi - 1)b_2 = -\frac {f(Y)} 2$ admits a solution
$b_2$ in $\mathbf{\tilde A}^+$, so that if we set
$$x=-\frac {f(Y)} X -\frac {f(Y)} 2 \in \mathbf{\tilde A}_L$$
and choose a solution $b \in \mathbf{\tilde A}$ of $(\varphi - 1)b =
x$, then $b \in \Fil^0B_{crys}$.

So $\tilde b + b \in \Fil^0B_{crys}$ and
$$(\varphi - 1)(\tilde b + b) = (\frac 1 t - \frac 1 X - \frac 1 2)f(Y).$$
And we have the following lemma:

\begin{lem}\label{equation_mu}
Solutions of the equation
\begin{equation}
(\varphi -1)(\mu) = \left(\frac 1 t - \frac 1 X - \frac 1
2\right)f(Y) \label{equa_mu}
\end{equation}
in $\Fil^0 B_{crys}$ lie in $\Q_p+\Fil^1 B_{crys}$ and are invariant
under the action of ${G_L}$.
\end{lem}
{\it Proof of the lemma: } Consider
\begin{eqnarray*}
u & = & t\left(\frac 1 t - \frac 1 X - \frac 1 2\right)f(Y) =
\left(1 - \frac t X - \frac t 2\right)f(Y)\\
& = & -\sum_{n\geq 2} \frac {(-X)^n} {n+1}f(Y) + \sum_{n\geq 2}
\frac {(-X)^n} {2n} f(Y)
\end{eqnarray*}
then letting $\mu' = t\mu$, Equation \eqref{equa_mu} becomes
\begin{equation}
\left(\frac \varphi p -1\right)(\mu') = u \label{equa_mu'}
\end{equation}
but the sequences $\frac {(-X)^n} {n+1}f(Y)$ and $\frac {(-X)^n}
{2n} f(Y)$ converge to $0$ in $B_{crys}$ and
$$ \left(\frac \varphi p \right)^k \left( \frac {X^n} {n+1} \right) =  \frac {((1+X)^{p^k}-1)^n} {(n+1)p^k}$$
but
$$((1+X)^{p^k}-1) = \sum_{1 \leq r \leq p^k} \frac {p^k !}{(p^k-r)!} \frac {X^r}{r!} \in p^k A_{crys}$$
so
$$ \left(\frac \varphi p \right)^k \left( \frac {X^n} {n+1} \right) \in \frac {p^{k(n-1)}}{n+1} A_{crys}$$
converges to $0$ uniformly in $n$ in $B_{crys}$. The same holds for
$\left(\frac \varphi p \right)^k \left( \frac {X^n} {n+1} \right).$
So we get a solution $-\sum_{n\geq 0}\left(\frac \varphi
p\right)^nu$ of Equation \eqref{equa_mu'} in
$(\Fil^2B_{crys})^{G_L}$ thus a solution of Equation \eqref{equa_mu}
in $(\Fil^1B_{crys})^{G_L}$. And the fact that $$
(\Fil^0B_{crys})_{\varphi = 1} = \Q_p$$ proves the lemma.\Epartproof

So $b + \tilde b \in (\Fil^0B_{crys})^{G_L}$, thus, for all $h \in
G_L$,
$$(h-1)(-b) = (h-1) \tilde b = \psi_{\alpha}(h).$$
We conclude that there exist a unique $z \in \mathbf{\tilde A}_L(1)$
and $y \in \mathbf{\tilde A}_L(1)$ unique modulo $(\gamma-1)\Z_p(1)$
such that $\kappa(\alpha)$ is the image in $H^1(K, \Z_p(1))$ of the
triple
$$(x,y,z)\in Z^1(C_{\varphi, \gamma, \tau}(\mathbf{\tilde A}_L(1)))$$
where $x = -(\frac 1 X +\frac 1 2) f(Y)\otimes\varepsilon$. Namely,
we know that there exists such a triple $(x',y',z')$, and
$$x' - x \in (\varphi -1)\mathbf{\tilde A}_L(1)$$
which shows the existence, and $x$ being fixed, the unicity modulo
$\alpha(\Z_p)$ (where $\alpha$ is the first map in the Herr complex
$C_{\varphi, \gamma, \tau}(M)$, cf. section \ref{section_coho_gal}).

We get the more precise result:
\begin{prop}\label{hilbert}
Let $F(Y) \in \left(W[[Y]][\frac 1 Y]\right)^\times$. Then the image
of $F(\pi)$ by Kummer's map corresponds to the class of a triple
$$\left(-f(Y)\left(\frac 1 X +\frac 1 2\right), y, z\right)\otimes \varepsilon $$
with $y,z \in W[[X,Y]]$. This triple is congruent modulo
$XYW[[X,Y]]$ to
$$\left(-\frac{f(Y)} X-\frac{f(Y)} 2,\ 0 , \ Yd_{\log} F(Y)\right)\otimes \varepsilon $$
where $d_{\log}$ stands for the logarithmic derivative.
\end{prop}
\emph{Proof:} We have to show the congruences.

Remark that
\begin{eqnarray*}
\gamma\left(\frac {1 \otimes \varepsilon} X \right) & = & \frac {\chi(\gamma)\otimes \varepsilon}{\chi(\gamma)X + \frac{\chi(\gamma)(\chi(\gamma)-1)}2 X^2 +X^3u(X)}\\
& = & \left(\frac 1 X - \frac {(\chi(\gamma)-1)} 2 +
Xv(X)\right)\otimes \varepsilon
\end{eqnarray*}
so that
$$(\gamma -1)x \in XYW[[X,Y]](1)$$
where $\varphi^n$ is topologically nilpotent thus $\varphi-1$ is
invertible. Then it comes
$$y \in \Z_p(1) + XYW[[X,Y]](1).$$
Moreover, let $\tilde \gamma$ lift $\gamma$ in $G_K$, we still have
$$ (\tilde \gamma -1)(\tilde b \otimes \varepsilon) = \psi_\alpha(\tilde \gamma)$$
where, because of $ii)$ of Theorem \ref{coho_gal} on the one hand,
and Lemma \ref{equation_mu} above on the other hand,
$$(\tilde \gamma -1)(\tilde b \otimes \varepsilon + b \otimes \varepsilon) =
\psi_\alpha(\tilde \gamma) + (\tilde \gamma -1)(b\otimes
\varepsilon) = y \in \Fil^1B_{crys}(1)$$ which shows that
$$y \in XYW[[X,Y]](1).$$

\bigskip

We proceed as well for $z$:
$$(\tau-1)f(Y) = (f(Y(1+X))-f(Y)) = \sum_{n\geq 1} \frac {(XY)^n}{n!} f^{(n)}(Y) \equiv XY f'(Y) \mod (XY)^2.$$
Remark moreover
$$ \left(Y \frac d{dY}\right) \circ \frac \varphi p = \varphi \circ \left(Y \frac d{dY}\right)$$
so that
$$(\tau-1)f(Y) \equiv X (1-\varphi)\left(Y d_{\log}F(Y)\right)  \mod (XY)^2$$
and thus
$$(\tau-1)x \equiv (\varphi - 1) (Y d_{\log}F(Y)\otimes \varepsilon) \mod XYW[[X,Y]](1)$$
which shows
\begin{equation} \label{congr_z}
z \in Y d_{\log}F(Y)\otimes \varepsilon + \Z_p(1) + XYW[[X,Y]](1).
\end{equation}
And if $\tilde \tau$ lifts $\tau$ in $G_K$,
$$ (\tilde \tau -1)(\tilde b + b) = \psi_\alpha(\tilde \tau) - \log \frac {F(Y(1+X))}{F(Y)}/t +
(\tilde \tau-1) b \in \Fil^1 B_{crys}$$ so that
$$z = \psi_\alpha(\tilde \tau) + (\tilde \tau-1) b \in \log \frac {F(Y(1+X))}{F(Y)}/t + \Fil^1 B_{crys}$$
which, combined with \eqref{congr_z}, proves the desired result.
\Eproof

\section{Formal Groups}

In this section, we will prove the Br\"uckner-Vostokov explicit
formula for formal groups. In \cite{abr_hilb_forml}, Abrashkin
proved this formula under the condition that the $p^M$-th roots of
unity belong to the base field, which turns out not to be necessary.
To prove this formula without this assumption, we will explicitly
compute the Kummer map linked to the Hilbert symbol of a formal
group in terms of its ($\varphi$, $\Gamma$)-module, then compute the
cup-product with the usual Kummer map and the image of this
cup-product through the reciprocity isomorphism, which gives the
desired formula.

\subsection{Notation and backgrounds on formal groups} Consider
$G$ a connected smooth formal group over $W=W(k)$, the ring of Witt
vectors with coefficients in the finite field $k$. Denote by $K_0$
the fraction field of $W$ and $K$ a totally ramified extension of
$K_0$. Under these hypotheses, one can associate (cf.
\cite{FontaineDivisible}) with $G$ a formal group law which
determines $G$. Let us recall what it is.

\subsubsection{Formal group laws} Fix $p$ an odd prime and $d>0$ a number. Write
$\mathbf X = (X_1, \dots ,X_d)$, $\mathbf Y = (Y_1, \dots, Y_d)$ and
$\mathbf Z = (Z_1, \dots ,Z_d)$.
\begin{dfn}
A (commutative) \emph{formal group law} $\mathbf F$ of dimension $d$
on a commutative ring $R$ is the data of a $d$-uple of formal power
series
$$\mathbf F(\mathbf X,\mathbf Y) = (F_i(X_1, \dots ,X_d,Y_1,
\dots, Y_d))_{1\leq i\leq d} \in (R[[\mathbf X, \mathbf Y]])^d$$
satisfying
\begin{enumerate}
\item $\mathbf F(\mathbf X,\mathbf 0) = \mathbf F(\mathbf 0,\mathbf X) = \mathbf X,$
\item $\mathbf F(\mathbf X,\mathbf F(\mathbf Y,\mathbf Z)) = \mathbf F(\mathbf F(\mathbf X,\mathbf Y),\mathbf Z),$
\item $\mathbf F(\mathbf X,\mathbf Y) = \mathbf F(\mathbf Y,\mathbf X).$
\end{enumerate}
\end{dfn}
For a given formal group law $\mathbf F$, there exists a $d$-uple
$\mathbf f \in (R[[\mathbf X]])^d$ such that
$$\mathbf F(\mathbf X,\mathbf  f(\mathbf X)) = \mathbf F(\mathbf f(\mathbf X),\mathbf X) = \mathbf 0$$
so that on a given area where $\mathbf F$ and $\mathbf f$ converge
(for instance $\maxim_R^d$, when $R$ is a local ring with maximal
ideal $\maxim_R$, complete for the $\maxim_R$-adict opology),
$\mathbf F$ defines a commutative group structure, with identity
element $\mathbf 0$, the inverse of $\mathbf x$ being $\mathbf
f(\mathbf x)$. We then denote the group law by
$$\mathbf x+_F \mathbf y := \mathbf F(\mathbf x,\mathbf y).$$
Let $\mathbf G$ be another formal group law over $R$ of dimension
$d'$. Then a morphism from $\mathbf F$ to $\mathbf G$ is a $d'$-uple
of formal series $\mathbf h(\mathbf X) \in (R[[\mathbf X]])^{d'}$
with no constant term such that
$$\mathbf h(\mathbf F(\mathbf X,\mathbf Y)) = \mathbf G(\mathbf h(\mathbf X),\mathbf h(\mathbf Y)).$$
A morphism $\mathbf h$ is an isomorphism if $d=d'$ and if there
exists $\mathbf g(\mathbf X) \in (R[[\mathbf X]])^d$ with no
constant term satisfying
$$\mathbf h\circ \mathbf g(\mathbf X) = \mathbf g \circ \mathbf h(\mathbf X) =\mathbf X$$
or equivalently if $d\mathbf h(\mathbf 0) \in GL_d(R)$. Call
$\mathbf h$ a strict isomorphism when the normalization $d\mathbf
h(\mathbf 0)=I_d$ holds, \emph{i.e.}, if for all $1\leq i \leq d$,
$h_i(\mathbf X) \equiv X_i \mod \textrm{deg } 2$.

When $R$ is an algebra over $\Q$, then any formal group $\mathbf F$
admits a unique strict isomorphism, denoted by $\log_F$ from
$\mathbf F$ to the additive formal group law $\mathbf F_a : (\mathbf
X,\mathbf Y) \mapsto \mathbf X+\mathbf Y$. Call this isomorphism the
\emph{vectorial logarithm} of $\mathbf F$.\\
Coordinate maps of $\log_F$ form a basis of the \emph{logarithms} of
$\mathbf F$, the morphisms from $\mathbf F$ to the additive group on
$R$.

\subsubsection{$p$-adic periods}

Let us recall the notation of the first part: $K$ is a finite
extension of $\Q_p$ with residue field $k$ and $K_0 = W(k)[\frac 1
p]$. Fix $M\in \N$.
\begin{dfn}
If the isogeny $p\id_G : G \rightarrow G$ is finite and flat over
$W$ of degree $p^h$, then $G$ is said to be \emph{of finite height}
and $h$ is called the \emph{height} of $G$.
\end{dfn}

Let $G$ be a formal group over $W$ of dimension $d$ and finite
height $h$. Define
$$G[p^n] = \ker (p^n\id_G : G \rightarrow G)$$
the sub-formal group of $p^n$-torsion points of $G$ and denote
$$T(G) =\lim_{\leftarrow} G[p^n](\overline K )$$
the Tate module of $G$. Suppose moreover that
$$G[p^M](\overline K ) = G[p^M](K),$$
that is, suppose $p^M$-torsion points of $G$ lie in $K$.\\
Then $T(G)$ is a free $\Z_p$-module of rank $h$ and
$G[p^M](\overline K )= G[p^M](K)$ is isomorphic as a $\Z_p$-adic
representation of $G_K$ to $(\Z/p^M\Z)^h$.\\
The space of pseudo-logarithms of $G$ (on $K_0$) is the quotient
$$\left\{F \in K_0[[\mathbf X]] \ | \ F(\mathbf X +_G \mathbf Y) - F(\mathbf
X) -F(\mathbf Y) \in \mathcal O_{K_0}[[\mathbf X,\mathbf Y]]\otimes
\Q_p\right\} \ / \mathcal O_{K_0}[[\mathbf X]]\otimes \Q_p.$$ Denote
it by $H^1(G)$. It is a $K_0$-vector space of dimension $h$. The
space of logarithms of $G$ is
$$\Omega(G) = \left\{F \in {K_0}[[\mathbf X]] \ | \ F(\mathbf X +_G \mathbf Y) = F(\mathbf
X)  + F(\mathbf Y)\right\}.$$ It is naturally a sub-${K_0}$-vector
space of $H^1(G)$ of dimension $d$. Moreover, $H^1(G)$ admits the
filtration
$$\Fil^0(H^1(G)) =H^1(G), \ \ \Fil^1(H^1(G)) = \Omega(G), \ \ \Fil^2(H^1(G)) = 0.$$
With its filtration, and the Frobenius:
$$\varphi : F(\mathbf X) \mapsto F^{\varphi}(\mathbf X^p),$$
$H^1(G)$ is called the \emph{Dieudonn\'e module} of $G$.\\
In \cite{FontaineDivisible}, Fontaine showed there exists a pairing
$$H^1(G) \times T(G) \rightarrow B_{crys}^+$$
explicitly described by Colmez in \cite{ColmezPeriodesVarAb}.\\
It is defined as follows: let $\overline F \in H^1(G)$, and
$o=(o_s)_{s \geq 0} \in T(G)$ ; choose for all $s$ a lift $\hat o_s
\in W(\maximE)^d$ of $o_s$, \emph{i.e.} satisfying $\theta(\hat o_s)
= o_s$. Then the sequence $p^sF(\hat o_s)$ converges to an element
$\int_o d\overline F$ in $B_{crys}^+$ independent of the choice of
lifts $\hat o_s$ and $F$. Moreover, this pairing is compatible with
actions of Galois and $\varphi$ and with filtrations: if $F$ is a
logarithm, then
$\int_o dF \in \Fil^1B_{crys}^+$.\\
This pairing permits to identify $H^1(G)$ with
$\Hom_{G_{K_0}}(T(G),B_{crys}^+)$ with the filtration induced by the
one of $B_{crys}^+$. In order to work at an entire level, let us
introduce a lattice of $H^1(G)$, the $W$-module
$$D_{crys}^*(G) =\Hom_{G_{K_0}}(T(G), A_{crys})$$ endowed with the filtration
and the Frobenius $\varphi$ induced by those on $A_{crys}$. The
functor $D_{crys}^*$ is a contravariant version of the crystalline
functor of Fontaine's theory. The filtration is of length $1$ and we
denote
$$D^0(G) = D_{crys}^*(G) = \Hom_{G_{K_0}}(T(G),
A_{crys})$$ and
$$D^1(G) = \Fil^1D_{crys}^*(G) =\Hom_{G_{K_0}}(T(G), \Fil^1A_{crys}).$$
Then $D^1(G)$ is a direct factor of $D^0(G)$ of rank $d$. Fix then a
basis $\{l_1, \dots, l_d\}$ of $D^1(G)$ completed into a basis
$$\{l_1, \dots, l_d, m_1, \dots, m_{h-d}\}$$
of $D^0(G)$.\\
For all $1 \leq i \leq d$, $\varphi(l_i)$ takes values in
$$\varphi(\Fil^1 A_{crys})^d \subset (p
A_{crys})^d$$ so, $\frac{\varphi} p (l_i)$ takes values in $D^0(G)$.
Moreover, \cite{FontaineDivisible} and \cite{FontaineLaffaille} show
on the one hand that $\varphi$ is topologically nilpotent on
$D^0(G)$ (because $G$ is connected) and on the other hand that the
filtered module $D^0(G)$ satisfies
$$D^0(G) = \varphi D^0(G) + \frac{\varphi} p D^1(G).$$
So, define $\tilde \varphi $ the endomorphism of $D^0$ by
\begin{eqnarray*}
\tilde \varphi (l_i) & = & \frac{\varphi} p (l_i) \ \forall \ 1 \leq
i \leq d, \ \ \ \ \textrm{     and}\\
\tilde \varphi (m_i) & = & \varphi (m_i) \ \forall \ 1 \leq i \leq
h-d,
\end{eqnarray*}
then its matrix $\mathcal{E} \in GL_h(W)$.\\
Let $\mathbf{l} = {}^t(l_1,\dots, l_n)$ and $\mathbf{m} =
{}^t(m_1,\dots, m_{h-n})$, then
$$ \begin{pmatrix}
\frac \varphi p (\mathbf l) \\
\varphi (\mathbf m)
\end{pmatrix}
= \mathcal{E} \begin{pmatrix}
\mathbf l \\
\mathbf m
\end{pmatrix}.$$
So, we can write a block decomposition
$$\mathcal{E}^{-1} = \begin{pmatrix}
A & B \\
C & D \end{pmatrix}$$ so that
$$ \mathbf l = A \frac \varphi p
(\mathbf l) + B \varphi (\mathbf m)\textrm{ and } \mathbf m = C
\frac \varphi p (\mathbf l) + D \varphi (\mathbf m).$$ But $\varphi$
is topologically nilpotent on $D^0(G)$, and we can write
\begin{equation}\label{eqn_lm}
\mathbf l = \sum_{u \geq 1} F_u \frac {\varphi^u (\mathbf l)} p, \ \
\mathbf m = \sum_{u \geq 1} F'_u \frac {\varphi^u (\mathbf l)} p
\end{equation}
where
$$F_1 = A, \ F_2 = B\varphi(C), \ F_u = B\left(\prod_{1\leq k \leq
u-2}\varphi^k(D)\right)\varphi^{u-1}(C) \textrm{ for } u
>2,$$ and
$$F'_1 = C, \  F'_2 = D\varphi(C), \ F'_u = (\prod_{0\leq k \leq u-2}\varphi^k(D))\varphi^{u-1}(C).$$

\bigskip

Define a $\Z_p$-linear operator
$$\mathcal A = \sum_{u \geq1}F_u\varphi^u$$
on $K_0[[\mathbf X]]^d$. The vectorial formal power series
$$l_{\mathcal A}(\mathbf X) = \mathbf X + \sum_{m \geq 1} \frac
{\mathcal A^m(\mathbf X)}{p^m}$$ gives then the vectorial logarithm
of a formal group $F$ from which we can recover the formal group law
$\mathbf F$ by:
$$ \mathbf F(\mathbf X, \mathbf Y) = l_{\mathcal A}^{-1}\left( l_{\mathcal
A}(\mathbf X) + l_{\mathcal A}(\mathbf Y) \right).$$ In
\cite{HondaFG}, Honda introduced the \emph{type} of a logarithm. A
logarithm $\log$ is of type $u \in M_d(W)[[\varphi]]$ if $u$ is
special, i.e. $u \equiv pI_d \mod \varphi$ and if
$$ u(\log) \equiv 0 \mod p.$$
We remark that $pI_d - \mathcal A$ is special and that, by
construction, $l_{\mathcal A}$ is of type $pI_d - \mathcal A$.
Moreover, $\mathbf l$ is also of type $pI_d - \mathcal A$ because of
Equation \eqref{eqn_lm}.

Furthermore, Honda showed in \cite[Theorem 2]{HondaFG} that two
formal groups with vectorial logarithms of the same type are
isomorphic over $W$. Thus, we can replace the study of the formal
group $G$ by the one of $F$, which is easier because we know an
explicit expression of its logarithms, which gives us a control on
denominators.

\subsection{Properties of the formal group $F$}\label{RappelsF}

In this section, the reader can refer to \cite{abr_hilb_forml} from
which we recall principal constructions.

\bigskip

Let us first describe the Dieudonn\'e module of $F$.

We already know a basis of the logarithms, the coordinate power
series of the vectorial series
$$l_{\mathcal A}(\mathbf X) = \mathbf X + \sum_{m \geq 1} \frac
{\mathcal A^m(\mathbf X)}{p^m}.$$ Complete it into a basis of
$H^1(F)$ by putting
$$m_{\mathcal A}(\mathbf X) =\sum_{u \geq 1} F'_u \frac {\varphi^u (l_{\mathcal A}(\mathbf X))} p.$$
Let $o=(o_s)_{s \geq 0} \in T(F)$. For all $s\geq 0$, choose a lift
$\hat o_s \in W(\maximE)^d$ of $o_s$, that is, with $\theta(\hat
o_s) = o_s$. Then the following lemma says that the sequence
$p^s\id_F \hat o_s$ converges in $W^1(\maximE)^d$ towards an element
$j(o)$ independent of the choice of lifts:

\begin{lem}
\begin{enumerate}
\item The series $l_{\mathcal A}$ defines a continuous one-to-one homomorphism
of $G_{K}$-modules
$$l_{\mathcal A} \ : \ F(W(\maximE)) \rightarrow A_{crys}^d\otimes_{\Z_p}\Q_p.$$
Moreover, the restriction of $l_{\mathcal A}$ to $F(W^1(\maximE))$
takes values in $(\Fil^1A_{crys})^d.$
\item The endomorphism $p\id_F$ of $F(W(\maximE))$ is topologically nilpotent. The convergence of $p\id_F$ to zero
is uniform on $F(W^1(\maximE))$.
\item The map $j \ : \ T(F) \rightarrow W^1(\maximE)^d$ is well
defined and provides a continuous one-to-one homomorphism of
$G_{K}$-modules $j \ : \ T(F) \rightarrow F(W^1(\maximE))$.
\end{enumerate}
\end{lem}
\emph{Proof:}
Point $1.$ is Lemma 1.5.1 of \cite{abr_hilb_forml}.\\
Point $2.$ follows from the fact that $W^1(\maximE)= \omega
W(\maximE)$ with $\omega = X/\varphi^{-1}(X) \in W(\maximE) +
p\mathbf {\tilde A}^+$ and that the series corresponding to $p\id_F$
can be written as
$$p\id_F \mathbf X = p\mathbf X + \textrm{higher degrees.}$$
Let us recall briefly the proof of Point $3.$:\\
For all $s\geq 0$,
$$\theta(p^s\id_F\hat o_s)=o_0=0$$
so that $p^s\id_F\hat o_s \in F(W^1(\maximE))$. On the other hand,
for all $s\geq 0$,
$$p\id_F \hat o_{s+1} \equiv \hat o_s \mod F(W^1(\maximE))$$
thus
$$p^{s+1}\id_F\hat o_{s+1} \equiv p^s\id_F\hat o_s \mod p^s\id_F\left(F(W^1(\maximE))\right)$$
And Point $2.$ provides the convergence of the sequence $(p^s\id_F\hat o_s)_s$.\\
The fact that the convergence is given without compatibility
condition on the lifts shows the independence of the limit with
respect to the choice of these lifts. Namely, let $(\hat o_s)_{s\geq
0}$ and $(\hat o_s')_{s\geq 0}$ be two given lifts of $(o_s)_{s\geq
0}$, then for any lift $(\hat o_s'')_{s\geq 0}$ where
$$\forall s\geq 0,\ \hat o_s'' = \hat o_s\textrm{ or } \hat o_s',$$
we still have the convergence of
$(p^s\id_F\hat o_s'')_s$, from which we deduce that the limits are the same.\\
The remainder is straightforward.\Eproof

\bigskip

Composing the vectorial logarithm $l_{\mathcal A}$ with $j$ gives a
$G_K$-equivariant injection that we will denote by $\mathbf l$ from
$T(F)$ into $(\Fil^1A_{crys})^d$. This map satisfies then for any
$o$ in $T(F)$:
$$\mathbf l (o) = l_{\mathcal A} (\lim_{s \to \infty} p^s \id_F \hat o_s) =
 \lim_{s \to \infty} p^s l_{\mathcal A}(\hat o_s).$$
Put now
$$\displaystyle{\mathbf m = \sum_{u \geq 1}
F'_u \frac {\varphi^u (\mathbf l)} p},$$ then $\begin{pmatrix}
\mathbf l \\
\mathbf m
\end{pmatrix}$ provides a basis of $D^0(F)$ with $ \mathbf l$ a basis of
$D^1(F)$. The map
$$\begin{pmatrix} \mathbf l \\ \mathbf m \end{pmatrix} : T(F) \rightarrow A_{crys}^h$$
then factorizes through
 $$\begin{pmatrix}
l_{\mathcal A} \\
m_{\mathcal A}
\end{pmatrix} \ : \ F(W^1(\maximE)) \rightarrow A_{crys}^h.$$
This map is the period pairing. Recall (cf. \cite{abr_hilb_forml},
Remark 1.7.5) that this map takes values in $\mathbf{\tilde
A}^+[[X^{p-1}/p]]$. It is also a consequence of Wach's computation
for potentially crystalline representations (cf.
\cite{Wach_pot_cris}).

\bigskip

Fix now a basis $(o^1, \dots, o^h)$ of $T(F)$. We can then introduce
the period matrix
$$ \mathcal V = \begin{pmatrix}
\mathbf l (o^1) & \dots & \mathbf l (o^h) \\
\mathbf m (o^1) & \dots & \mathbf m (o^h)
\end{pmatrix} \in M_h(\mathbf{\tilde A}^+[[X^{p-1}/p]]) \cap GL_h(\Frac \mathbf{\tilde A}^+[[X^{p-1}/p]]).$$
It satisfies
$$\begin{pmatrix} I_d\frac \varphi p & 0\\
0 & I_{h-d}\varphi
\end{pmatrix} \mathcal V = \mathcal {E V}.$$
Remark that the inverse of $\mathcal V$ is then the change of basis
matrix from the basis $(o^1, \dots, o^h)$ to a basis of
$$D_{crys}(T(F)) = (T(F)\otimes_{\Z_p}A_{crys})^{G_{K_0}},$$
the  covariant version of the crystalline module of Fontaine's
theory associated with $T(F)$.

\bigskip

Let $u \in T(F) \otimes A_{crys}$, and $U$ be the coordinate vector
of $u$ in the basis $(o^1, \dots, o^h)\mathcal V^{-1}$ of
$D_{crys}(T(F))$, then we can compute the coordinates of
$$\varphi(u) = (o^1, \dots, o^h)\varphi(\mathcal
V^{-1})\varphi(U).$$ We know that
$$\varphi(\mathcal V) = \begin{pmatrix} pI_d & 0 \\ 0 & I_{h-d}
\end{pmatrix} \begin{pmatrix} I_d\frac \varphi p & 0 \\
0 & I_{h-d}\varphi \end{pmatrix} \mathcal V = \begin{pmatrix} pI_d &
0 \\ 0 & I_{h-d}
\end{pmatrix}\mathcal {E V}$$ so that
$$\varphi(\mathcal V^{-1}) =\mathcal V^{-1}
\mathcal E^{-1}\begin{pmatrix} p^{-1}I_d & 0 \\ 0 & I_{h-d}
\end{pmatrix}$$ and coordinates of $\varphi(y)$ in the basis
$(o^1, \dots, o^h)\mathcal V^{-1}$ are then
$$\mathcal E^{-1}\begin{pmatrix} I_d\frac \varphi p & 0 \\
0 & I_{h-d}\varphi \end{pmatrix} U.$$
Keeping this in mind, the following lemma shows that $\begin{pmatrix} \frac {\mathcal A} p & 0 \\
0 & I_{h-d} \end{pmatrix}$ acts as the Frobenius on
$D_{crys}(T(F))$.

\begin{lem}\label{frob_form}
The following equality holds:
$$\mathcal E^{-1}\begin{pmatrix} \frac
\varphi p \circ l_{\mathcal A}\\ \varphi \circ m_{\mathcal
A}\end{pmatrix}
= \begin{pmatrix} \frac {\mathcal A} p \circ l_{\mathcal A} \\
m_{\mathcal A}
\end{pmatrix}$$

\end{lem}
\emph{Proof:} Compute:
$$A \frac \varphi p (l_{\mathcal A}) + B \varphi (m_{\mathcal A}) =
A \frac \varphi p (l_{\mathcal A}) + \sum_{u \geq 1} B\varphi F'_u
\frac {\varphi^u (l_{\mathcal A})} p = \frac {\mathcal A } p
(l_{\mathcal A})$$ for $B\varphi F'_u = F_{u+1}$ for all $u \geq 1$.
And:
$$C \frac \varphi p (l_{\mathcal A}) + D \varphi (m_{\mathcal A}) =
C \frac \varphi p (l_{\mathcal A}) + \sum_{u \geq 1} D\varphi F'_u
\frac {\varphi^u (l_{\mathcal A})} p = m_{\mathcal A}$$ since
$D\varphi F'_u = F'_{u+1}$ for all $u \geq 1$. \Eproof

Abrashkin also computed the cokernel of injection $j$ (cf.
\cite[Proposition 2.1.]{abr_hilb_forml}):
\begin{prop}\label{s_ex_j}
There is an equality $$(\mathcal A- p)\circ l_{\mathcal
A}(F(W(\maximE))) = (\mathcal A- p)\circ l_{\mathcal
A}(F(W^1(\maximE)))$$ and the following sequence is exact:
$$\xymatrix{0 \ar[r] & T(F) \ar[r]^-{j} & F(W^1(\maximE))\ar[r]^-{\left(\frac {\mathcal A}p-1\right)\circ l_{\mathcal A}} &
W(\maximE)^d \ar[r] & 0\\}$$
\end{prop}

\paragraph{Remark} Beware that if $x \in F(W(\maximE))$,
$$ \varphi(l_{\mathcal A})(x) = \varphi(l_{\mathcal A}(x)) $$
and then
$$ \mathcal A(l_{\mathcal A})(x) = \mathcal A(l_{\mathcal A}(x)) $$
hold if $\varphi(x)=x^p$ (e.g. when $x$ is a Teichm\"uller
representative) but not in general ! On the left side, $\varphi$ and
$\mathcal A$ act on $W[[\mathbf X]]$, whereas they act on $A_{crys}$
on the right side.

\bigskip

Abrashkin showed furthermore (cf. \cite[Lemma
1.6.2.]{abr_hilb_forml})
\begin{lem}
$F(\maximE)$ is uniquely $p$-divisible.
\end{lem}
This provides a continuous one-to-one $G_K$-equivariant morphism
$$\delta \ : \ F(\maximE) \rightarrow F(W(\maximE))^{(\mathcal A -p)\circ l_{\mathcal A}=0}$$
defined as follows: let $x \in F(\maximE)$, then because of the
lemma, for all $s\geq 0$ there exists a unique $x_s \in F(\maximE)$
such that
$$p^s\id_Fx_s=x.$$
Thus the sequence $(p^s\id_F[x_s])_s$ converges to an element $\delta(x)$ in $F(W(\maximE))$.\\
$\delta$ is a morphism since
$$\delta(x+_Fy) = \lim_s p^s\id_F[x_s+_Fy_s] =  \lim_s p^s\id_F([x_s]+_F[y_s]+_Fu_s)$$
with $u_s \in pW(\maximE)$ where the convergence of $p^s\id_F$ towards zero is uniform.\\
Moreover, for $\mathcal A \circ l_{\mathcal A}$  coincides with
$\mathcal A(l_{\mathcal A})$ on Teichm\"uller representatives, we
get the last point:
$$(\mathcal A -p)\circ l_{\mathcal A}(\delta(x)) = (\mathcal A -p)(l_{\mathcal A})(\delta(x)) = 0.$$
Finally, remark
$$\theta(\delta(x)) = \theta([x]).$$
Namely, for all $s\geq 0$,
$$\theta(p^s\id_F[x_s]) = p^s\id_F\theta([x_s])  =\theta([x]).$$

\subsection{The ring $\mathcal G_{[b,a]}$ and some subrings.}
\subsubsection{Introducing the objects}

Fix $e$ the absolute ramification index of $K$.

In \cite{Berger_equadiff}, Berger introduced for $s\geq r \geq 0$
the ring $\mathbf{\tilde A}_{[s,r]}$, the $p$-adic completion of the
ring
$$\mathbf{\tilde A}^+\left[\frac p {Y^{rep/(p-1)}},\frac {Y^{sep/(p-1)}} p\right].$$
Let us then introduce for $a>b \geq 0$, the ring
$$\mathcal G_{[b,a]} := \mathbf{\tilde A}^+\left[\left[ \frac {Y^{ae}} p, \frac p {Y^{be}}\right]\right]$$
which for integers $a$ and $b$ admits the description
$$\mathcal G_{[b,a]} = \left\{ \sum_{n\in \Z}a_nY^n\ | \ a_n \in \mathbf{\tilde A}^+\left[ \frac 1 p \right], \left.
\begin{array}{ll}
ae v_p(a_n) +n \geq 0 & \textrm{ for } n\geq 0\\
be v_p(a_n) +n \geq 0 & \textrm{ for } n\leq 0
\end{array}\right.\right\}.$$
Note that the expression $\sum_{n\in \Z}a_nY^n$ for an element of
$\mathcal G_{[b,a]}$ is not unique. The ring $\mathcal G_{[b,a]}$ is
naturally, for $a>\alpha\geq \beta > b$ a subring of $\mathbf{\tilde
A}_{[\alpha(p-1)/p,\beta(p-1)/p]}$. We even have inclusions
$$\mathbf{\tilde A}_{[a(p-1)/p,b(p-1)/p]} \subset
\mathcal G_{[b,a]} \subset \mathbf{\tilde
A}_{[\alpha(p-1)/p,\beta(p-1)/p]}.$$ Endow then $\mathcal G_{[b,a]}$
with the induced topology, which is well defined since inclusions
$$\mathbf{\tilde A}_{[r_1,s_1]} \hooklongrightarrow \mathbf{\tilde A}_{[r_2,s_2]}$$
for $r_1 \leq r_2 \leq s_2 \leq s_1$ are continuous.

This topology then admits as a basis of neighborhoods of zero
$$ \left\{ \left\{ \sum_{n > N} a_n \left(\frac{Y^{ae}} p\right)^n + \sum_{n > N} b_n \left(\frac p {Y^{be}} \right)^n ;
a_n, b_n \in \mathbf{\tilde A}^+ \right\} + p^k \mathcal G_{[b,a]}
\right\}_{N,k \in \N}.$$

For $a_1>a_2>b_2>b_1 \geq 0$ we still have continuous injections
$$\mathcal G_{[a_1,b_1]} \hooklongrightarrow  \mathcal G_{[a_2,b_2]}.$$
Denote then for $a \geq b \geq 0$ by $\mathcal G_{[b,a[}$ the
$p$-adic completion of $\bigcup_{\alpha>a} \mathcal G_{[b,\alpha]}.$
For integers $a$ and $b$, it admits the description:
$$\mathcal G_{[b,a[} = \left\{ \sum_{n\in \Z}a_nY^n; a_n \in \mathbf{\tilde A}^+\left[ \frac 1 p \right], \left.
\begin{array}{ll}
ae v_p(a_n) +n >0 & \textrm{for } n > 0 \textrm{ and } \\
ae v_p(a_n) +n \mathop{\longrightarrow}\limits_{n \rightarrow + \infty} +\infty, & \\
be v_p(a_n) +n \geq 0 & \textrm{for } n\leq 0
\end{array}\right.\right\}.$$
Because of the inclusion
$$\bigcup_{\alpha>a} \mathcal G_{[b,\alpha]} \hooklongrightarrow
\bigcup_{\alpha>a, \beta > b} \mathbf{\tilde
A}_{[\alpha(p-1)/p,\beta(p-1)/p]}$$ we endow $\mathcal G_{[b,a[}$
with the topology induced by the $p$-adic toplogy of the $p$-adic
completion of $\bigcup_{\alpha>a, \beta > b} \mathbf{\tilde
A}_{[\alpha(p-1)/p,\beta(p-1)/p]}$.

Let us also introduce for $b \geq 0$,
$$ \mathcal G_{[b,\infty[} :=  \bigcap_{a>b} \mathcal G_{[b,a]} =
\mathbf{\tilde A}^+\left[\left[ \frac p {Y^{eb}} \right]\right]
\subset \mathbf{\tilde A}$$ which is for $b$ integer
$$\mathcal G_{[b,\infty[} =   \left\{ \sum_{n\leq 0}a_nY^n; a_n \in \mathbf{\tilde A}^+, be v_p(a_n) +n \geq 0 \textrm{ for } n \leq 0\right\}.$$

Remark that the Frobenius
$$ \varphi_{\mathcal G}\left( \sum_{n <0}a_nY^{aen} +\sum_{n\geq 0}a_nY^{ben}\right) =
\sum_{n <0}\varphi(a_n)Y^{paen} + \sum_{n \geq
0}\varphi(a_n)Y^{pben}$$ defines a one-to-one morphism from
$\mathcal G_{[b,a]}$ (respectively $\mathcal G_{[b,a[}$) into
$\mathcal G_{[pb,pa]}$ (respectively $\mathcal G_{[pb,pa[}$).

Introduce for integers $a$ and $b$ the subring of $\mathcal
G_{[b,a]}$:
\begin{eqnarray*}
\mathcal G_{Y, [b,a]} & := & W[[Y]]\left[\left[ \frac {Y^{ae}} p, \frac p {Y^{be}}\right]\right]\\
& = & \left\{ \sum_{n\in \Z}a_nY^n\ ; \ a_n \in K_0, \left.
\begin{array}{ll}
ae v_p(a_n) + n \geq 0 & \textrm{ for } n\geq 0\\
be v_p(a_n) + n \geq 0 & \textrm{ for } n\leq 0
\end{array}\right.\right\}
\end{eqnarray*}
and $\mathcal G_{Y, [b,a[}$ the subring of $\mathcal G_{[b,a[}$
admitting the description
$$ \mathcal G_{Y, [b,a[} = \left\{ \sum_{n\in \Z}a_nY^n\ ; \ a_n \in K_0, \left. \begin{array}{ll}
ae v_p(a_n) + n > 0 & \textrm{ for } n> 0 \textrm{ and }\\
ae v_p(a_n) + n \mathop{\longrightarrow}\limits_{n \rightarrow + \infty} +\infty, & \\
be v_p(a_n) + n \geq 0 & \textrm{ for } n\leq 0
\end{array}\right.\right\}.$$

Finally, for $b \geq 0$,
\begin{eqnarray*}
\mathcal G_{Y, [b,\infty[} & := & \bigcap_{a>b} \mathcal G_{Y, [b,a]} = W[[Y]]\left[\left[ \frac p {Y^{eb}} \right]\right] \\
& = &  \left\{ \sum_{n\leq 0}a_nY^n; a_n \in K_0, be v_p(a_n) +n
\geq 0 \right\}.
\end{eqnarray*}

Contrary to the above situation, the expression $\sum_{n\in
\Z}a_nY^n$ is unique as it is shown in the following

\begin{lem}\label{alg_G_intersections}
\begin{enumerate}
\item In $\mathcal G_{[b,a]}\left[\frac 1 p\right]$, one has $\mathcal G_{[0,a]}\left[\frac 1 p\right] \bigcap \mathcal G_{[b,\infty[}\left[\frac 1 Y\right] = \mathbf{\tilde A}^+$.
\item An element of $\mathcal G_{Y, [b,a]}$ or $\mathcal G_{Y, [b,a[}$ can uniquely be written as
$\sum_{n\in \Z}a_nY^n$ with $a_n\in K_0$.
\item For $a \geq \alpha \geq \beta \geq b$, and $)$ designating $]$ or $[$, one has
$$ \mathcal G_{Y, [\beta,\alpha]}\left[\frac 1 p\right] \bigcap \mathcal G_{[b,a)} = \mathcal G_{Y, [b,a)}.$$
\end{enumerate}
\end{lem}
\emph{Proof:} The first point can be shown in Berger's rings
$\mathbf{\tilde A}_{[s,r]}$, in fact in the ring $\mathbf{\tilde
A}_{[s,\infty[}\left[\frac 1 Y \right] + \mathbf{\tilde
A}_{[0,r]}\left[\frac 1 p \right]$. Any element of this ring is of
the form $\sum_{n\in \N} p^n (\frac {x_n}{Y^k}-\frac{y_n}{p^l})$
with $x_n \in \mathbf{\tilde A}^+\left[\frac p
{Y^{rep/(p-1)}}\right]$ and $y_n \in \mathbf{\tilde A}^+\left[\frac
{Y^{sep/(p-1)}} p\right]$. Such an element is zero when
$$p^l\sum_{n\in \N} p^n x_n = Y^k\sum_{n\in \N} p^n y_n \in \mathbf{\tilde A}_{[s,\infty[}\bigcap \mathbf{\tilde A}_{[0,r]}.$$
The condition is that for all $N \in \N$,
$$\sum_{n<N} p^n (p^lx_n-Y^ky_n) \in p^N\mathbf{\tilde A}_{[s,r]}.$$
That is
$$p^l\sum_{n<N} p^n x_n \in \mathbf{\tilde A}^+ + p^N\mathbf{\tilde A}^+\left[\frac p {Y^{rep/(p-1)}}\right]$$
and then
$$\sum_{n<N} p^n x_n \in \mathbf{\tilde A}^+ + p^{N-l}\mathbf{\tilde A}^+\left[\frac p {Y^{rep/(p-1)}}\right],$$
and similarly
$$\sum_{n<N} p^n y_n \in \mathbf{\tilde A}^+ + p^{N-l/se}\mathbf{\tilde A}^+\left[\frac {Y^{sep/(p-1)}} p\right].$$
Thus the limit $p^l\sum_{n\in \N} p^n x_n = Y^k\sum_{n\in \N} p^n
y_n$ lies in $p^l\mathbf{\tilde A}^+ \bigcap Y^k\mathbf{\tilde A}^+
= p^lY^k\mathbf{\tilde A}^+$, hence
$$\sum_{n\in \N} p^n \frac {x_n}{Y^k}= \sum_{n\in \N} p^n \frac{y_n}{p^l} \in \mathbf{\tilde A}^+$$
as claimed.

Because of the first point, it is enough to prove the second one for
$\sum_{n < 0}a_nY^n$ and $\sum_{n > 0}a_nY^n$. It is to prove that
such a series converges to zero if and only if all the $a_n$
actually are zero. In the first case, it is a series converging in
$\mathbf{\tilde A}$ and the natural map $\mathcal G_{Y,
[\beta,\infty[} \rightarrow \mathbf{\tilde A}$ is a continuous
one-to-one morphism. Successive approximations modulo $p^n$ and
modulo $Y^k$ then provide the result. On the other side, $\mathcal
G_{Y, [0,\alpha)}$ is naturally a subring of the separable
completion of $\mathbf{\tilde A}\left[\frac 1 p \right]$
for the $Y$-adic topology. The result then follows from successive reductions modulo $Y^k$.\\
We use similar techniques to show the last point. Again, because of
the first one, it suffices to prove both $\mathcal G_{Y,
[\beta,\infty[}\left[\frac 1 p\right] \bigcap \mathcal
G_{[b,\infty[} = \mathcal G_{Y, [b,\infty[}$ and $\mathcal G_{Y,
[0,\alpha]}\left[\frac 1 p\right] \bigcap \mathcal G_{[0,a)} =
\mathcal G_{Y, [0,a)}.$ First consider then $x = \sum_{n \leq
0}a_nY^n \in \frac 1 {p^\lambda}\mathcal G_{Y, [\beta,\infty[}$ with
$\beta e v_p(a_n) +n+\lambda \geq 0$ for all $n$. We suppose
furthermore that $x$ belongs to $\mathcal G_{[b,\infty[}$, that is,
it can be written as $ \sum_{n\in \N} b_n\frac {p^n}{Y^{ebn}}$ with
$b_n \in \mathbf{\tilde A}^+$. The identity
$$\sum_{n \leq 0}a_nY^n =\sum_{n\in \N} b_n\frac {p^n}{Y^{ebn}}$$
makes sense in $\frac 1 {p^\lambda}\mathcal G_{[\beta,\infty[}$,
thus in $\mathbf{\tilde A}$. Denote by $n_0$ the highest integer,
supposing it exists, satisfying $be v_p(a_{n_0}) +n_0 < 0$. We can
then suppose that the identity above is of the form
$$\sum_{n \leq n_0}a_nY^n =\sum_{n\in \N} b_n\frac {p^n}{Y^{ebn}}.$$
Multiplying by $Y^{bev_p(a_{n_0})}$ and reducing modulo
$p^{v_p(a_{n_0})}$ yields then to
$$\sum_{n = n_0'}^{n_0} a_nY^{n+bev_p(a_{n_0})} \equiv \sum_{n=0}^{v_p(a_{n_0})} b_n p^nY^{eb(v_p(a_{n_0})-n)} \mod p^{v_p(a_{n_0})}$$
but the right term is entire (it belongs to $\mathbf{\tilde A}^+$) and not the left one, whence a contradiction.\\
Now, consider as before an identity of the form
$$\sum_{n \geq 0}a_nY^n =\sum_{n\in \N} b_n\frac {Y^{ean}}{p^n}$$
and denote by $n_0$ the lowest integer satisfying $ae v_p(a_{n_0})
+n_0 < 0$. It can be reduced to an identity of the form
$$\sum_{n \geq n_0}a_nY^n =\sum_{n\in \N} b_n\frac {Y^{ean}}{p^n}.$$
Multiplying by $p^{-v_p(a_{n_0})}$ and reducing modulo $Y^{n_0+1}$
yields to
$$ p^{-v_p(a_{n_0})}a_{n_0}Y^{n_0} \equiv \sum_{0 \leq n \leq \frac{n_0}{ea}} b_n p^{-v_p(a_{n_0})-n}Y^{ean} \mod Y^{n_0+1}$$
and the contradiction comes from the inequality
$$n\leq \frac{n_0}{ea} < -v_p(a_{n_0})$$
hence the right term is divisible by $p$, and not the left one.\\
The case of $\mathcal G_{Y, [0,\alpha]} \bigcap \mathcal G_{[0,a[} =
\mathcal G_{Y, [0,a[}$ follows from a similar argument. \Eproof

\paragraph{Remark}
As said before, periods of formal groups lie in the ring
$$\mathbf{\tilde A}^+[[X^{p-1}/p]] = \mathbf{\tilde A}^+[[Y^{pe}/p]] = \mathcal G_{[0,p]}.$$
We can also recover some well known rings by
\begin{eqnarray*}
\mathbf{\tilde A}^+ & = & \mathcal G_{[0,\infty[}\\
\mathbf B_{rig}^{\dagger} & = & \bigcup_{b>0} \mathcal G_{[
b,\infty[}.
\end{eqnarray*}

\subsubsection{Some topological precisions }

\begin{lem}
\begin{enumerate}
\item Finite sums
$$\left\{ \sum_{n=0}^N a_n\left(\frac {Y^{ea}} p \right)^n + b_n\left(\frac p {Y^{eb}} \right)^n ;
a_n, b_n \in \mathbf{\tilde A}^+, N \in \N \right\}$$ form a dense
subset of $\mathcal G_{[b,a]}$. The same holds for the sub-algebra
$$\mathcal G_{[b,a]}\bigcap \tilde{\mathbf A}\left[\frac 1 p \right] =
\left\{  \sum_{n=0}^N a_n\left(\frac {Y^{ea}} p \right)^n + \sum_{n
\in \N} b_n\left(\frac p {Y^{eb}} \right)^n ; a_n, b_n \in
\mathbf{\tilde A}^+, N \in \N \right\}.$$
\item The topology of $\mathcal G_{[b,a]}$ is weaker than the $p$-adic topology.
\item $\mathcal G_{[b,a]}$ is Hausdorff and complete.
\item The ring $\mathcal G_{[b,a]}$ is local with maximal ideal
$$\maxim_{[b,a]} = \left\{\sum_{n >1}a_n\left(\frac {Y^{ea}} p \right)^n + b_n\left(\frac p {Y^{eb}} \right)^n ;
a_n, b_n \in \mathbf{\tilde A}^+ \right\} + W(\maximE)$$ and residue
$\overline k$.
\item Any element of $\maxim_{[b,a]}$ is topologically nilpotent.
\item Powers of the ideal
$$\maxim^1_{[b,a]} = \left\{\sum_{n >1}a_n\left(\frac {Y^{ea}} p \right)^n + b_n\left(\frac p {Y^{eb}} \right)^n ;
a_n, b_n \in \mathbf{\tilde A}^+ \right\} + Y^{e(a-b)}\mathbf{\tilde
A}^+$$ form a basis of neighborhoods of $0$ consisting in ideals of
$\mathcal G_{[b,a]}$.
\item The ring $\mathcal G_{[b,a[}$ is local with maximal ideal
$\maxim_{[b,a[}$ the $p$-adic completion of
$$ \bigcup_{\alpha>a} \maxim_{[b,\alpha]}$$
and with residue field $\overline k$.
\item Any element of $\maxim_{[b,a[}$ is topologically nilpotent.
\end{enumerate}
\end{lem}
\emph{Proof:} Let us introduce the notation
$$\mathcal G_{[b,a]}^{>N} = \left\{ \sum_{n > N} a_n \left(\frac{Y^{ae}} p\right)^n + \sum_{n > N} b_n
\left(\frac p {Y^{be}} \right)^n ;a_n, b_n \in \mathbf{\tilde A}^+
\right\}\subset \mathcal G_{[b,a]}.$$ Recall that a basis of
neighborhoods of zero in $\mathcal G_{[b,a]}$ is given by
$$ \left\{\mathcal G_{[b,a]}^{>N} + p^k\mathcal G_{[b,a]}\right\}_{N,k \in \N}.$$
This shows the first two points. The fact that $\mathcal G_{[b,a]}$
is Hausdorff follows from that $\mathbf{\tilde A}_{[s,r]}$ is (cf.
\cite{Berger_equadiff}).

The following shows that the topology on $\mathcal G_{[b,a]}$ is
metrizable, and one can immediately see from the form of
neighborhoods of zero that any series with a general term going to
$0$ converges. This shows that $\mathcal G_{[b,a]}$ is complete.

We will prove Points $4.$, $5.$ et $6.$ simultaneously: we first
show $\maxim_{[b,a]}$ is an ideal, then that any element of
$\maxim_{[b,a]}$ has a power belonging to $\maxim^1_{[b,a]}$ and we
make powers of $\maxim^1_{[b,a]}$ explicit, which allows to
conclude.

Let
$$x =  \sum_{n<0} a_n \left(\frac {Y^{ea}} p \right)^{-n} + \sum_{n\geq 0} a_n \left(\frac p {Y^{eb}} \right)^n$$
we say that $x$ is the element of $\mathcal G_{[b,a]}$ associated
with the sequence $(a_n)_{n\in \Z} \in \left(\mathbf{\tilde
A}^+\right)^\Z$. Let $y$ be the element associated with another
sequence $(b_n)_{n\in \Z}$, then write the product of two elements
$x$ and $y$ of $\mathcal G_{[b,a]}$:
$$xy = \sum_{n<0} c_n \left(\frac {Y^{ea}} p \right)^{-n} + \sum_{n\geq 0} c_n \left(\frac p {Y^{eb}} \right)^n $$
is associated with the sequence
\begin{equation}\label{equa_c_n}
c_n = \left\{ \begin{array}{ll}\displaystyle{
\sum_{k >0} Y^{e(a-b)k} ( a_{k+n}b_{-k} + a_{-k}b_{k+n}) + \sum_{k =0}^n a_{k}b_{n-k}} & \textrm{ if } n\geq 0\\
\displaystyle{ \sum_{k >0} Y^{e(a-b)k} ( a_{k}b_{n-k} +
a_{n-k}b_{k}) + \sum_{k =0}^{-n} a_{-k}b_{n+k}} & \textrm{ if }
n\leq 0.
\end{array}\right.
\end{equation}
This yields to
\begin{equation}\label{equa_c_0}
c_0 = \sum_{n \in \Z} Y^{e(a-b)|n|} a_{n}b_{-n}
\end{equation}
and shows that $\maxim_{[b,a]}$ is an ideal.

Suppose $x \in \maxim_{[b,a]}$. Because of the previous computation,
one can define for all $k\in \N$ a sequence $(c_{n,k})_{n\in\Z}$
such that $x^k$ is associated with $(c_{n,k})_{n\in\Z}$. The fact
that there exists a $k$ such that $x^k \in \maxim_{[b,a]}^1$ is
equivalent to that the rest $\overline c_{0,k} \in \mathbf{\tilde
E}^+$ of $c_{0,k}$ modulo $p$ has a valuation greater or equal to
$a-b$. But because of Equality \eqref{equa_c_0},
$$v_{\mathbf E}(\overline c_{0,k}) \geq \min(a-b, kv_{\mathbf E}(\overline{a}_0))$$
which shows $x^k \in \maxim_{[b,a]}^1$ for $k$ large enough.

Let us show now that $\maxim^k_{[b,a]}$ consists in elements
associated with sequences $(a_n)_{n \in \Z}$ such that
$$\forall n \in \Z, \ v_{\mathbf E}(\overline a_n) \geq g_{a,b}^k(n) $$
where
$$g_{a,b}^k(n) = \left\lfloor \frac {(k-|n|+1)_+} 2 \right\rfloor (a-b) = \left\{ \begin{array}{ll}
 \left\lfloor \frac {k-|n|+1} 2 \right\rfloor (a-b) & \textrm{ if } |n|\leq k\\
 0 & \textrm{ otherwise} \end{array}\right.$$
satisfying the induction relation
\begin{equation}\label{rec1}
g_{a,b}^{k+1}(n) = \left\{ \begin{array}{ll}
 g_{a,b}^k(n-1) + a-b & \textrm{ if } -k-1 \leq n \leq 0 \\
 g_{a,b}^k(n+1) + a-b & \textrm{ if } 0 \leq n \leq k+1 \\
 0 & \textrm{ otherwise} \end{array}\right.
\end{equation}
or equivalently
\begin{equation}\label{rec2}
g_{a,b}^{k+1}(n)=\left\{ \begin{array}{ll}
 g_{a,b}^k(n+1) & \textrm{ if } n < 0 \\
 g_{a,b}^k(n-1) & \textrm{ if } n > 0
  \end{array}\right..
\end{equation}
Remark also that $g_{a,b}^k$ is even and decreasing on $\N$.

Let then $x \in \maxim^k_{[b,a]}$ be associated with a sequence
$(a_n)_{n \in \Z}$ satisfying the previous induction relation, let
$y\in \maxim^k_{[b,a]}$ be associated with $(b_n)_{n \in \Z}$ and
$xy\in \maxim^{k+1}_{[b,a]}$ be associated with $(c_n)_{n \in \Z}$
which we compute as before.

Equations \eqref{equa_c_n} show the relation for $n \geq 0$ (case
$n<0$ provides the same computation):
$$v_{\mathbf E}(\overline c_n) \geq \inf \left\{ \begin{array}{ll} (a-b)r + g_{a,b}^k(n+r), & \textrm{for } r>0,\\
(a-b)r + g_{a,b}^k(-r), & \textrm{for } r>0, \\
g_{a,b}^k(r), & \textrm{for } 0\leq r < n, \\
g_{a,b}^k(n)+a-b &
\end{array} \right\}$$
which gives because $g_{a,b}^k$ is even and decreasing
$$v_{\mathbf E}(\overline c_n) \geq \inf \left\{ \begin{array}{ll} (a-b)r + g_{a,b}^k(n+r), & \textrm{for } r>0,\\
g_{a,b}^k(n-1) & \\
g_{a,b}^k(n)+a-b &
\end{array} \right\}.$$
But
$$(a-b)r + g_{a,b}^k(n+r) = (a-b)\left(r +\left\lfloor \frac {(k-|n+r|+1)_+} 2 \right\rfloor \right)$$
is strictly increasing in $r$ and
$$(a-b) + g_{a,b}^k(n+1) \geq g_{a,b}^{k+1}(n)$$
because of \eqref{rec1}. Likewise,
$$ g_{a,b}^k(n)+a-b \geq g_{a,b}^{k+1}(n-1) \geq g_{a,b}^{k+1}(n)$$
and finally, according to \eqref{rec2},
$$g_{a,b}^k(n-1) = g_{a,b}^{k+1}(n).$$
The minimum is then equal to $g_{a,b}^{k+1}(n)$, which lets us
conclude on the description of $\maxim^k_{[b,a]}$.

Point $6.$ follows from this description, and proves $5.$ and $4.$

At last, $7.$ is a consequence of $8.$, which is left to show.
Remark that any element $x\in \maxim_{[b,a[}$ can be written as
$$x = x_0+px_1 \ ; \ x_0 \in \maxim_{[b,\alpha]}, x_1 \in \mathcal G_{[b,a[}$$
for some $\alpha>a$. We have to show
$$p^kx_0^nx_1^k \mathop{\longrightarrow}\limits_{k,n \rightarrow +\infty} 0.$$
When $k$ goes to infinity, it is clear. For the case where $n$ goes
to infinity, remark that the convergence of $x_0^n$ to $0$ in
$\mathcal G_{[b,\alpha]}$ implies for $n$ large enough that $x_0^n$
belongs to $p^N\mathcal G_{[b,\alpha']} + \mathcal
G_{[b,\infty[}^{>N}$ for $\alpha > \alpha' > a$, from which we
conclude. \Eproof

\paragraph{Remark} The preceding lemma makes $ \mathcal G_{[b,a)}$ into a complete valuation ring with the valuation given by
$$v_{[b,a)}(x) = \lim_{n \to \infty} \frac {k_n} n$$
where
$$k_n = \sup\{k \in \N, x^n\in \maxim^k_{[b,a)}\}.$$

The following lemma provides a link between algebras $\mathcal
G_{[b,a]}$ and Fontaine's rings.

\begin{lem}\label{alg G}
\begin{enumerate}
\item $\mathcal G_{[0,p]}$  injects continuously in $A_{crys}$.
\item Frobenius $\varphi$ of $A_{crys}$ and $\varphi_{\mathcal G}$ coincide on $\mathcal G_{[0,p]}$.
\item Any non zero element of $\mathcal G_{Y, [0,p]}$ is invertible
in $\mathcal G_{Y, [p-1,p-1[}\otimes_{\Z_p}\Q_p$.
\item The series defining $t/X$ converges in $\mathcal G_{[0,p]}$ where it is invertible.
\end{enumerate}
\end{lem}
\emph{Proof:} The first point consists in showing that $\frac
{Y^{pen}} {p^n} \in A_{crys}$ for all $n$ and converges to $0$. Let
$E_{\pi}$ be an Eisenstein polynomial for $\pi$, it is of degree $e$
and $E_\pi(Y)$ generates $W^1(\mathbf{\tilde E}^+)$ so that
$A_{crys}$ is the $p$-adic completion of $\mathbf{\tilde A}^+[\frac
{E_\pi(Y)^n} {n!}]$ and it is obvious that $\frac {Y^{pen}} {p^n}$
belongs to this ring and $p$-adically converges to $0$.

The second point is an immediate consequence of the first one.

Now let $x \in \mathcal G_{Y, [0,p]}$, then there exists a sequence
$(a_n)_{n\in \N} \in \left(\mathbf{\tilde A}^+\left[\frac 1 p
\right]\right)^\N$ such that
$$x = \sum_{n\in \N} a_n Y^n$$
with
$$\forall n\in \N \ ; \ epv_p(a_n)+n \geq 0.$$
Then,
$$\forall n\in \N \ ; \ e(p-1)v_p(a_n)+n \geq \frac n p$$
and for all $x$ non zero, $e(p-1)v_p(a_n)+n$ goes to $+\infty$ when
$n\rightarrow + \infty$, it reaches its minimum $K$ a finite number
of times and we fix $n_0$ the greater integer with $K =
e(p-1)v_p(a_{n_0})+n_0$, so that
\begin{eqnarray}
e(p-1)v_p(a_n/a_{n_0})+ n-n_0 \geq 0 & \textrm{ if } n \leq n_0 \label{equa_val_1}\\
e(p-1)v_p(a_n/a_{n_0})+ n-n_0 > 0 & \textrm{ if } n > n_0
\label{equa_val_2}
\end{eqnarray}
and
$$\forall n > n_0 \ ; \ e(p-1)v_p(a_n/a_{n_0})+ n-n_0 > \frac n p - K$$
hence it comes
$$ \mathop{\lim \inf}\limits_{n \rightarrow \infty} \frac {e(p-1)v_p(a_n/a_{n_0})+ n-n_0} {n-n_0} \geq \frac 1 p,$$
which,  combined with \eqref{equa_val_2}, shows the existence of
some $0< \lambda < 1$ such that
$$e(p-1)v_p(a_n/a_{n_0})+ n-n_0 \geq  \lambda(n-n_0) $$
hence
$$e\frac {p-1} {1-\lambda} v_p(a_n/a_{n_0})+ n-n_0 \geq 0.$$
This shows that for $a = \frac {p-1} {1-\lambda} > p-1$,
$$ \sum_{n> n_0} \frac {a_n}{a_{n_0}} Y^{n-n_0} \in \maxim_{[0,a]}.$$
Inequality \eqref{equa_val_1} shows furthermore
$$\sum_{n=0}^{n_0-1} \frac {a_n}{a_{n_0}} Y^{n-n_0}   \in \maxim_{[p-1,\infty[}$$
and finally
$$\sum_{n \neq n_0} \frac {a_n}{a_{n_0}} Y^{n-n_0} \in \maxim_{[p-1,a]}.$$
Then,
$$x =  a_{n_0} Y^{n_0}(1+\epsilon) \ ; \ \epsilon \in \maxim_{[p-1,a]}$$
is invertible in $\mathcal G_{Y, [p-1,a]}\otimes_{\Z_p}\Q_p \subset
\mathcal G_{Y, [p-1,p-1[}\otimes_{\Z_p}\Q_p$.

Remark finally
$$X = [\varepsilon -1] + pv = Y^{ep/(p-1)}u+pv \ ; \ u,v \in \mathbf{\tilde A}^+$$
so that
$$X^{p-1} = Y^{ep}u'+pv' \ ; \ u',v' \in \mathbf{\tilde A}^+$$
from which we deduce that for $s$ prime to $p$,
\begin{eqnarray*}
\frac {X^{p^rs-1}}{p^rs} & = & \frac {X^{p^r(s-1)}}{s} \frac {X^{p^r-1}}{p^r} \\
& =& X^{p^r(s-1)} \sum_{k=0}^{p^r-1} \frac {Y^{pek}}p^{\frac {p^r-1}{p-1}-k} {p^r}u_k\\
& =& X^{p^r(s-1)} \sum_{k=0}^{p^r-1} \frac {Y^{pek}}p^{\frac
{p^r-1}{p-1}-r} {p^k}u_k
\end{eqnarray*}
where $u_k \in \tilde{\mathbf A}^+$. But
$$\frac {p^r-1}{p-1} \geq r$$
so that for all $n \geq1$,
$$X^{n-1}/n \in \mathcal G_{[0,p]}$$
and $\frac {p^r-1}{p-1}-r$ goes to $+\infty$ with $r \to \infty$,
which shows that $X^{n-1}/n$ converges $p$-adically to $0$ in
$\mathcal G_{[0,p]}$ and completes the proof of the lemma.\Eproof

\subsection{The Hilbert symbol of a formal group}

\subsubsection{The pairing associated with the Hilbert symbol}

In this paragraph we express the Hilbert symbol of $F$ in terms
of the Herr complex attached to $F[p^M]$.\\
Let us recall that the Hilbert symbol of a formal group is defined
as the pairing:
$$\left. \begin{array}{rclcl}
K^* & \times & F(\maxim_K) & \rightarrow & F[p^M]\\
(\alpha & , & \beta) & \mapsto & (\alpha, \beta)_{F,M} = r(\alpha)
(\beta_1) -_F \beta_1
\end{array} \right.$$ where $\beta_1 \in F(\maxim_{\C_p})$ satisfies $p^M
\id_F \beta_1 = \beta$ and $r : K^* \rightarrow G_K^{\textrm{ab}}$
is the reciprocity map of local class field theory.\\
In fact, we will be interested in the pairing
$$\left. \begin{array}{rclcl}
F(\maxim_K) & \times & G_K & \rightarrow & F[p^M]\\
(\beta & , & g) & \mapsto & (\beta, g]_{F,M} = g\beta_1 -_F \beta_1
\end{array} \right.$$ where $\beta_1 \in F(m_{\C_p})$ satisfies $p^M
\id_F \beta_1 = \beta$. Then
$$(\beta, r(\alpha)]_{F,M} = (\alpha, \beta)_{F,M}.$$
Put
$$\mathcal R(F) = \{(x_i)_{i\geq 0} \in F(\maxim_{\C_p})
\textrm{ such that } x_0 \in F(\maxim_K) \textrm{ and }
(p\id_F)x_{i+1} = x_i \ \forall i \geq 0 \}$$ then the Hilbert
symbol is a mod $p$ reduction of the pairing
$$\left. \begin{array}{rclcl}
\mathcal R(F) & \times & G_K & \rightarrow & T(F)\\
(x & , & g) & \mapsto & (x, g]_{\mathcal R(F)} = (gx_i -_F x_i)_i
\end{array} \right.$$
with  $((x, g]_{\mathcal R(F)})_M = (x_0,g]_{F,M}$
for any $x = (x_i) \in \mathcal R(F)$.\\
We can see this pairing as the connection map
$$ F(\maxim_{K}) \rightarrow H^1(K, T(F))$$
in the long exact sequence associated with the short exact one:
$$ 0 \rightarrow T(F) \rightarrow   \lim_{\leftarrow} F(\maxim_{\C_p}) \rightarrow F(\maxim_{\C_p}) \rightarrow 0$$
where the transition maps in the inverse limit are $p\id_F$ and the
last map is the projection on the first component. The ring
$\mathcal R(F)$ is then the preimage of $F(\maxim_K)$ by surjection
$\displaystyle{\lim_{\leftarrow} F(\maxim_{\C_p})
\rightarrow F(\maxim_{\C_p})}$.\\
Let now $x \in F(\maximE)$ be such that $\theta([x])\in
F(\maxim_K)$. Then for all $g \in G_K$,
$$(g-1)\delta(x)\in F(W^1(\maximE))^{(\mathcal A-p)\circ l_{\mathcal A}=0} \simeq T(F)$$
where $\delta$ is the map defined at the end of Paragraph
\ref{RappelsF}. The following diagram is commutative

$$\xymatrix{
F(W(\maximE))_K^{(\mathcal A-p)\circ l_{\mathcal A}=0}\times G_K
\ar[r] &
F(W^1(\maximE))^{(\mathcal A-p)\circ l_{\mathcal A}=0} & \\
F(\maximE)_K \times G_K  \ar^{\delta \times \id}[u] \ar_{\iota \times \id}[d] \\
\mathcal R(F) \times G_K \ar[r] & T(F) \ar_{j}[uu]  \\
  }$$
where $\iota(x) = (\theta \circ \delta (p^{-s}\id_F(x)))_s =
(\theta([p^{-s}\id_F(x))])_s$ and $F(\maximE)_K$ (respectively
$F(W(\maximE))_K$) stands for the set of $x\in F(\maximE)$
(respectively $F(W(\maximE))$) with $\theta([x])\in K$ (respectively
$\theta(x)\in K$) and where the first pairing is simply
$$(u,g) \mapsto (g- 1)u.$$

\bigskip

Fix now $\alpha \in F(\maxim_K)$ and a lift $\xi$ of $\alpha$ in
$F(\maximE)$ which then satisfies
$$\theta([\xi])=\alpha.$$
We get the equality
$$j((\iota(\xi), g]_{\mathcal R(F)}) = (g-1)\delta(\xi) $$
for all $g \in G_K$.\\
Choose now $\beta \in F(YW[[Y]])$ such that
$$\theta(\beta) = \alpha =\theta([\xi]).$$
Then for all $h \in G_L$,
$$(h-1)(\delta(\xi)-_F\beta) =j((\iota(\xi), h]_{\mathcal R(F)}).$$
Moreover, $\delta(\xi)-_F\beta \in F(W^1(\maximE))$ thus
$$l_{\mathcal A}(\delta(\xi)-_F\beta) \in (\Fil^1A_{crys})^d$$
and
$$m_{\mathcal A}(\delta(\xi)-_F\beta) =\sum_{u \geq 1} F'_u \frac {\varphi^u
(l_{\mathcal A}(\delta(\xi)-_F\beta))} p$$ converges in
$A_{crys}^{h-n}$. Put now
$$\Lambda = \mathcal
V^{-1}\begin{pmatrix} l_{\mathcal A}(\delta(\xi)-_F\beta) \\
m_{\mathcal A}(\delta(\xi)-_F\beta) \end{pmatrix} \in A_{crys}^h.$$
These are the coordinates of an element $\lambda$ in $D_{crys}
(T(F)) \otimes A_{crys}$ in the basis $(o^1, \dots, o^h)$. And, for
all $h \in G_L$,
\begin{equation}
(h-1)\lambda = (\iota(\xi), h]_{\mathcal R(F)} + \left((h-1)\mathcal
V^{-1}\right)\begin{pmatrix} l_{\mathcal A}(\delta(\xi)-_F\beta) \\
m_{\mathcal A}(\delta(\xi)-_F\beta) \end{pmatrix}. \label{eqnlambda}
\end{equation}
\subsubsection{The approximated period matrix }

Let us now explicitly compute the Hilbert symbol of $F$, \emph{i.e.}
the image of $\iota(\xi)$ in $H^1 (K, F[p^M])$ which coincides with
the one of $\alpha$. For that, we have to give a triple in the first
homology group of the Herr complex of $F[p^M]$ corresponding to a
cocycle representing the image of $\iota(\xi)$. Recall that if such
a triple is written as $(x,y,z)$, then the associated cocycle is
$$g \mapsto (g-1)(-b) + \gamma^n \frac {\tau^m-1}{\tau - 1}z + \frac{\gamma^n - 1 }{\gamma - 1}y$$
where $g_{|_{\Gamma}} = \gamma^n\tau^m$ and $b \in F[p^M]\otimes
\tilde{\mathbf A}$ is a solution of
$$(\varphi - 1) b = x.$$
In particular, the image of $h \in G_L$ through this cocycle is
$(h-1)(-b)$. Let us start with finding $b \in T(F)\otimes
\mathbf{\tilde A}$ such that
$$\forall h \in G_L, \ \ (h-1) b \equiv -(\iota(\xi),h]_{\mathcal R(F)} \mod p^M.$$
Equality  \eqref{eqnlambda} incites to build $b$ as an approximation
of $-\lambda$. In fact, we will build $x$ by approximating
$(\varphi-1)(-\lambda)$, whose coordinates in the basis $(o^1,
\dots, o^h)$ are
$$\mathcal V^{-1}\begin{pmatrix} (\frac {\mathcal A} p -1)\circ l_{\mathcal A}(\beta) \\ 0
\end{pmatrix}.$$
Indeed, Lemma \ref{frob_form} shows that the action of Frobenius
$\varphi$ is written in the basis $(o^1, \dots, o^h)\mathcal
V^{-1}$:
$$\begin{pmatrix} \frac {\mathcal A} p & 0 \\ 0 & I_{h-d} \end{pmatrix}.$$
Because $(o^1 = (o_n^1)_n, \dots, o^h)$ is the fixed basis of
$T(F)$, $(o_M^1, \dots, o_M^h)$ is a basis of $F[p^M]$ and we
further fix $\hat o_M^1, \dots,\hat o_M^h$ elements in $F(YW[[Y]])$
such that for all $i$,
$$\theta(\hat o_M^i) = \hat o_M^i(\pi) = o_M^i.$$
Define then the matrix
$$\mathcal V_Y = \begin{pmatrix}
p^Ml_{\mathcal A} (\hat o_M^1) & \dots & p^Ml_{\mathcal A}(\hat o_M^h) \\
p^Mm_{\mathcal A} (\hat o_M^1) & \dots & p^Mm_{\mathcal A}(\hat
o_M^h)
\end{pmatrix} $$
whose coefficients belong to $A_{crys}$, and more precisely to
$W[[Y]]\left[\left[\frac{Y^{pe}} p\right]\right] = \mathcal
G_{Y,[0,p]}$. From Lemma \ref{alg G}, $\mathcal V_Y$ is invertible
in $\mathcal G_{Y, [p-1,p-1[}\otimes \Q_p$.
\begin{lem}\label{frob_inv}
\begin{enumerate}
\item $X\mathcal V_Y^{-1}$ has coefficients in
$\mathcal G_{[0,p]} + \frac {p^{M}} {Y^{e/(p-1)}} \maxim_{[ \frac 1
{p-1},\infty[}\subset \mathcal G_{[\frac 1 {p-1},p]}$ and thus
$$\varphi_{\mathcal G}(X\mathcal V_Y^{-1})\in \mathcal G_{[p/(p-1),p]}.$$
\item The matrix $\mathcal V_Y^{-1}$ has coefficients in $\frac 1 {Y^{ep/(p-1)}} \mathcal G_{[1,p]}$, then in
$\frac 1 {Y^{\lceil ep/(p-1) \rceil }} \mathcal G_{Y,[1,p]}$ and
$$ \mathcal V_Y^{-1} \equiv \mathcal V^{-1} \mod \frac {p^{M}}  {Y^{e(p+1)/(p-1)}} \maxim_{[1,p]}.$$
\item The principal part $\mathcal V_Y^{(-1)}$ of $\mathcal V_Y^{-1}$ has
$p$-entire coefficients and its derivative $\frac d {dY}\mathcal
V_Y^{(-1)}$ has coefficients in $p^M \mathbf{\tilde A}$.
\item The matrix $X\mathcal V_Y^{(-1)}$ has coefficients in $\mathbf{\tilde A}^+ +p^M\mathbf{\tilde A}$.
\end{enumerate}
\end{lem}
\emph{Proof:} We use the strategy of Paragraph 3.4. in
\cite{abr_hilb_forml}. Let us recall that Abrashkin there showed
\begin{eqnarray*}
p^Ml_{\mathcal A} (\hat o_M^i) & \in &
\left(E_\pi(Y)YW[[Y]] + \frac {E_\pi(Y)^p} p W[[Y]]\left[\left[\frac{Y^{ep}} p\right]\right]\right)^n\\
p^Mm_{\mathcal A} (\hat o_M^i) & \in & \left(YW[[Y]] + \frac
{Y^{ep}} p W[[Y]]\left[\left[\frac{Y^{ep}}
p\right]\right]\right)^{h-n}
\end{eqnarray*}
and
\begin{eqnarray*}
\mathbf l(o^i)-p^Ml_{\mathcal A} (\hat o_M^i) & \in &
p^M\left(E_\pi(Y)W(\maximE) + \frac {E_\pi(Y)^p} p \mathbf{\tilde A}^+ \left[\left[\frac{Y^{ep}} p\right]\right]\right)^n\\
\mathbf m(o^i) - p^Mm_{\mathcal A} (\hat o_M^i) & \in &
p^M\left(W(\maximE) + \frac {Y^{ep}} p \mathbf{\tilde A}^+
\left[\left[\frac{Y^{ep}} p\right]\right]\right)^{h-n}.
\end{eqnarray*}
Let $\mathcal V^D$ be the matrix of the group dual to $F$. It
satisfies the relation:
$$ {}^t\mathcal V^D \  \mathcal V = tI_h.$$
And one can then write
$$ {}^t\mathcal V^D \ \mathcal V_Y \equiv  tI_h \mod p^M\left(E_\pi(Y)W(\maximE)+ \frac {E_\pi(Y)^p} p \mathbf{\tilde A}^+
\left[\left[\frac{Y^{ep}} p\right]\right] \right)$$ in particular,
$$ {}^t\mathcal V^D \ \mathcal V_Y \equiv  tI_h  \mod p^M\left(Y^eW(\maximE)+ \frac {Y^{ep}} p \mathbf{\tilde A}^+
\left[\left[\frac{Y^{ep}} p\right]\right] \right).$$ Remark, because
of Lemma \ref{alg G}, that the element $t/X$ converges in $\mathcal
G_{[0,p]}^*$, and
$$X = \omega [\varepsilon^{1/p}-1] = E_\pi(Y)Y^{e/(p-1)}v \ ; \ v \in \mathcal G_{[\frac 1 {p-1},\infty[}^*,$$
so that
$$ {}^t\mathcal V^D \ \mathcal V_Y =  t(I_h - p^Mu)$$
with
\begin{eqnarray*}
u \in  \frac {E_\pi(Y)} t W(\maximE)+ \frac {Y^{ep}} {pt}\mathcal
G_{[0,p]} & \subset & \frac 1 {Y^{e/(p-1)}} \maxim_{[\frac 1
{p-1},p]}+ \frac
{Y^{e\frac{p^2-2p}{p-1}}} p \mathcal G_{[\frac 1 {p-1},p]} \\
& \subset & \frac 1 {Y^{e/(p-1)}} \maxim_{[\frac 1 {p-1},p]} \subset
\frac 1 p \maxim_{[\frac 1 {p-1},p]}
\end{eqnarray*} thus
$$ \mathcal V_Y^{-1} = \frac 1 t \left(\sum_{n \in \N} p^{Mn}u^n\right) {}^t\mathcal V^D  \in
\frac 1 t \mathcal G_{[\frac 1 {p-1},p]}$$ and then we deduce the
first point ; and even
$$ \mathcal V_Y^{-1} \equiv   \mathcal V^{-1} \mod \frac {p^{M}}  {tY^{e/(p-1)}} \maxim_{[\frac 1 {p-1},p]}$$
or
$$ \mathcal V_Y^{-1} \in \frac 1 t \mathcal G_{[0,p]} + \frac {p^{M}} {Y^{e/(p-1)}} \maxim_{[\frac 1 {p-1},p]}.$$

Recall
$$t = E_\pi(Y)\varphi^{-1}(X)u' \ ; \ u' \in \mathcal G_{[0,p]}^*$$
and remark that because $E_\pi$ is an Eisenstein polynomial,
$E_\pi(Y)$ and $Y^e$ are associated in $\mathcal G_{[1,\infty[}$ ;
finally, with the above computation, we deduce that $t$ and
$Y^{ep/(p-1)}$ are associated in $\mathcal G_{[1,p]}$. Then
$$ \mathcal V_Y^{-1}  \in \frac 1 {Y^{ep/(p-1)}} \mathcal G_{[1,p]} + \frac {p^{M}}  {Y^{e(p+1)/(p-1)}} \maxim_{[1,p]}
\subset \frac 1 {Y^{ep/(p-1)}} \mathcal G_{[1,p]}.$$ So,  $Y^{\lceil
ep/(p-1)\rceil} \mathcal V_Y^{-1}$ has coefficients in $\mathcal
G_{Y,[p-1,p-1[}\left[\frac 1 p \right] \bigcap \mathcal G_{[1,p]} =
\mathcal G_{Y,[1,p]}$
because of Lemma \ref{alg_G_intersections}.\\
Let us further deduce that $\mathcal V_Y^{(-1)}$ has $p$-entire
coefficients. It is to show that any element
$$x = \sum_{n\in \Z}a_n Y^n \in \frac 1 {Y^{\lceil ep/(p-1)\rceil}} \mathcal G_{Y,[1,p]}$$ satisfies $a_n \in W$
for all $n\leq 0$. But that means that
$$Y^{\lceil ep/(p-1) \rceil}x = \sum_{n\in \Z}a_n Y^{n+\lceil ep/(p-1) \rceil} \in \mathcal G_{Y,[1,p]}$$
and thus if $v_p(a_n)\leq 1$, the following inequality holds
$$ n+ep/(p-1) \geq ep $$
hence
$$ n\geq \frac{(p-2)ep} {p-1} >0$$
and $\mathcal V_Y^{(-1)}$ has $p$-entire coefficients.

For the third point, let us recall the argument of Lemma 4.5.4 in
\cite{abr_hilb_forml}. Write
$$\frac d {dY}\mathcal V_Y^{(-1)} =
-\mathcal V_Y^{(-1)}\left(\frac d {dY}\mathcal V_Y \right)\mathcal
V_Y^{(-1)} $$ and because differentials of $l_{\mathcal A}$ and
$m_{\mathcal A}$ have coefficients in $W$, one gets
$$\frac d {dY}\mathcal V_Y  \in p^M M_h(W[[Y]])$$
so that
$$\frac d {dY}\mathcal V_Y^{(-1)} \in \mathcal G_{Y, [p-1,p-1[}\otimes_{\Z_p}\Q_p \bigcap
\frac 1 {Y^{2ep/(p-1)}} \mathcal G_{[1,p]}$$ and the same argument
as above permits to conclude (we get then the inequality $ n\geq
\frac{(p-3)ep} {p-1} \geq 0$).

Finally, the proof of Point $4.$ is the same as the one of
Proposition 3.7, Point \emph{d)} in \cite{abr_hilb_forml}. Let us
write it in the following way: we know on the one hand that
$\mathcal V_Y^{(-1)}$ and then also $X\mathcal V_Y^{(-1)}$ has
$p$-entire coefficients, so they have coefficients in $\mathcal
G_{[p-1,\infty[}\left[\frac 1 Y\right]$ and that $\mathcal U =
X(\mathcal V_Y^{-1} - \mathcal V_Y^{(-1)})$ has coefficients in
$\mathcal G_{[0,p-1[}\left[\frac 1 p\right]$. On the other hand,
Lemma \ref{frob_inv} tells
$$X\mathcal V_Y^{-1}\in M_h\left(\mathcal G_{[0,p]} + p^{M-1}\mathcal G_{[1/(p-1),\infty[}\right).$$
Remark
$$\mathcal G_{[1/(p-1),\infty[} = \mathbf{\tilde
A}^+\left[\left[\frac p{Y^{e/(p-1)}}\right]\right] = \mathbf{\tilde
A}^+ + \frac p{Y^{e/(p-1)}}\mathcal G_{[1/(p-1),\infty[}$$ Thus we
can write
$$X\mathcal V_Y^{-1} = M_1 + p^MM_2$$
with $M_1$ having coefficients in $\mathcal G_{[0,p]}$ and $M_2$ in
$\frac 1 {Y^{e/(p-1)}}\mathcal G_{[1/(p-1),\infty[} \subset
\mathbf{\tilde A}$. So
$$X\mathcal V_Y^{(-1)} - p^MM_2 = M_1 - \mathcal U$$
has coefficients in $\mathcal G_{[p-1,\infty[}\left[\frac 1 Y\right]
\bigcap \mathcal G_{[0,p-1[} \left[\frac 1 p\right]  =
\mathbf{\tilde A}^+$, as desired.\Eproof

Remark that if $x \in F(W(\maximE))$ it can be written as
$$ x = [x_0] +_F u$$
with $u \in F(pW(\maximE))$, and thus
\begin{eqnarray*}
\left(\frac {\mathcal A} p -1\right)\circ l_{\mathcal A}(x) & = &
 \left(\frac {\mathcal A} p -1\right)\circ l_{\mathcal A}([x_0]) + \left(\frac {\mathcal A} p -1\right)\circ l_{\mathcal A}(u) \\
& = & [x_0] + \left(\frac {\mathcal A} p -1\right)\circ l_{\mathcal
A}(u) \in W(\maximE)^d
\end{eqnarray*}
since $l_{\mathcal A}(u) \in pW(\maximE)^d$. In particular
$$\left(\frac {\mathcal A} p -1\right)\circ l_{\mathcal A}(\beta) \in W(\maximE)^d,$$
so that
$$\mathcal V_Y^{(-1)}\begin{pmatrix} \left(\frac {\mathcal A} p -1\right)\circ
l_{\mathcal A}(\beta) \\ 0 \end{pmatrix} \in \mathbf{\tilde A}^h.$$

\subsubsection{An explicit computation of the Hilbert symbol}

We come now to the proposition that explicitly gives the desired
triple. The basic ingredient is Proposition 3.8 of
\cite{abr_hilb_forml} which provides the $x$ coordinate of the
triple and allows us to prove that $y$ is zero. However, in order to
get $z$, we have to carry Abrashkin's computations to the higher
order. Indeed, we already know that $z$ belongs to $W(\maximE)$, but
we need to specify its value modulo $XW(\maximE)$.

Let us recall the results we are going to use

\begin{prop}\label{hilbert_abr}
Let $U$ be the principal part of $\mathcal V_Y^{(-1)}\begin{pmatrix}
\left(\frac {\mathcal A} p -1\right)\circ l_{\mathcal A}(\beta) \\ 0
\end{pmatrix}$ and $\hat x = (o^1, \dots, o^h)U$. Then
\begin{enumerate}
\item $U \in (W[[\frac 1 Y]] \cap \mathbf{\tilde A})^h$
\item Let $\hat b \in T(F)\otimes \mathbf{\tilde A}$ be a solution of $(\varphi-1)\hat b
=\hat x$ then for any $g \in G_K$,
$$(g-1)\hat b \equiv (\beta(\pi),g]_{F,M} \mod p^M\mathbf{\tilde A} + W(\maximE).$$
\end{enumerate}
\end{prop}

\bigskip

\emph{Proof:}\\
The first point can be shown like Point 3. of Lemma \ref{frob_inv} above.\\
The second point can be viewed as a reformulation of Proposition 3.8
of \cite{abr_hilb_forml}.
Let us give another proof.\\
Let us recall from Lemma \ref{frob_inv} that
$$ \mathcal V_Y^{-1} \equiv \mathcal V^{-1} \mod \frac {p^{M}}  {Y^{e(p+1)/(p-1)}} \maxim_{[1,p]}$$
so that there is  $\delta \in
\frac{p^{M}}{Y^{e(p+1)/(p-1)}}\maxim_{[1,p]}$ such that
$$ X\mathcal V_Y^{-1} = X\mathcal V^{-1} + X\delta.$$
Write $\delta =\delta_1 + \delta_2$ with $\delta_1 \in
p^{M-1}Y^{e(p^2-2p-1)/(p-1)}\mathcal G_{[0,p]}$ and
$\delta_2 \in \frac {p^M}{Y^{e(p+1)/(p-1)}}\maxim_{[1,\infty[}$.\\
Let us recall that we write $\mathcal V_Y^{-1} = \mathcal V_Y^{(-1)}
+ \mathcal U$ so that
$$X\mathcal V_Y^{(-1)} - \delta_2 = X\mathcal V^{-1} + X\delta_1 - X \mathcal U$$
has coefficients in $\mathcal G_{[p-1,\infty[}\left[\frac 1 Y\right]  \bigcap \mathcal G_{[0,p-1[}\left[\frac 1 p\right] = \mathbf {\tilde A}^+$.\\
Then, if $\mathcal B$ is a matrix with coefficients in $\mathbf
{\tilde A}$ such that
\begin{equation}\label{approx_B}
(\varphi -1)\mathcal B = \left(\mathcal V_Y^{(-1)} -
\delta_2\right)\begin{pmatrix} \left(\frac {\mathcal A} p
-1\right)\circ l_{\mathcal A}(\beta) \\ 0
\end{pmatrix},
\end{equation}
write as in Paragraph \ref{KummersMap},
$$(\varphi - \omega)(X_1\mathcal B) = \left(X\mathcal V_Y^{(-1)} - X\delta_2\right)\begin{pmatrix}
\left(\frac {\mathcal A} p -1\right)\circ l_{\mathcal A}(\beta) \\ 0
\end{pmatrix}$$
has coefficients in $\mathbf {\tilde A}^+ $ so that, by successive
approximations modulo $p^k$ and since $\mathbf {\tilde E}^+$ is
integrally closed, we get
$$\mathcal B \in \frac 1 {X_1} \mathbf {\tilde A}^+ \subset \Fil^0B_{crys}.$$
Still write
$$\Lambda = \mathcal V^{-1}\begin{pmatrix} l_{\mathcal A}(\delta(\xi)-_F\beta) \\
m_{\mathcal A}(\delta(\xi)-_F\beta) \end{pmatrix} \in
\left(\Fil^0A_{crys}\right)^h.$$ We compute
$$(\varphi-1)\left(\mathcal B-\Lambda\right) = \left(\delta_1 - \mathcal U\right)\begin{pmatrix}
\left(\frac {\mathcal A} p -1\right)\circ l_{\mathcal A}(\beta) \\ 0
\end{pmatrix}$$
Since $\delta_1' = \delta_1 \begin{pmatrix} \left(\frac {\mathcal A}
p -1\right)\circ l_{\mathcal A}(\beta) \\ 0 \end{pmatrix}$ has
coefficients in $Y\mathcal G_{[0,p]}$, the series $-\sum_{n\in
\N}\varphi^n(\delta_1')$ converges to an element $\Delta_1 \in
Y\mathcal G_{[0,p]}$ satisfying
$$(\varphi-1)(\Delta_1) = \delta_1'.$$
Likewise $\delta_2' = \mathcal U \begin{pmatrix} \left(\frac
{\mathcal A} p -1\right)\circ l_{\mathcal A}(\beta) \\ 0
\end{pmatrix}$ has coefficients in $YW[[Y]] + \frac {Y^{ep-\lceil
ep/(p-1) \rceil}} p \mathcal G_{Y,[0,p]}$ so that the series
$-\sum_{n\in \N}\varphi^n(\delta_2')$ converges to an element
$\Delta_2$ with coefficients in $YW[[Y]] + \frac {Y^{ep-\lceil
ep/(p-1) \rceil}} p \mathcal G_{Y,[0,p]}$ satisfying
$$(\varphi-1)(\Delta_2) = \delta_2'.$$
Finally,
$$(\varphi-1)\left(\mathcal B-\Lambda - \Delta_1 + \Delta_2\right) = 0$$
with $\mathcal B-\Lambda - \Delta_1 + \Delta_2$ having coefficients
in $\Fil^0 B_{crys}$. And the fact that
$$(\Fil^0B_{crys})_{\varphi = 1} = \Q_p$$
shows
$$\mathcal B-\Lambda - \Delta_1 + \Delta_2 \in \Q_p.$$
Then, for $g \in G_K$, $(g-1)\left(\mathcal B-\Lambda - \Delta_1 +
\Delta_2\right) = 0 $ so that
\begin{eqnarray*}
(g-1)\left(\mathcal B\right) & = & (g-1)\left(\Lambda + \Delta_1 - \Delta_2\right) \\
(g-1)\left(\mathcal B\right) & = & (g-1)\left(\mathcal V^{-1}\begin{pmatrix} l_{\mathcal A}(\delta(\xi)-_F\beta) \\
m_{\mathcal A}(\delta(\xi)-_F\beta) \end{pmatrix}\right) + (g-1)\left(\Delta_1 - \Delta_2\right) \\
(g-1)\left(\mathcal B\right) & = & \left((g-1)\mathcal V^{-1}\right)g\begin{pmatrix} l_{\mathcal A}(\delta(\xi)-_F\beta) \\
m_{\mathcal A}(\delta(\xi)-_F\beta) \end{pmatrix} \\ && +  \mathcal V^{-1} (g-1)\left(\begin{pmatrix} l_{\mathcal A}(\delta(\xi)-_F\beta)\\
m_{\mathcal A}(\delta(\xi)-_F\beta) \end{pmatrix}\right) +
(g-1)\left(\Delta_1 - \Delta_2\right)
\end{eqnarray*}
Now, we remark that
\begin{eqnarray*}
\mathbf l(o^i)- g \mathbf l(o^i) & \in &
p^M\left(E_\pi(Y)W(\maximE) + \frac {E_\pi(Y)^p} p \mathbf{\tilde A}^+ \left[\left[\frac{Y^{ep}} p\right]\right]\right)^n\\
\mathbf m(o^i) - g\mathbf m(o^i) & \in & p^M\left(W(\maximE) + \frac
{Y^{ep}} p \mathbf{\tilde A}^+ \left[\left[\frac{Y^{ep}}
p\right]\right]\right)^{h-n}
\end{eqnarray*}
so that $g\mathcal V$ enjoys the same approximation properties as
$\mathcal V_Y$, hence $(g-1)\mathcal V^{-1}$ has coefficients in
$\frac {p^{M}}  {Y^{e(p+1)/(p-1)}} \maxim_{[1,p]}$. Thus
coefficients of
$\left((g-1)\mathcal V^{-1}\right)g\begin{pmatrix} l_{\mathcal A}(\delta(\xi)-_F\beta) \\
m_{\mathcal A}(\delta(\xi)-_F\beta) \end{pmatrix} +
(g-1)\left(\Delta_1 - \Delta_2\right)$ lie in $\frac {p^{M}}
{Y^{e(p+1)/(p-1)}} \maxim_{[1,\infty[} + \frac Y p \mathcal
G_{[0,p]}$. Let us recall finally that
$$\mathcal V^{-1} (g-1)\left(\begin{pmatrix} l_{\mathcal A}(\delta(\xi)-_F\beta) \\
m_{\mathcal A}(\delta(\xi)-_F\beta) \end{pmatrix}\right)$$ are the
coordinates of $(\iota(\xi), g]_{\mathcal R(F)}$ in the basis
$(o_1,\dots, o^h)$ of $T(F)$, and thus it is congruent to the
coordinates of $(\beta(\pi),g]_{F,M}$ modulo $p^M\Z_p$. We get
$$(g-1)\left(\mathcal B\right) - (\iota(\xi), g]_{\mathcal R(F)} \in \frac {p^{M}}{Y^{e(p+1)/(p-1)}} \maxim_{[1,\infty[} + \frac Y p \mathcal G_{[0,p]}.$$
Recall now that coefficients of $\mathcal B$ lie in
$$\frac 1 {X_1}\mathbf{\tilde A}^+ \subset \frac 1 {Y^{e/(p-1)}}\mathcal G_{[\frac 1 {p-1},\infty[}.$$
Gathering all information, we deduce the existence of $u_1\in \frac
{1}  {Y^{e(p+1)/(p-1)}} \maxim_{[1,\infty[}$ and $u_2 \in  \frac Y p
\mathcal G_{[0,p]}$ such that
$$(g-1)\left(\mathcal B\right) - (\iota(\xi), g]_{\mathcal R(F)} -p^Mu_1 = u_2$$
has coefficients in
$$\frac 1 {Y^{e(p+1)/(p-1)}}\mathcal G_{[\frac 1 {p-1},\infty[} \bigcap  \frac Y p \mathcal G_{[0,p]} = Y\mathbf{\tilde A}^+ \subset W(\maximE)$$
so that $(g-1)\left(\mathcal B\right)$ is congruent to coordinates
of $(\beta(\pi),g]_{F,M}$ modulo $p^M\mathbf{\tilde A} +
W(\maximE)$. To finish the proof, just recall Equality
\eqref{approx_B}: there is $\delta_2 \in \frac
{p^M}{Y^{e(p+1)/(p-1)}}\maxim_{[1,\infty[} \subset p^M\mathbf{\tilde
A}$ such that
$$(\varphi -1)\mathcal B = \left(\mathcal V_Y^{(-1)} - \delta_2\right)\begin{pmatrix}
\left(\frac {\mathcal A} p -1\right)\circ l_{\mathcal A}(\beta) \\ 0
\end{pmatrix}$$
And surjectivity of $\varphi-1$ on $\mathbf{\tilde A}$ permits to
conclude.\Eproof
\paragraph{Remark}
It is possible to get rid of $A_{crys}$ in the proof by studying the
action of $(\varphi -1)$ on $\mathcal G_{[0,p]}\left[\frac 1
Y\right]$.

\par

We will use this result in the following specified form.

\begin{prop}\label{hilbert_formel} Let $\beta \in F(YW[[Y]])$ and $\alpha = \theta(\beta) = \beta(\pi) \in F(\maxim_K)$.
Put
$$x = (o^1, \dots, o^h)\mathcal V_Y^{(-1)}\begin{pmatrix} (\frac {\mathcal A} p -1)\circ l_{\mathcal
A}(\beta) \\ 0 \end{pmatrix} \in \tilde D_L(T(F))$$ then there
exists
$$z \in \tilde D_L(T(F)) \cap T(F) \otimes W(\maximE)$$
unique modulo $p^M$ such that the class of the triple $(x,0,z)$
corresponds to the image of $\alpha$ by the Kummer map $F(\maxim_K)
\rightarrow H^1(K, F[p^M])$.\\
Moreover, $z$ is congruent to
$$XY\mathcal V_Y^{(-1)}\displaystyle{\frac d {dY}}\begin{pmatrix} l_{\mathcal
A}(\beta) \\ m_{\mathcal A}(\beta) \end{pmatrix} \mod X
W(\maximE).$$
\end{prop}

\subsubsection{Proof of Proposition \ref{hilbert_formel}}

We use Proposition \ref{hilbert_abr}, and remark that
$$\hat x -x \in T(F)\otimes YW[[Y]] \subset (\varphi - 1)(T(F)\otimes
YW[[Y]]).$$ So, if $b\in T(F)\otimes \mathbf{\tilde A}$ satisfies
$(\varphi-1) b = x$, then for any $g \in G_K$,
$$(g-1) b \equiv (\alpha,g]_{F,M} \mod p^M\mathbf{\tilde A} + W(\maximE).$$
Thus for any $h \in G_L$,
$$(h-1)b \equiv (\alpha,h]_{F,M} \mod p^MT(F)$$
for $(h-1)b \in \ker(\varphi -1) = T(F)$.

 We deduce that there exist $y, z \in \tilde D_L(T(F))$ unique modulo $p^M$ such that
the class of the triple $(x,y,z)$ corresponds to the image of
$\alpha$ in $H^1(K, F[p^M])$ ; indeed let $(x_1, y_1, z_1)$ be such
a triple, and $b_1\in T(F)\otimes \mathbf{\tilde A}$ a solution of
$(\varphi - 1) b_1 = x_1$ then for all $h \in G_L$,
$$(h-1)(b_1-b) \equiv 0 \ \ \mod p^M,$$
$$\textrm{thus, } \ b_1-b \in \tilde D_L(F[p^M]),$$
which shows that the class of
$$(x,y_1+(\gamma-1)(b-b_1), z_1+(\tau-1)(b-b_1))$$
corresponds to the same class as $(x_1,y_1,z_1)$ and, if $x$ is fixed, this triple is unique.\\
Let us now determine $y$: let $\tilde \gamma$ lift $\gamma$ then
$$(\tilde \gamma-1)(-b) + y = (\alpha,\tilde \gamma]_{F,M}
\equiv (\tilde \gamma-1)(-b) \mod p^M\mathbf{\tilde A} +
W(\maximE)$$ hence, for $(\tilde \gamma-1)(-b)\in T(F)$,
$$y \in T(F) \otimes W(\maximE) \cap T(F) = \{0\}.$$
Likewise, let $\tilde \tau$ lift $\tau$ then
$$(\tilde \tau-1)(-b) + z = (\alpha,\tilde \tau]_{F,M} \equiv (\tilde
\tau-1)(-b) \mod p^M\mathbf{\tilde A} + W(\maximE)$$
hence $z \in T(F)\otimes W(\maximE)$.\\
Thus $z$ belongs to $T(F) \otimes W(\maximE)$ and satisfies
$$(\tau-1)x=(\varphi-1)z.$$
This uniquely determines $z$ since $\varphi-1$ is injective on $T(F)
\otimes W(\maximE)$. In order to specify $z$, we need the following
lemma:

\begin{lem}\label{lem_mat_periodes}
\begin{enumerate}
\item For all $U \in W[[Y]]$, the following congruence holds
$$ (\tau -1)\mathcal V_Y^{(-1)}U \equiv XY\mathcal V_Y^{(-1)}\frac
{dU}{dY} \mod XW(\maximE)+p^M \mathbf{\tilde A}.$$
\item There exists $u \in \maxim_{[ p/(p-1),p]}$ such that
$$ \varphi_{\mathcal G}(X\mathcal V_Y^{-1}) =
(\varphi_{\mathcal G}(X) \mathcal V_Y^{-1}\mathcal E^{-1}
+p^Mu)\begin{pmatrix} \frac 1 p I_d & 0\\ 0 & I_{h-d} \end{pmatrix}
$$
\end{enumerate}
\end{lem}
\emph{Proof of the lemma: } For Point $2.$, we first specify $(\tau
-1)\mathcal V_Y^{(-1)}$. Remark that if $f(Y)$ is a series in
$W\{\{Y\}\} \cap \mathbf{\tilde A}$,
$$(\tau -1)f(Y) = \sum_{n \geq 1} \frac{(XY)^n}{n!}f^{(n)}(Y)$$
Thus for $\mathcal V_Y^{(-1)}$:
$$(\tau -1)\mathcal V_Y^{(-1)} = XY\frac{d}{dY}\mathcal V_Y^{(-1)}+
\frac{(XY)^2}{2}\frac{d^2}{dY^2}\mathcal V_Y^{(-1)} + \sum_{n \geq
3} \frac{(XY)^n}{n!}\frac{d^n}{dY^n}\mathcal V_Y^{(-1)}$$ We then
have to estimate $\displaystyle{\frac{(XY)^n}{n!} \frac{d^n}{dY^n}
\mathcal V_Y^{(-1)}}.$ Lemma \ref{frob_inv} shows
 $$\frac d  {dY} \mathcal V_Y^{-1} = p^M \mathcal V_Y^{-1} \widetilde W \mathcal
 V_Y^{-1},$$
where $\widetilde W \in W[[Y]]$ and the principal part of $\mathcal
V_Y^{-1} \widetilde W  \mathcal  V_Y^{-1}$ is entire. Thus, on the
one hand
$$XY\frac{d}{dY}\mathcal V_Y^{(-1)}+ \frac{(XY)^2}{2} \frac{d^2}{dY^2}\mathcal V_Y^{(-1)}
 \in p^M \mathbf{\tilde A}$$
and on the other hand one can write
$$\frac{d^n}{dY^n}\mathcal V_Y^{-1} = \sum_{k=1}^n  p^{Mk} w_{n,k}$$
where the $w_{n,k}$ are sums of terms of the form
$$\mathcal V_Y^{-1}\widetilde W_{n,1} \mathcal V_Y^{-1}\widetilde W_{n,2} \dots \widetilde W_{n,k} \mathcal V_Y^{-1},$$
where the $\widetilde W_{n,i} \in W[[Y]]$ are derivatives of $\widetilde W$.\\
Recall that $\mathcal V_Y^{-1}$ has coefficients in
$\displaystyle{\frac 1 X\left( \mathcal G_{[0,p]} +
p^{M-1}\maxim_{[1/(p-1),p]}\right)},$ then
$$\mathcal V_Y^{-1}\widetilde W_{n,1} \mathcal
V_Y^{-1}\widetilde W_{n,2} \dots \widetilde W_{n,k} \mathcal
V_Y^{-1} \in M_h\left(\frac 1 {X^{k+1}} \mathcal G_{[0,p]} +
p^{M-1}\left(\frac 1 {X^{k+1}}\maxim_{[1/(p-1),p]}\right)\right).$$
Suppose $1<k<n-1$. Since
$$v_p(n!)\leq \lfloor n/(p-1) \rfloor = n',$$ there exists $u \in \Z_p$ such that
$$p^{Mk}\frac{(XY)^n}{n!} = Y^n  X^{k+2} u \frac {X^{n-k-2}}{p^{n'-Mk}}.$$
Since $p>2$ and $k>1$,
$$(n'-Mk)(p-1) \leq n-k-2 \textrm{ and }\frac {X^{n-k-2}}{p^{n'-Mk}} \in W\left[\left[\frac{X^{p-1}}{p}\right]\right]
\subset \mathcal G_{[0,p]}.$$ Thus,
$$p^{Mk}\frac{(XY)^{n-1}}{n!}w_{n,k} \in \maxim_{[0,p]} + p^{M-1}\maxim_{[1/(p-1),p]}.$$
Let now $k=1$, write
$$w_{n,1} = \mathcal V_Y^{-1}\frac{d^{n-1}}{dY^{n-1}}\widetilde W\mathcal V_Y^{-1} = (n-1)! \mathcal
V_Y^{-1} \widetilde W' \mathcal V_Y^{-1}$$ and
$$\frac{(XY)^{n-1}}{n!}p^Mw_{n,1} = \frac{(XY)^{n-1}}{n}\mathcal V_Y^{-1} \widetilde W' \mathcal V_Y^{-1} =
\frac{X^{n-p}}{n}Y^nX^{p-1}\mathcal V_Y^{-1} \widetilde W' \mathcal
V_Y^{-1}$$ has coefficients in  $\maxim_{[0,p]} +
p^{M-1}\maxim_{[1/(p-1),p]}$
as before.\\
Let $k=n$, one has obviously
$$w_{n,n} = n! \mathcal V_Y^{-1}\widetilde W_{n,1} \mathcal V_Y^{-1}\widetilde W_{n,2} \dots
\widetilde W_{n,n} \mathcal V_Y^{-1}$$ where for all $1 \leq i \leq
n$, $\widetilde W_{n,i}=\widetilde W$, so that
$$\frac{(XY)^{n}}{n!}p^{Mn}w_{n,n} \in p^{Mn}\frac 1 X\left(\mathcal G_{[0,p]} + p^{M-1}\maxim_{[1/(p-1),p]}\right).$$
Finally, for $k=n-1$, since $v_p(n!)\leq n/(p-1) \leq n-1$,
$$\frac{(XY)^{n}}{n!}p^{M(n-1)}w_{n,n-1} \in Y^{n}\left(\mathcal G_{[0,p]} + p^{M-1}\mathcal G_{[1/(p-1),p]}\right).$$
The same argument as $1.$ then shows that for all $n>2$,
$$\frac{(XY)^n}{n!} \frac{d^n}{dY^n} \mathcal V_Y^{(-1)} \in XW(\maximE) + p^M
\mathbf{\tilde A}$$ hence
$$(\tau-1)\mathcal V_Y^{(-1)} \in XW(\maximE) + p^M \mathbf{\tilde A}.$$
Point $2.$ then follows from the equality
$$(\tau -1)\mathcal V_Y^{(-1)}U = \left((\tau -1)\mathcal V_Y^{(-1)}\right)\tau U + \mathcal
V_Y^{(-1)}(\tau -1)U$$ and the congruence
$$(\tau -1)U \equiv XY \frac {dU}{dY}\mod XW(\maximE).$$

Now, let us carry on computations of Lemma \ref{frob_inv}:
$$X\mathcal V_Y^{-1} = \frac X t (I_h + p^{M-1} u_1) ^t\mathcal V^{D}$$
with $u_1 \in \maxim_{[\frac 1 {p-1},p]}$. And since $\mathcal V^D$
has coefficients in $\mathcal G_{[0,p]} \subset A_{crys}$ where
$\varphi_{\mathcal G}$ and $\varphi$ coincide, the following holds
in $\mathcal G_{[p/(p-1),p]}$:
\begin{eqnarray*}
\varphi_{\mathcal G}\left(X\mathcal V_Y^{-1}\right) & = &
\varphi_{\mathcal G}\left(\frac X
t\right)\left(I_h+p^M\varphi_{\mathcal
G}(v_1)\right)\varphi_{\mathcal
G}\left({}^t\mathcal V^D\right)\\
 & = & \varphi_{\mathcal G}\left(\frac X t\right)\left(I_h+p^M\varphi_{\mathcal G}(v_1)\right)p \ ^t\mathcal
V^D\mathcal{E}^{-1}\begin{pmatrix} \frac 1 p I_d & 0\\ 0 & I_{h-d}
\end{pmatrix}\\
& = & \varphi_{\mathcal G}\left(\frac X
t\right)(I_h+p^M\varphi_{\mathcal G}(v_1))(I_h-p^Mv)pt\mathcal
V_Y^{-1}\mathcal{E}^{-1}\begin{pmatrix} \frac 1 p I_d & 0\\ 0 &
I_{h-d}
\end{pmatrix}\\
& = & \varphi_{\mathcal G}(X)\left(\mathcal
V_Y^{-1}\mathcal{E}^{-1}+ p^M\tilde v \right)\begin{pmatrix} \frac 1
p I_d & 0\\ 0 & I_{h-d}
\end{pmatrix}
\end{eqnarray*}
where
$\tilde v= \left(\varphi_{\mathcal G}(v_1)- v - p^M\varphi_{\mathcal G}(v_1)v\right)\mathcal V_Y^{-1}\mathcal{E}^{-1}$.\\
Let us clarify these computations:
$$v, v_1 \in \frac 1 p\maxim_{[1/(p-1),p]},$$
so that
$$\varphi_{\mathcal G}(v_1) \in \frac 1 p\maxim_{[p/(p-1),p]}.$$
Thus
$$p^Mv\varphi_{\mathcal G}(v_1) \in \frac 1 p\maxim_{[p/(p-1),p]}$$
and finally,
$$\varphi_{\mathcal G}(v_1)- v - p^M\varphi_{\mathcal G}(v_1)v \in \frac 1 p\maxim_{[p/(p-1),p]}.$$
Hence, since $\mathcal V_Y^{-1} \in \frac 1 {Y^{ep/(p-1)}}\mathcal
G_{[1,p]}$,
$$p\tilde v \in \frac 1 {Y^{ep/(p-1)}} \maxim_{[ p/(p-1),p]}.$$
Finally since
$$\varphi_{\mathcal G}(X) \in pX\mathcal G_{[0,p]}$$
and
$$ X \in Y^{ep/(p-1)} \mathcal G_{[p/(p-1),\infty[}$$
we deduce the result. \Epartproof\\
Remark that
$$\varphi\left(XY \circ \frac d {dY}\right) = \frac{\varphi(X)} p Y\frac d {dY} \circ \varphi$$
and
$$ u \frac d {dY} \begin{pmatrix} l_{\mathcal A}(\beta) \\ m_{\mathcal A}(\beta)\end{pmatrix} \in  \maxim_{[ p/(p-1),p]}$$
so that we compute modulo $p^M\maxim_{[ p/(p-1),p]}$:
\begin{eqnarray*}  \varphi_{\mathcal G}  \left(XY\mathcal V_Y^{-1}\frac d {dY}
\begin{pmatrix} l_{\mathcal A}(\beta) \\ m_{\mathcal A}(\beta)\end{pmatrix}\right) & \equiv &
\frac{\varphi(X)} p Y\mathcal V_Y^{-1}\frac d {dY}\mathcal{E}^{-1}
\begin{pmatrix} \frac 1 p I_d & 0 \\ 0 & I_{h-d} \end{pmatrix}
\varphi  \begin{pmatrix} l_{\mathcal A}(\beta) \\ m_{\mathcal A}(\beta)\end{pmatrix}\\
& \equiv & \frac{\varphi(X)} p Y\mathcal V_Y^{-1}\frac d
{dY}\begin{pmatrix} \frac {\mathcal A} p \circ l_{\mathcal A}(\beta)
\\ m_{\mathcal A}(\beta) \end{pmatrix}.
\end{eqnarray*}
This yields to
\begin{eqnarray*}
& & \varphi \left(X Y \mathcal V_Y^{(-1)}\frac d {dY} \begin{pmatrix} l_{\mathcal A}(\beta) \\ m_{\mathcal A}(\beta)\end{pmatrix} \right)\\
& = & \varphi_{\mathcal G} \left(X Y \mathcal V_Y^{-1}\frac d {dY}
\begin{pmatrix} l_{\mathcal A}(\beta) \\ m_{\mathcal
A}(\beta)\end{pmatrix} \right) + \varphi_{\mathcal G} \left(X Y
\left( \mathcal V_Y^{(-1)}-\mathcal V_Y^{-1}\right)\frac d {dY}
\begin{pmatrix} l_{\mathcal A}(\beta) \\ m_{\mathcal A}(\beta)\end{pmatrix} \right)\\
&=& \textstyle \frac{\varphi(X)} p Y\mathcal V_Y^{-1} \frac d
{dY}\begin{pmatrix} \frac {\mathcal A} p \circ l_{\mathcal A}(\beta)
\\ m_{\mathcal A}(\beta) \end{pmatrix}
+ p^{M}u + \varphi_{\mathcal G} \left(X Y \left( \mathcal V_Y^{(-1)}-\mathcal V_Y^{-1}\right)\frac d {dY} \begin{pmatrix} l_{\mathcal A}(\beta) \\
 m_{\mathcal A}(\beta)\end{pmatrix} \right)
\end{eqnarray*}
with $u \in \maxim_{[p/(p-1),p]}$. Write $u = u_1 +u_2$ with $u_1
\in \frac {X^{p-1}} p \mathcal G_{[0,p]}$, thus $p^Mu_1 \in X
\maxim_{[0,p]}$ and $u_2 \in \maxim_{[ p/(p-1),\infty[}$. In
addition, $\mathcal V_Y^{(-1)}-\mathcal V_Y^{-1} \in \mathcal
G_{[0,p]}\otimes \Q_p$ hence
$$\varphi_{\mathcal G} \left(X Y \left( \mathcal V_Y^{(-1)}-\mathcal V_Y^{-1}\right)\frac d {dY} \begin{pmatrix} l_{\mathcal A}(\beta)
\\ m_{\mathcal A}(\beta)\end{pmatrix} \right) \in X \mathcal G_{[0,p]}\otimes \Q_p.$$
Write moreover
\begin{eqnarray*}
\frac{\varphi(X)} p Y\mathcal V_Y^{-1}\frac d {dY}\begin{pmatrix}
\frac {\mathcal A} p\circ l_{\mathcal A}(\beta) \\ m_{\mathcal
A}(\beta) \end{pmatrix}
& = & X Y\mathcal V_Y^{(-1)}\frac d {dY}\begin{pmatrix} \frac {\mathcal A} p\circ l_{\mathcal A}(\beta) \\ m_{\mathcal A}(\beta) \end{pmatrix}\\
& + & \left(\frac{\varphi(X)} p -X \right) Y\mathcal V_Y^{-1}\frac d {dY}\begin{pmatrix} \frac {\mathcal A} p \circ l_{\mathcal A}(\beta) \\ m_{\mathcal A}(\beta) \end{pmatrix}\\
& + & X Y\left(\mathcal V_Y^{-1} - \mathcal V_Y^{(-1)}\right) \frac
d {dY}\begin{pmatrix} \frac {\mathcal A} p \circ l_{\mathcal
A}(\beta) \\ m_{\mathcal A}(\beta) \end{pmatrix}
\end{eqnarray*}
where
$$\displaystyle{X Y\left(\mathcal V_Y^{-1} - \mathcal
V_Y^{(-1)}\right) \frac d {dY}\begin{pmatrix} \frac {\mathcal A} p
\circ l_{\mathcal A}(\beta) \\ m_{\mathcal A}(\beta)
\end{pmatrix}} \in \left(X\maxim_{[0,p]} \otimes \Q_p\right)^h$$
and
$$\displaystyle{\left(\frac{\varphi(X)} p -X \right) Y\mathcal
V_Y^{-1}\frac d {dY}\begin{pmatrix} \frac {\mathcal A} p \circ
l_{\mathcal A}(\beta) \\ m_{\mathcal A}(\beta) \end{pmatrix}} \in
X\left(\maxim_{[0,p]} + p^{M-1}\mathcal G_{[
p/(p-1),\infty[}\right)^h$$ and it can then be written as $M_1 +
M_2$ with $M_1$ in $X\maxim_{[0,p]}$ and $M_2$ in $p^M
\mathbf{\tilde A}$.

Eventually,
$$ \varphi \left(X Y \mathcal V_Y^{(-1)}\frac d {dY} \begin{pmatrix} l_{\mathcal A}(\beta) \\ m_{\mathcal A}(\beta)\end{pmatrix} \right) -
X Y\mathcal V_Y^{(-1)}\frac d {dY}\begin{pmatrix} \frac {\mathcal A}
p \circ l_{\mathcal A}(\beta) \\ m_{\mathcal A}(\beta) \end{pmatrix}
- p^M M_0$$ belongs to $M_h(X\mathcal G_{[0,p]}\otimes \Q_p)$ for
some $M_0 \in \mathbf{\tilde A}$. Then, since $X\mathcal
G_{[0,p]}\otimes \Q_p \cap\mathbf{\tilde A} = X \mathbf{\tilde
A}^+$, we deduce that
$$\varphi \left(X Y \mathcal V_Y^{(-1)}\frac d {dY} \begin{pmatrix} l_{\mathcal A}(\beta) \\ m_{\mathcal A}(\beta)\end{pmatrix} \right) =
X Y\mathcal V_Y^{(-1)}\frac d {dY}\begin{pmatrix} \frac {\mathcal A}
p\circ l_{\mathcal A}(\beta) \\ m_{\mathcal A}(\beta) \end{pmatrix}
\mod X \mathbf{\tilde A}^+ + p^M \mathbf{\tilde A}$$ which lets us
prove the proposition with the computation modulo $X \mathbf{\tilde
A}^+ + p^M \mathbf{\tilde A}$
\begin{eqnarray*}
(\varphi-1) \left(X Y \mathcal V_Y^{(-1)}\frac d {dY}
\begin{pmatrix} l_{\mathcal A}(\beta) \\ m_{\mathcal
A}(\beta)\end{pmatrix} \right)
& \equiv & X Y\mathcal V_Y^{(-1)}\frac d {dY}\begin{pmatrix} \left(\frac {\mathcal A} p -1\right)\circ l_{\mathcal A}(\beta) \\ 0 \end{pmatrix} \\
& \equiv & (\tau - 1) x
\end{eqnarray*}
and the fact that the equation $(\varphi-1)Z = \alpha \in X
\mathbf{\tilde A}^+ + p^M \mathbf{\tilde A}$ admits a solution $Z
\in X \mathbf{\tilde A}^+ + p^M \mathbf{\tilde A}$. \Eproof

\subsection{The explicit formula}
\subsubsection{Statement of the theorem} We come now to the proof of
the main theorem, the explicit formula for the Hilbert symbol of a
formal group.

\begin{thm}\label{main}
Let $\beta \in F(YW[[Y]])$ and $\alpha \in (W[[Y]][\frac 1
Y])^\times$. Denote
$$ \fcnl(\alpha) = \left(1-\frac \varphi p \right)\log \alpha(Y) = \frac 1 p \log \frac {\alpha(Y)^p} {\alpha^\varphi(Y^p)} \in W[[Y]].$$
Then the Hilbert symbol $(\alpha(\pi),\beta(\pi))_{F,M}$ has
coefficients in the basis $(o^1_M, \dots, o^h_M)$:
$$(\Tr_{W/\Z_p} \circ \Res_Y)\mathcal V_Y^{-1}\left(\begin{pmatrix} \left(1 - \frac {\mathcal A} p\right)\circ
l_{\mathcal A}(\beta) \\ 0
\end{pmatrix}d_{\log}\alpha(Y)-\fcnl(\alpha) \frac d {dY}
\begin{pmatrix} \frac {\mathcal A} p \circ l_{\mathcal A}(\beta) \\ m_{\mathcal A}(\beta)\end{pmatrix}\right).$$
\end{thm}

\subsubsection{Proof of Theorem \ref{main}}

We use the fact that if $\eta \in H^1(K,\Z/p^M\Z)$ and $r(x) \in
G_K^{ab}$ is the image by the reciprocity isomorphism $x \in K$ then
$$ \inv_K(\partial x \cup \eta) = \eta(r(x)).$$
From Proposition \ref{hilbert}, $\partial \alpha(\pi)$ corresponds
to a triple $(x,y,z)$ congruent modulo $XYW[[X,Y]]$ to
$$\left(-\frac{s(Y)} X  -\frac{s(Y)} 2 ,\ 0 ,
\ Yd_{\log} S(Y)\right)\otimes \varepsilon.$$ We compute its
cup-product with the image
$$ \left(x', 0, z'\right)$$
in $H^1(K, F[p^M])$ of $\theta(\beta)$ given by Proposition
\ref{hilbert_formel} where we recall that
$$x'=\mathcal V_Y^{(-1)}\begin{pmatrix} (\frac {\mathcal A} p -1)\circ l_{\mathcal
A}(\beta) \\ 0 \end{pmatrix}$$ and
$$z' \equiv XY\mathcal V_Y^{(-1)}\displaystyle{\frac d {dY}}\begin{pmatrix} l_{\mathcal
A}(\beta) \\ m_{\mathcal A}(\beta) \end{pmatrix} \mod X
W(\maximE).$$  We get the triple $(a,b,c)$ where:
$$a = y \mathcal V_Y^{-1}\begin{pmatrix} (\frac {\mathcal A} p -1) \circ l_{\mathcal A}(\beta) \\ 0 \end{pmatrix}
\in W(\maximE)$$ because Proposition \ref{hilbert} says that $y \in
XYW[[X,Y]]$ and Lemma \ref{frob_inv} that $XY\mathcal V_Y^{(-1)}$
has coefficients in $W(\maximE)+p^M\mathbf{\tilde A}$. Moreover,
$$c =- y \otimes \gamma z' + \sum_{n \geq 1} \binom {\chi(\gamma)} n
\sum_{k=1}^{n-1} C_{n-1}^k (\tau -1)^{k-1}z \otimes \tau^k (\tau
-1)^{n-1-k}z' \in W(\maximE)$$ because $y, z, z' \in W(\maximE)$.
Finally,
$$b = z \otimes \tau x' - x \otimes \varphi z'$$
and
$$z \otimes \tau x' = (\tau - 1) (\log (S(Y))/t +\mu)\ \tau
\left(\mathcal V_Y^{(-1)}\begin{pmatrix} (\frac {\mathcal A} p
-1)\circ l_{\mathcal A}(\beta) \\ 0 \end{pmatrix}\right) \otimes
\varepsilon.$$ On the one hand
$$(\tau - 1) (\log (S(Y))/t +\mu)
\equiv Yd_{\log} F(Y) \mod XYW[[X,Y]]$$ and on the other hand, Lemma
\ref{lem_mat_periodes} says $\tau \left(\mathcal
V_Y^{(-1)}\begin{pmatrix} (\frac {\mathcal A} p -1)\circ l_{\mathcal
A}(\beta) \\ 0 \end{pmatrix}\right)$ is congruent modulo
$XYW[[X,Y]]$ to
$$\mathcal V_Y^{(-1)}\begin{pmatrix} (\frac {\mathcal A} p -1)\circ l_{\mathcal
A}(\beta) \\ 0 \end{pmatrix} + XY\mathcal V_Y^{(-1)}\frac d
{dY}\begin{pmatrix} (\frac {\mathcal A} p -1)\circ l_{\mathcal
A}(\beta) \\ 0 \end{pmatrix}.$$ Thus, since $XY\mathcal V_Y^{(-1)}$
has coefficients in $W(\maximE)+p^M\mathbf{\tilde A}$,
$$z \otimes \tau x' \equiv Y\mathcal V_Y^{(-1)}\begin{pmatrix} (\frac {\mathcal A} p -1)\circ
l_{\mathcal A}(\beta) \\ 0 \end{pmatrix} d_{\log} S(Y) \mod
W(\maximE).$$ Finally,
$$- x \otimes \varphi z' = \left(-\frac{s(Y)} X
-\frac{s(Y)} 2\right) z' \otimes \varepsilon$$ and since
$$z' \equiv XY\mathcal V_Y^{(-1)}\displaystyle{\frac d {dY}}\begin{pmatrix}
l_{\mathcal A}(\beta) \\ m_{\mathcal A}(\beta) \end{pmatrix} \mod X
W(\maximE),$$ we get the congruence
$$- x \otimes \varphi z' \equiv Y s(Y)\mathcal V_Y^{(-1)}\displaystyle{\frac d
{dY}}\begin{pmatrix} l_{\mathcal A}(\beta) \\ m_{\mathcal A}(\beta)
\end{pmatrix} \mod W(\maximE).$$
The triple $(a,b,c)$ is eventually congruent modulo $W(\maximE)$ to
$$\left(0,Y\mathcal V_Y^{(-1)}\left(\begin{pmatrix} (\frac {\mathcal A} p -1)\circ
l_{\mathcal A}(\beta) \\ 0 \end{pmatrix} d_{\log} S(Y) +
\displaystyle{\frac d {dY}}\begin{pmatrix} l_{\mathcal A}(\beta) \\
m_{\mathcal A}(\beta)
\end{pmatrix}\frac 1 p \log \frac {S(Y)^p}{S(Y^p)}\right) \otimes \varepsilon, 0\right).$$
The theorem follows then from the lemma:

\begin{lem}
Let $C =C_{\varphi, \gamma, \tau}(\mathbf{\tilde A}_L(1))$ be the
complex computing Galois cohomology of $\Z_p(1)$.
\begin{enumerate}
\item Let $f(Y) = \sum_{n>0}{\frac {a_n}{Y^n}} \in M_h(\mathbf{\tilde A})$ be the principal part of a
series $\mathcal V_Y^{(-1)}g(Y)$ with $g(Y)$ having coefficients in
$W[[Y]]$. Then there exists a triple $(x_1,x_2,0)$ with coefficients
in $\in W(\maximE)$ such that
$$(x_1,x_2 + f(Y)\otimes \varepsilon,0) \in B^2(C).$$
In other words its image in $H^2(K,\Z_p(1))$ is zero.
\item Let $(x,y,z) \in Z^2(C)$ with $x,y,z \in W(\maximE)(1)$ then
$$(x,y,z) \in B^2(C).$$
\item Let $w\in W$ then
$$(0,w\otimes \varepsilon,0) \in Z^2(C)$$
and its image through the reciprocity isomorphism is
$\Tr_{W/\Z_p}(w)$.
\end{enumerate}
\end{lem}
\emph{Proof of the Lemma:} Put
$$w_n = \frac 1 {Y^n\left((1+X)^{-n} -1 \right)} + \frac 1 {2Y^n} \in \mathbf{\tilde A}_L.$$
Then
\begin{equation}\label{relw_n1}
(\tau-1) w_n = \frac 1 {Y^n} + \frac 1 2 (\tau-1)\frac 1 {Y^n}
\end{equation}
and
\begin{eqnarray*}
\gamma \left(\frac \varepsilon {Y^n\left((1+X)^{-n} -1 \right)}\right) & = & \frac {\chi(\gamma)\varepsilon} {Y^n\left((1+X)^{-\chi(\gamma)n} -1 \right)}\\
& = & \chi(\gamma)\delta^{-1} \left(\frac \varepsilon
{Y^n\left((1+X)^{-n} -1 \right)}\right).
\end{eqnarray*}
The Taylor expansion
$$ \delta^{-1} = \chi(\gamma) - \frac {\chi(\gamma)(\chi(\gamma)-1)} 2 (\tau-1) + (\tau-1)^2g(\tau-1)$$
where $g(\tau-1)$ is a power series in $\tau-1$ yields to the
relation
\begin{equation}\label{relw_n2}
(\gamma-1) w_n\otimes \varepsilon = g(\tau - 1) (\tau - 1) \frac 1
{Y^n}.
\end{equation}
From Lemma \ref{lem_mat_periodes}, we know $(\tau-1) \mathcal
V^{(-1)}U$ for $U\in W[[Y]]$ has coefficients in $W(\maximE)$.

Relation \eqref{relw_n1} then shows that
$$ (\tau - 1) \sum_{n>0} a_nw_n = f(Y) \mod W(\maximE)$$
and Relation \eqref{relw_n2} that
$$ (\gamma - 1) \sum_{n>0} a_nw_n = 0 \mod W(\maximE)$$
which proves that the coboundary image of triple $(\sum_{n>0}
a_nw_n,0,0)$ in $H^2(C)$ has the desired form, hence $1.$

To show $2.$ we have to solve for $x,y,z \in W(\maximE)(1)$ the
system
\begin{eqnarray*}
x & = & (\gamma -1)u + (1-\varphi)v\\
y & = & (\tau -1)u + (1-\varphi)w\\
z & = & (\tau^{\chi(\gamma)} -1)v + (\delta-\gamma)w
\end{eqnarray*}
Consider therefore $v,\ w \in W(\maximE)(1)$ solutions of
\begin{eqnarray*}
x & = & (\varphi-1)v\\
y & = & (\varphi-1)w
\end{eqnarray*}
which exist, and are unique since $\varphi-1$ is bijective on
$W(\maximE)(1)$. Then, by combining these equations with the ones of
the system, we get
$$(\varphi-1)((\tau^{\chi(\gamma)} -1)v + (\delta-\gamma)w) =-(\tau^{\chi(\gamma)} -1)x - (\delta-\gamma)y
=(\varphi-1)z.$$ But since $z$ and $(\tau^{\chi(\gamma)} -1)v +
(\delta-\gamma)w$ are elements of $W(\maximE)(1)$ where $(\varphi
-1)$ is injective, the equality
$$z =(\tau^{\chi(\gamma)} -1)v + (\delta-\gamma)w$$
holds ; $(x,y,z)$ is then a coboundary, the image of $(0, v, w)$.\\
Finally, for Point $3.$, remark that
$$(0,w\otimes\varepsilon,0) = (0,0,1\otimes \varepsilon) \cup (w,0,0).$$
Proposition \ref{hilbert} says $(0,0,1\otimes \varepsilon)$ is the
image through the Kummer map of $\pi$ a uniformizer of $K$. (To see
this, take $F(Y)=Y$.) In addition $(0,w\otimes \varepsilon,0)$
corresponds  from Theorem \ref{coho_gal} to the character $\eta$ of
$G_K$ defined in the following way: choose $b\in \mathbf{\tilde A}$
such that $(\varphi-1)b = w$, then
$$\forall g \in G, \ \eta(g) = (1-g)b.$$
Remark that since $w \in W$, we can choose $b \in W^{nr}$ and that
the image through the Kummer map of a uniformizer is the Frobenius
$\Frob_K$, thus the image through reciprocity isomorphism of
$(0,w\otimes \varepsilon,0)$ is
$$(1- \Frob_K)b = (1-\varphi^{f_K})b = (1+\varphi + \dots + \varphi^{f_K-1})w=\Tr_{W/\Z_p}w$$
where $f_K=f(K/\Q_p)$, which proves the lemma. \Epartproof

We prove then the theorem by remarking, from the congruence shown
above, that the triple $(a,b,c)$ can be written as a sum of a triple
$(0,g(Y),0)$ where $g$ is the negative part of a vector series in
$Y$ and then is zero in $H^2(K,\Z/p^M\Z)$, of a triple with
coefficients in $W(\maximE)(1)$, then also a coboundary because of
the lemma above and finally a triple $(0,w\otimes \varepsilon,0)$
where $w$ is the constant term of the vector series
$$Y\mathcal V_Y^{(-1)}\left(\begin{pmatrix} (\frac {\mathcal A} p -1)\circ
l_{\mathcal A}(\beta) \\ 0 \end{pmatrix} d_{\log} \alpha(Y) +
\displaystyle{\frac d {dY}}\begin{pmatrix} l_{\mathcal A}(\beta) \\
m_{\mathcal A}(\beta) \end{pmatrix}\frac 1 p \log \frac
{\alpha(Y)^p}{\alpha(Y^p)}\right) $$ hence the residue of
$$\mathcal V_Y^{-1}\left(\begin{pmatrix}
\left(\frac {\mathcal A} p -1\right)\circ l_{\mathcal A}(\beta) \\ 0
\end{pmatrix} d_{\log} \alpha(Y) + \displaystyle{\frac d
{dY}}\begin{pmatrix} l_{\mathcal A}(\beta) \\ m_{\mathcal A}(\beta)
\end{pmatrix}\frac 1 p \log \frac {\alpha(Y)^p}{\alpha(Y^p)}\right).$$
The only term with a non zero contribution is then the residue, and
that contribution is, according to the lemma, given by the trace,
which completes the proof. \Eproof

\addcontentsline{toc}{section}{References}
\bibliographystyle{amsalpha}

\def\cprime{$'$}
\providecommand{\bysame}{\leavevmode\hbox
to3em{\hrulefill}\thinspace}
\providecommand{\MR}{\relax\ifhmode\unskip\space\fi MR }
\providecommand{\MRhref}[2]{%
  \href{http://www.ams.org/mathscinet-getitem?mr=#1}{#2}
} \providecommand{\href}[2]{#2}

\bigskip

Floric Tavares Ribeiro\\
D\'epartement de Math\'ematiques, Universit\'e de Franche-Comt\'e\\
16 route de Gray, 25000 Besan\c con, France\\
\emph{email:} ftavares@univ-fcomte.fr

\end{document}